\documentclass{article}
\usepackage{amsfonts}
\usepackage{graphicx}
\usepackage[dvips]{color}
\usepackage{amssymb}
\usepackage{url}

\usepackage[section]{placeins}	
\usepackage{floatrow}
\usepackage{caption}
\title{Brouwer and Euclid}        
\author{Michael Beeson}        
\date{\today
}          
\def\axioms{\begin{tabular*}{0.8\textwidth}{@{\extracolsep{\fill}}lr}}
\def\endaxioms{\end{tabular*}\smallskip}
\def\satisfies{{\models}}   

\def\implies{\rightarrow}

\def\seq{\Rightarrow}

\def\R{{\mathbb R}}

\def\F{{\mathbb F}}

\def\M{{\textsf M}}

\def\B{{\bf B}}
\def\T{{\bf T}}   

\medskip
\newtheorem{theorem}{Theorem}
\newtheorem{definition}{Definition}
\newtheorem{lemma}{Lemma}
\newtheorem{corollary}{Corollary}
\def\newaxioms#1{\begin{tabular*}{0.8\textwidth}{@{}l@{\hskip-#1cm}r}}
\def\endaxioms{\end{tabular*}\smallskip}

\usepackage{graphicx}
\usepackage{pstricks}
\psset{unit=3cm}

\def\OneSixteenInnerFigure{%
\pspicture(0,-0.15)(2,0.75)
\pspolygon[fillstyle=solid,fillcolor=yellow](1.8,0)(0,0)(2,0.8)(1.4,0.2)(1.8,0)
\psline[linecolor=red,linestyle=dashed](1,0.4)(1.8,0)  
\psline(0,0)(1.8,0)  
\qline(0,0)(0.8,0.8)  
\qline(0.8,0.8)(1.2,0)  
\psdot(0.8,0.8)   
\put(0.65,0.77){$A$}
\psdot(1.2,0)  
\put(1.17,-0.13){$C$}
\psdot(1.8,0)  
\put(1.75,-0.13){$D$}
\psline[linecolor=blue](1.2,0)(2 ,0.8) 
\psline[linecolor=blue](0,0)(1.19,0.477)  
\psline[linecolor=blue](1.21,0.483)(2,0.8)  
\psdot(1,0.4)  
\put(0.83,0.4){$E$}
\psdot(2,0.8)  
\put(2,0.67){$F$}
\pscircle(1.4,0.2){0.03} 
\put(1.5,0.19){$H$}
\psdot(0,0)  
\put(-0.04,-0.13){$B$}
\endpspicture}

\def\InnerOuterPaschFigure{%
\hskip0.9cm   
\pspicture(3.6,1.8)
\psdot(0,0.6)
\put(0,0.47){$a$}
\pscircle(1,0.6){0.03}
\put(1,0.47){$x$}
\qline(0,0.6)(0.97,0.6)
\psdot(1.3,0.6)
\put(1.35,0.55){$q$}
\qline(1.03,0.6)(1.3,0.6)
\psdot(1.1,1.4)
\put(1.1,1.45){$c$}
\qline(0,0.6)(1.1,1.4)  
\qline(1.1,1.4)(1.375,0.3) 
\psdot(1.375,0.3)
\put(1.43,0.25){$b$}
\psdot(0.52,0.98)
\put(0.43,1.03){$p$}
\qline(1.375,0.3)(1.015,0.584)  
\qline(0.52,0.98)(0.979,0.62)  
\psdot(2,0.3)
\put(1.95,0.19){$a$}
\psdot(3.5,0.3)
\put(3.45,0.19){$q$}
\qline(2,0.3)(2.47,0.3)
\qline(2.53,0.3)(3.5,0.3)
\pscircle(2.5,0.3){0.03}
\put(2.45,0.19){$x$}
\psdot(2.7,1.5)
\put(2.65,1.57){$b$}
\qline(2.7,1.5)(3.5,0.3)  
\psdot(3.0,1.05)
\put(3.07,1.02) {$c$}
\qline(2.7,1.5)(2.5,0.33)  
\qline(2,0.3)(3.0,1.05)   
\psdot(2.566,0.722)
\put(2.6,0.6){$p$}
\endpspicture}

\def\PositiveHypotenuseFigure{%
\pspicture(0,-0.4)(2,1.3)
\pspolygon[fillstyle=solid,fillcolor=yellow](0,0)(2,0)(0.8,1.2)(0.6,0.4)(0,0)
\psdot(0,0)  
\put(-0.04,-0.11){$a$}
\psdot(1,0)  
\put(0.97,-0.11){$b$}
\psdot(2,0)  
\put(1.94,-0.11){$d$}
\qline(0,0)(2,0)  
\psdot(0.5,0) 
\put(0.49,-0.11){$m$}
\psdot(0.8,1.2) 
\put(0.83,1.2){$u$}
\psdot(0.7,0.8) 
\put(0.64,0.85){$\ell$}
\psdot(1.6,0.4) 
\put(1.65,0.4){$e$}
\psdot(1.2,0.8) 
\put(1.25,0.8){$f$}
\psdot(0.2,0.8) 
\put(0.18,0.85){$v$}
\psdot(-0.4,0.4) 
\put(-0.5,0.4){$p$}
\psdot(0.4,-0.4)  
\put(0.43,-0.47){$q$}
\psline(-0.4,0.4)(0.8,1.2) 
\psline(-0.4,0.4)(0.4,-0.4) 
\psline(0.4,-0.4)(1.6,0.4)  
\psline(0.2,0.8)(1.2,0.8)  
\psline(0.57,0.4)(-0.4,0.4)  
\psline(0.575,0.425)(0.2,0.8)  
\psline(0.625,0.375)(1,0)  
\psline(0.63,0.4)(1.6,0.4)  
\psline(0.625,0.425)(1.2,0.8)  
\pscircle[fillcolor=white](0.6,0.4){0.03}  
\put(0.44,0.43){$c$}
\psline(0.59,0.37)(0.4,-0.4)  
\endpspicture}

\def\LowerDimensionFigure{%
\pspicture(1,1)(0,-0.07)
\psdot(0,0)
\put(0,-0.1){$\beta$}
\psdot(1,0)
\put(1,-0.1){$\gamma$}
\psdot(0.5,0.866)
\put(0.4,0.87){$\alpha$} 
\qline(0,0)(1,0)
\qline(1,0)(0.5,0.866)
\qline(0.5,0.866)(0,0)
\psdot(0.75,0.433)
\put(0.83,0.39){$c_2$}
\put(0.47,0.18){$c_4$}
\psline[linestyle=dashed](0,0)(0.75,0.433)
\psdot(0.25,0.433)
\put(0.1,0.43){$c_1$}
\psline[linestyle=dashed](1,0)(0.25,0.433)
\psdot(0.5,0.288)
\psdot(0.5,0)
\put(0.47,-0.1){$c_3$}
\endpspicture}

\def\SixtyDegreesFigure{%
\pspicture(1,0)(3,1)
\psdot(1,0)
\put(1,-0.1){$a$}
\psdot(2,0)
\put(2,-0.1){$c$}
\psarc[linecolor=green]{->}(2,0){0.1}{60}{180}
\psdot(3,0)
\put(3,-0.1){$e$}
\psline[linecolor=blue](1,0)(2,0)
\qline(2,0)(3,0)
\psdot(1.5,0.866)
\put(1.4,0.87){$f$}
\psdot(2.5,0.866)
\put(2.55,0.87){$g$}
\psline[linecolor=yellow](1.5,0.866)(2.5,0.866)  
\psline[linecolor=yellow](1.5,0.866)(3,0)   
\psdot(2.25,0.433)  
\put(2.3,0.43){$m$}
\qline(1,0)(1.5,0.866)
\qline(1.5,0.866)(2,0)
\psline[linecolor=blue](2,0)(2.5,0.866)
\qline(2.5,0.866)(3,0)
\psdot(1.75,0.433)
\put(1.83,0.39){$x$}
\psline[linestyle=dashed,linecolor=red](1,0)(2.5,0.866)
\endpspicture}
 
\def\OneFiftyFigure{%
\psset{unit=1cm}
\pspicture(-2,-1)(5.8,3)
\pspolygon[linecolor=lightgray](-1.732,1)(-1.732,-1)(0,0)
\pspolygon[linecolor=lightgray](-1.732,0)(-0.866,0.5)(-0.866,-0.5)
\pspolygon[linecolor=lightgray](0,0)(1,1.732)(2,0)
\pspolygon[linecolor=lightgray](1,1.732)(2,0)(3,1.732)
\pspolygon[linecolor=lightgray](3,1.732)(4.732,0.732)(4.732,2.732)
\pspolygon[linecolor=lightgray](0.5,0.866)(1.5,0.866)(1,0)
\pspolygon[linecolor=lightgray](1.5,0.866)(2,1.732)(2.5,0.866)
\pspolygon[linecolor=lightgray](3.866,2.232)(3.866,1.232)(4.732,1.732)
\psdot(0,0)
\put(-0.1,-0.35){$b$}
\psdot(-1.732,0)
\put(-2,-0.2){$a$}
\psdot(1.5,0.866)
\put(1.3,0.6){$c$}
\psdot(4.732,1.732)
\put(4.85,1.7){$e$}
\psdot(2,0)
\put(2.0,-0.35){$d$}
\qline(0,0)(2,0)
\psline[linecolor=red,linestyle=dashed](-1.732,0)(4.598,1.732) 
\psline[linecolor=blue](-1.732,0)(0,0)  
\psdot(4.732,2.732)
\put(4.85,2.7){$f$}
\psline[linecolor=blue](0,0)(4.732,2.732)  
\psarc[linecolor=green]{->}(0,0){0.25}{30}{180}
\psset{unit=3cm}
\endpspicture}

\def\SaccheriHelperFigure{%
\psset{unit=2.7cm}
\pspicture(0.4,-0.5)(1.9,1.5)
\qline(0.5,1)(2,1)
\qline(0.5,1.5)(2,0)
\qline(2,0)(2,1)
\qline(0.5,1)(0.5,1.5)
\psdot(0.5,1)
\put(0.39,0.95){$a$}
\psdot(0.5,1.5)
\put(0.39,1.45){$b$}
\psdot(2,1)
\put(2.04,0.95){$d$}
\psdot(2,0)
\put(2.04,-0.05){$c$}
\psdot(1,1)
\put(1.0,1.05){$m$}
\endpspicture
\pspicture(0,-0.5)(2,1.5)
\pspolygon[fillstyle=solid,fillcolor=yellow](0.5,1.5)(1,1)(2,1)(0.5,-0.3)
\qline(0.5,1)(2,1)
\qline(0.5,1.5)(2,0)
\qline(2,0)(2,1)
\psline[linecolor=gray](0.5,-0.3)(2,1)
\qline(0.5,1)(0.5,1.5)
\psline[linecolor=gray](0.5,-0.3)(0.5,1)
\psdot(0.5,1)
\put(0.39,0.95){$a$}
\psdot(0.5,1.5)
\put(0.39,1.45){$b$}
\psdot(2,1)
\put(2.04,0.95){$d$}
\psdot(2,0)
\put(2.04,-0.05){$c$}
\psdot(1,1)
\put(1.0,1.05){$m$}
\psdot(1.462,0.54)
\put(1.54,0.52){$p$}
\psdot[linecolor=gray](0.5,-0.3)
\put(0.39,-0.31){$x$}
\psset{unit=3cm}
\endpspicture}

\def\PositiveApexFigure{%
\pspicture(0,0)(2,0.8)
\psline(0,0)(1.45,0)  
\qline(0,0)(1.2,0.8)  
\qline(1.2,0.8)(1.2,0)  
\psdot(1.2,0.8)   
\put(1.16,0.85){$a$}
\psdot(1.2,0)  
\put(1.17,-0.09){$c$}
\psdot(1.45,0)  
\put(1.45,-0.11){$d$}
\psdot(0,0)  
\put(-0.04,-0.11){$b$}
\psarc(0,0){1.45}{0}{34}
\psline[linecolor=red](1.2,0.8)(1.45,0)
\endpspicture}

 \def\CrossbarFigure{%
\pspicture(0,0)(1.7,0.8)
\psline(0,0)(1.45,0)  
\qline(0,0)(1.2,0.8)  
\psdot(0.67,0.45)  
\put(0.64,0.5){$a$}
\qline(0.67,0.45)(1.2,0)  
\psdot(0.962,0.2)   
\put(0.962,0.24){$e$}
\pscircle(1.36,0.28){0.03}   
\put(1.23,0.3){$w$}
\qline(0,0)(1.33,0.275) 
\qline(1.2,0.8)(1.35,0.305)  
\qline(1.37,0.26)(1.45,0)  
\psdot(1.2,0.8)   
\put(1.16,0.85){$u$}
\psdot(1.2,0)  
\put(1.17,-0.09){$c$}
\psdot(1.45,0)  
\put(1.45,-0.11){$v$}
\psdot(0,0)  
\put(-0.04,-0.11){$b$}
\endpspicture}

 \def\CrossbarProofFigure{%
\pspolygon[fillcolor=yellow,fillstyle=solid,linestyle=none]
(0,0)(1.2,0)(1.45,0)(0.67,0.45)(0.962,0.2)
\pspicture(0,0)(1.7,0.7)
\psline(0,0)(1.45,0)  
\qline(0,0)(1.2,0.8)  
\psdot(0.67,0.45)  
\put(0.64,0.5){$a$}
\qline(0.67,0.45)(1.2,0)  
\psdot(0.962,0.2)   
\put(0.82,0.2){$e$}
\pscircle(1.07,0.22){0.03} 
\put(1.03,0.28){$f$}
\qline(0.67,0.45)(1.05,0.23)   
\qline(0,0)(1.05,0.22)   
\qline(1.09,0.21)(1.45,0)   
\qline(1.2,0.8)(1.45,0)  
\psdot(1.2,0.8)   
\put(1.16,0.85){$u$}
\psdot(1.2,0)  
\put(1.17,-0.09){$c$}
\psdot(1.45,0)  
\put(1.45,-0.11){$v$}
\psdot(0,0)  
\put(-0.04,-0.11){$b$}
\endpspicture}

\def\CrossbarProofFigureTwo{%
\pspicture(0,0)(1.5,0.7)
\pspolygon[fillcolor=yellow,fillstyle=solid,linestyle=none]
(0,0)(0.67,0.45)(1.2,0.8)(1.45,0)(1.07,0.22)
\qline(0,0)(1.45,0)  
\qline(0.67,0.45)(1.45,0)  
\qline(0,0)(1.2,0.8)  
\psdot(0.67,0.45)  
\put(0.64,0.5){$a$}
\qline(0.67,0.45)(1.2,0)  
\psdot(0.962,0.2)   
\put(0.82,0.2){$e$}
\psdot(1.07,0.22) 
\put(1.03,0.28){$f$}
\pscircle(1.36,0.28){0.03}   
\put(1.23,0.3){$w$}
\qline(0,0)(1.33,0.275) 
\qline(1.2,0.8)(1.35,0.305)  
\qline(1.37,0.26)(1.45,0)  
\psdot(1.2,0.8)   
\put(1.16,0.85){$u$}
\psdot(1.2,0)  
\put(1.17,-0.09){$c$}
\psdot(1.45,0)  
\put(1.45,-0.11){$v$}
\psdot(0,0)  
\put(-0.04,-0.11){$b$}
\endpspicture}

\def\PlayfairFigure{%
\pspicture(0,0.15)(2.2, 1.05)
\qline(0.0,0.15)(2.0,0.15)
\qline(0.0,0.7)(2.0,0.7)
\qline(0.0,0.9)(2.0,0.5)
\psdot(1,0.7)
\put(1,0.78){$p$}
\put(-0.15,0.12) {$L$}
\put(-0.15,0.68) {$K$}
\put(-0.15,0.88) {$M$}
\endpspicture}

\def\EuclidParallelFigure{%
\pspicture(2.7, 0.9)
\qline(0.0,0.15)(2.47,0.15)
\qline(2.53,0.15)(2.7,0.15)
\qline(0.0,0.7)(2.7,0.7)
\qline(0.0,0.9)(2.473,0.163)
\qline(2.5275,0.14)(2.7,0.085)
\pscircle(2.5,0.15){0.03}
\psdot(0.67,0.7)
\put(0.67,0.78){$p$}
\psdot(1.29,0.515)
\put(1.28,0.39) {$a$}
\psdot(0.9,0.15)
\put(0.88,0.04) {$q$}
\psline[linestyle=dashed](0.67,0.7)(0.9,0.15)
\psdot(1.49,0.7)
\put(1.5,0.78) {$r$}
\qline(0.9,0.15)(1.49,0.7)
\put(-0.15,0.12) {$L$}
\put(-0.15,0.68) {$K$}
\put(-0.15,0.88) {$M$}
\endpspicture}

\def\BrouwerParallelFigure{%
\pspicture(2.7, 1.2)
\qline(0.0,0.15)(2.47,0.15)
\qline(2.53,0.15)(2.7,0.15)
\qline(0.0,0.7)(2.7,0.7)
\qline(0.0,0.9)(2.473,0.163)
\qline(2.5275,0.14)(2.7,0.085)
\pscircle(2.5,0.15){0.03}
\psdot(0.67,0.7)
\put(0.67,0.78){$p$}
\psdot(0.9,0.15)
\put(0.88,0.04) {$q$}
\psline[linestyle=dashed](0.67,0.7)(0.9,0.15)
\psdot(1.49,0.7)
\put(1.5,0.78) {$r$}
\qline(0.9,0.15)(1.49,0.7)
\put(-0.15,0.12) {$L$}
\put(-0.15,0.68) {$K$}
\put(-0.15,0.88) {$M$}
\endpspicture}

\def\EuclidParallelRawFigure{%
\pspicture(3.2, 1.2)
\pspolygon[fillstyle=solid,fillcolor=yellow](0.67,0.7)(0.78,0.425)(1.49,0.7)
\pspolygon[fillstyle=solid,fillcolor=yellow](0.07,0.15)(0.78,0.425)(0.9,0.15)
\qline(-0.15,0.15)(2.47,0.15)
\qline(2.53,0.15)(2.7,0.15)
\qline(-0.15,0.7)(2.7,0.7)
\qline(0.0,0.9)(2.473,0.163)
\qline(2.5275,0.14)(2.7,0.085)
\pscircle(2.5,0.15){0.03}
\psdot(0.67,0.7)
\put(0.67,0.78){$p$}
\psdot(1.29,0.515)
\put(1.28,0.39) {$a$}
\psdot(0.9,0.15)
\put(0.88,0.04) {$q$}
\psdot(0.07,0.15)
\put(0.07,0.04){$s$}
\qline(0.67,0.7)(0.9,0.15)
\psdot(1.49,0.7)
\put(1.5,0.78) {$r$}
\qline(1.49,0.7)(0.07,0.15)  
\psdot(0.78,0.425)
\put(0.65,0.3){$t$}
\qline(0.9,0.15)(1.49,0.7)  
\put(-0.3,0.12) {$L$}
\put(-0.3,0.68) {$K$}
\put(-0.15,0.88) {$M$}
\endpspicture}

\def\TarskiFiveSegmentFigure{%
\psset{unit=2.25cm}
\pspicture(2.5,0)(2.2, 1.0)
\qline(0.0,0.15)(2.0,0.15)
\qline(0.0,0.15)(1.0,0.85)
\psline[linestyle=dashed](1.0,0.85)(2.0,0.15)
\qline(1.0,0.85)(0.7,0.15)
\put(1,0.9){$d$}
\put(0.0,0){$a$}
\put(0.7,0){$b$}
\put(2.0,0){$c$}
\qline(2.5,0.15)(4.5,0.15)
\qline(2.5,0.15)(3.5,0.85)
\psline[linestyle=dashed](3.5,0.85)(4.5,0.15)
\qline(3.5,0.85)(3.2,0.15)
\put(3.5,0.9){$D$}
\put(2.5,0){$A$}
\put(3.2,0){$B$}
\put(4.5,0){$C$}
\psset{unit=3cm}
\endpspicture}

 \def\GuptaMidpointFigureOne{%
\psset{unit=2cm}
\pspicture(0,-0.2)(2,1.4)  
\pspolygon[fillstyle=solid,fillcolor=yellow,linestyle=none](1.414,1.414)(0,0)(2,0)(1.414,0.5857)(1.414,1.414)
\psdot(0,0)
\put(-0.15,-0.17){$C$}
\psdot(1.414,0)
\put(1.4,-0.17){$B$}
\psdot(1.414,1.414)
\put(1.48,1.4){$E$}
\psdot(2,0)
\put(2,-0.17){$D$}
\psdot(1,1)
\put(0.9,1.1){$A$}
\qline(0,0)(2,0)  
\qline(0,0)(1.414,1.414) 
\qline(1.414,0)(1.414,0.545) 
\qline(1.414,0.62)(1.414,1.414)  
\qline(1,1)(1.397,0.60)  
\qline(2,0)(1.43,0.57)  
\qline(1,1)(1.414,0) 
\qline(0,0)(1.39,0.575)
\qline(1.45,0.594)(2,0.8284)
\pscircle(1.414,0.5857){0.04}
\put(1.55, 0.5){$F$}
\psset{unit=3cm}
\endpspicture}

 \def\GuptaMidpointFigureTwo{%
\psset{unit=2cm}
\pspicture(0,-0.2)(2,1.4)  
\pspolygon[fillstyle=solid,fillcolor=yellow,linestyle=none](0,0)(1.414,1.414)(1.414,0)(1.207,0.5)(0,0)
\psdot(0,0)
\put(-0.15,-0.17){$C$}
\psdot(1.414,0)  
\put(1.4,-0.17){$B$}
\psdot(1.414,1.414)
\put(1.48,1.4){$E$}
\psdot(2,0)
\put(2,-0.17){$D$}
\psdot(1,1)
\put(0.9,1.1){$A$}
\qline(0,0)(2,0)  
\qline(0,0)(1.414,1.414) 
\qline(1.414,0)(1.414,1.414) 
\qline(1,1)(1.397,0.60)  
\qline(2,0)(1.43,0.57)  
\qline(1,1)(1.20,0.53) 
\qline(1.225,0.47)(1.414,0) 
\qline(0,0)(1.187,0.491)
\qline(1.222,0.512)(1.39,0.575)
\qline(1.44,0.594)(2,0.8284)
\pscircle(1.207,0.5){0.04}
\put(1.05,0.3){$M$}
\psdot(1.414,0.5857) 
\put(1.55, 0.5){$F$}
\psset{unit=3cm}
\endpspicture}

\def\LemmaCrossbarFigure{%
\pspicture(1.8,1.5)(0,0.3)
\pspolygon[fillstyle=solid,fillcolor=yellow, linestyle=none](0,0.6)(1.1,1.4)(1.375,0.3)(1,0.6)
\psdot(0,0.6) 
\put(0,0.47){$a$}
\pscircle(1,0.6){0.03}
\put(0.94,0.46){$f$}
\qline(0,0.6)(0.97,0.6)
\psdot(1.3,0.6) 
\put(1.35,0.55){$e$}
\qline(1.03,0.6)(1.3,0.6)
\psdot(1.1,1.4)
\put(1.1,1.45){$b$}
\qline(0,0.6)(1.1,1.4)  
\qline(1.1,1.4)(1.375,0.3) 
\psdot(1.375,0.3) 
\put(1.43,0.25){$c$}
\psdot(0.52,0.98)
\put(0.43,1.03){$d$}
\qline(1.375,0.3)(1.02,0.584)  
\qline(0.52,0.98)(0.979,0.62)  
\qline(0,0.6)(1.375,0.3)  
\endpspicture}

\def\DiagonalsOfParallelogramFigure{%
\pspicture(-0.1,-0.1)(2,1.5)
\pspolygon[fillcolor=yellow,fillstyle=solid,linestyle=none]
(0,0)(1.6,1.6)(1,0)(0.65,0.4)(0,0)
\psdot(0,0)  
\put(-0.05,-0.13){$D$}
\psdot(1,0) 
\put(0.95,-0.13){$C$}
\qline(0,0)(1,0)
\qline(0,0)(0.3,0.8)
\qline(0.3,0.8)(1.3,0.8) 
\psdot(0.3,0.8)  
\put(0.15,0.73){$A$}
\psdot(1.6,1.6) 
\put(1.64,1.53){$E$}
\psdot(1.3,0.8) 
\put(1.33,0.73){$B$}
\qline(1,0)(1.6,1.6)  
\qline(0.3,0.8)(1.6,1.6) 
\qline(0,0)(1.6,1.6) 
\psdot(0.8,0.8) 
\put(0.75,0.85){$F$}
\pscircle(0.65,0.4){0.03}  
\put(0.55,0.25){$M$}
\qline(0,0)(0.633,0.387) 
\qline(0.67,0.41)(1.3,0.8) 
\qline(0.67,0.385)(1,0) 
\qline(0.633,0.41)(0.3,0.8) 
\psdot(0.53,0.53) 
\put(0.37,0.48){$H$}
\endpspicture
}

\begin{document}
\maketitle
\abstract{We explore the relationship between Brouwer's intuitionistic mathematics 
and Euclidean geometry. Brouwer wrote a paper in 1949 called 
{\em The contradictority of elementary geometry}.  In 
that paper, he showed that a certain classical 
consequence of the parallel postulate implies Markov's principle,
which he found intuitionistically unacceptable.  But 
Euclid's geometry, having served as a beacon of clear and correct
reasoning for two millenia, is not so easily discarded.

Brouwer started from a ``theorem'' that is not in 
Euclid, and requires Markov's principle for its proof.  That means
that Brouwer's paper did not  address the question whether
Euclid's {\em Elements} really requires Markov's principle. 
In this paper we show that there is a coherent theory of 
``non-Markovian Euclidean geometry.''  
 We show in some detail that our theory is an adequate formal 
rendering of (at least) Euclid's Book~I, and suffices to define geometric arithmetic,  thus refining
the author's previous investigations (which include Markov's principle as an axiom).  

Philosophically, Brouwer's proof that his version of the parallel postulate implies 
Markov's principle could be read just as well as geometric evidence for the truth of Markov's principle,
if one thinks the geometrical ``intersection theorem''
with which Brouwer started is geometrically evident.      
}

\tableofcontents

\section{Introduction} 
Brouwer, in founding the philosophy of mathematics known as ``intuitionism'',
 rejected many of the mathematical results  that were obtained in the nineteenth 
century or before.  The rejected body of mathematics has become known as 
``classical mathematics''.  The word ``classical'' also brings to mind the 
world of ancient Greece, where Euclid and his predecessors laid the foundations of 
modern mathematical reasoning in the third century BCE.  Euclidean geometry was 
for two millenia the {\em sine qua non} of careful reasoning; every educated person
in Europe studied it; the American Declaration of Independence was modeled on 
Euclidean reasoning.  Brouwer did not directly challenge Euclid 
by name in any publication, but in 1949 he published a paper \cite{brouwer1949b} with the 
 title, {\em Contradictority of Elementary Geometry}.

Brouwer was never much of a diplomat.   If his personality had been more diplomatic,
he might have pointed out that, on a certain reading, 
certain theorems of elementary geometry related to the 
parallel postulate may 
seem (intuitionistically) contradictory;  but that every classical theorem permits various refinements,
once constructive distinctions are taken into account, and Euclid's parallel postulate and its consequences are not exceptions. 

What would Euclid have written,  if 
he had come after Brouwer, instead of before?   Of course, he would not have thrown 
up his hands, thinking geometry is contradictory,  and gone into investment banking instead.  We will show in this paper that, if one is careful
about the formulation of the axioms, Euclidean geometry is
perfectly consistent with Brouwer's intuitionism.  
 
What Brouwer calls a ``contradiction'' has two parts:
(i) Brouwer rejects a
certain property of the ordering of points on a line known 
as Markov's principle, and (ii)
Brouwer shows that a certain ``intersection theorem''
(a classical consequence of Euclid's parallel postulate)
implies Markov's principle, even though it appears not to
imply the law of the excluded middle.   This result 
would be better summarized by the statement
 
\begin{quote}
{\em Intuitionistic  ``elementary geometry''
must distinguish between different propositions
classically equivalent to Euclid's parallel postulate.} 
\end{quote}

In previous work \cite{beeson2016b, beeson2015b}, such 
distinctions were made, and a coherent theory of constructive 
geometry developed.  But that theory would not have met 
with Brouwer's approval, because Markov's principle is 
assumed as an axiom.  The reason for that was simple 
pragmatism:  it enables one to 
argue by contradiction and cases,
 as long as one is trying to prove betweenness
and congruence assertions about specific points, rather than
assertions that more points exist with certain properties.%
\footnote{For logicians: allowing Markov's principle
 enabled the double negation interpretation 
to work, allowing the ``importation'' of a certain class
of geometrical results, whose classical proofs are long
and complicated.} 

Thus, the door is open (and has been open for 68 years)
for a development of Euclidean geometry that 
explicitly avoids not only the law of the excluded 
middle, but also Markov's principle.   
We call this theory ``non-Markovian Euclidean geometry'', or for short
just ``non-Markovian geometry.''  Perhaps it should 
be called ``intuitionistic geometry'',  if one feels that 
the rejection of Markov's principle is fundamental to 
intuitionism.    

The key concepts of non-Markovian geometry
are the concepts of ``distinct points'' and ``positive angles.''
The concept that $a$ and $b$ are distinct is written $a \# b$, and is 
stronger than simple inequality $a \neq b$.   Intuitively, $a \# b$ means
that we have a positive lower bound on how far apart $a$ and $b$ are,
although in axiomatic geometry there is of course no notion of ``distance.''
The concept that angle $abc$ is ``positive'',  written $0 < abc$, means
intuitively that we have a lower bound on how different the directions
of the rays $bc$ and $ba$ are.  

Both these concepts will be {\em defined}
in terms of betweenness and congruence,  rather than be introduced as primitive.
In particular the correct definition of ``positive angle'' is not obvious
{\em a priori}.  In order to ensure that the axioms do not imply Markov's
principle immediately,  Pasch's axiom must be restricted to assume that 
certain angles are positive.  Then the definition of ``positive angle'' 
must be broad enough to permit the applications of Pasch that are needed
to ``bootstrap'' geometry.  But finally, we should be able to prove that 
a positive angle is simply the apex angle of an isosceles triangle (whose
three points are distinct).  Then Euclid Book~I can be proved, 
when angles are assumed to be positive (and have positive supplements)
 and triangles, by definition
of ``triangle'', have distinct vertices.  In the last section
we give a metamathematical theorem to the effect that this claim
can be extended to (at least) Books II and III as well.%
\footnote{This claim cannot be taken too literally, as of course
Euclid does contain gaps and errors.}

There have been some papers on related subjects in the 
 seven decades since Brouwer's rejection of ``elementary geometry'',
and the obligation arises to explain in what relation the present
work stands to those papers.  Almost all of those 
papers were  about projective geometry, 
or affine geometry, or Desarguesian geometry,  rather than Euclidean
geometry, and also were based on (or included) apartness, which is rejected 
in this work.  Here is a brief, possibly 
incomplete,  list of such papers.  The first was Heyting's 1925 thesis  
(published two years later as \cite{heyting1927} and 
again 34 years later as \cite{heyting1959}.)
Heyting's student van Dalen continued work on 
intuitionistic projective spaces in
\cite{vanDalen1963,vanDalen1990,vanDalen1996}; see also \cite{vonPlato1995} 
and \cite{mandelkern2007,mandelkern2013,mandelkern2014,mandelkern2016}.  
  Aside from the previous work of
the present author \cite{beeson2015b,beeson2016c},  the 
only previous paper on constructive Euclidean geometry was by 
Lombard and Vesley \cite{lombard-vesley}, who
followed Heyting in taking apartness as primitive.

\section{Versions of the parallel postulate}
The ``parallel postulate'' of Euclid has been reformulated many times in 
the history of geometry, as efforts to prove it led instead to many 
different equivalent propositions.   But not all these versions are 
equivalent using intuitionistic logic.  
Euclid's version requires two lines to meet, if two specified angles 
``make less than two right angles.'' 
A different version
due to the Englishman Playfair (1729), 
was popularized by Hilbert; it asserts the impossibility of two 
different lines parallel to a given line through the same point.
In so doing it makes no existential assertion, unlike Euclid's 
version, which asserts the existence of an intersection point.

 In this section,
we will review several versions of the parallel postulate, and then 
discuss Brouwer's 1949 proof. 
All these versions of the 
parallel postulate, including Euclid's own,  are consistent with 
intuitionistic logic.  Brouwer's paper shows that
his version of the parallel postulate implies a 
certain ordering principle, known as Markov's principle,  which 
Brouwer  believed to be contrary to the nature of the intuitionistic
continuum, for reasons far removed from Euclidean geometry.  

\subsection{Euclid's parallel axiom}
Euclid's postulate 5 is 
\begin{quote} 
{\em If a straight line falling on two straight lines make the interior angles on the 
same side less than two right angles, the two straight lines, if produced indefinitely, meet on that side on which are 
the angles less than the two right angles.}  
\end{quote}

\begin{figure}[h]
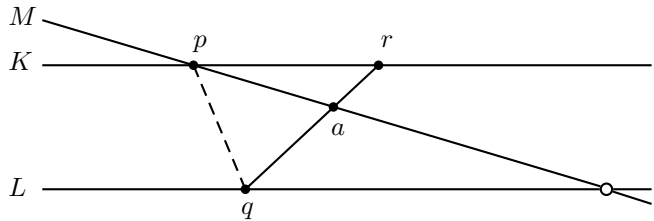
 
\center{\EuclidParallelFigure}  
\caption{Euclid~5:  $M$ and $L$ must meet on the right side, provided $\B(q,a,r)$ and $pq$ makes
alternate interior angles equal with $K$ and $L$.}
\label{figure:EuclidParallelFigure}
\end{figure}

We consider the formal expression of Euclid's parallel axiom.  Let $L$ and $M$ be two straight 
lines, and let $pq$ be the ``straight line falling on'' $L$ and $M$, with $p$ on $M$ and $q$ 
on $L$.  We think that what Euclid meant by ``makes the interior angles on the same side 
less than two right angles'' was that, if $K$ is another line through $p$,
making the interior angles with $pq$ equal to two right angles, then $M$ would lie in 
the interior of one of those interior angles (see Fig.~\ref{figure:EuclidParallelFigure}).

Euclid did not define ``angle'', and did not define ``lies in the
interior of an angle'', but these issues of precision have little to do 
with intuitionism.  Assuming for the moment that we understand the notions
of ``angle'' and ``alternate interior angle'', then we can state Euclid's
parallel axiom, using 
  three more points to ``witness'' that one ray of line $M$ emanating from $p$
lies in the interior of one of the interior angles made by $K$.   
Fig.~\ref{figure:EuclidParallelFigure}
illustrates the postulate.  The point asserted to exist is shown by a
small open circle (a convention we will follow throughout).  
 
\smallskip    

In this formulation, the point is asserted to exist ``on the right side'';
more precisely,  ``on the same side of $pq$ as $a$.''  
Euclid did not define ``on the same side of'',  although his postulate 
uses that phrase.   
That notion was, apparently, first defined by M. Pasch in 1882
(on page 27 of \cite{pasch1882}).  

\subsection{Playfair's axiom}
 ``Playfair's axiom'' is the version of the parallel axiom adopted 
 by Hilbert in \cite{hilbert1899}.   That version, unlike all the other versions, makes 
 no existence assertion at all, but only asserts that there cannot exist two different lines
 parallel to a given line through a given point.
 
 \begin{figure}[ht]
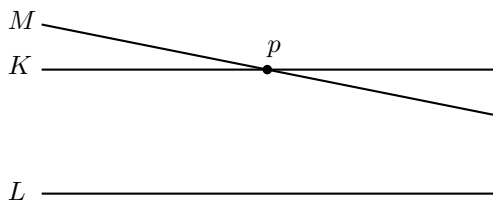
    
 \center{\PlayfairFigure}
\caption{Playfair:  if $K$ and $L$ are parallel, $M$ and $L$ are not parallel. }\label{figure:PlayfairFigureOne}
\end{figure}
 
\noindent
The conclusion of Playfair's axiom is that  $M$ and $L$ are not parallel.
By definition, parallel lines are lines that do not meet, so the conclusion is that $M$ and $L$
cannot fail to meet.   That is, not not there exists an intersection 
point.   Since $\neg \neg\, \exists$ is equivalent to $\neg \forall \neg$,
no existential quantifier is needed to express Playfair's axiom.

\subsection{Brouwer's intersection theorem}
\begin{figure}[h]
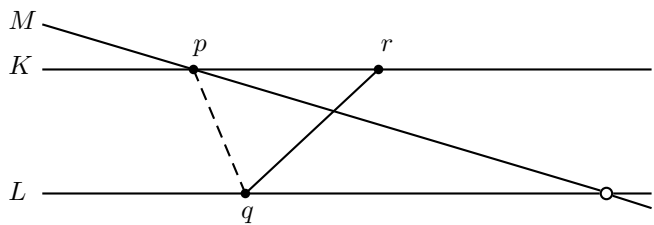
 
\center{\BrouwerParallelFigure}  
\caption{Brouwer's intersection theorem:  $M$ and $L$ must meet, provided  $pq$ makes
alternate interior angles equal with $K$ and $L$.}
\label{figure:BrouwerParallelFigure}
\end{figure}

Brouwer considers a version of the parallel postulate similar to 
Euclid's,  but without the ``witness'' $a$ in the hypothesis testifying to 
the side on which the angles are less than a right angle. 
It is similar to a variant of Playfair's axiom in which 
 ``$M$ and $L$ cannot fail to meet''  is replaced by  ``$M$ and $L$ meet.''  Rather
than calling it a ``postulate'',  Brouwer refers to the ``intersection theorem of Euclidean plane
geometry'', which he states as ``a common point can be found for any two lines $a$ and $\ell$
in the Euclidean plane which can neither coincide nor be parallel.''  
See Fig.~\ref{figure:BrouwerParallelFigure}.

This version of the parallel postulate was rediscovered (or re-invented?) in
\cite{beeson2010, beeson2016b},  where it was called the ``strong parallel postulate.''   There it was introduced in an axiomatic setting that 
included the stability of betweenness, or Markov's principle.
As we shall discuss below, this version of the parallel 
postulate,  in the absence of Markov's principle, needs a stronger
hypothesis to make constructive sense:  he should have required 
that $M$ and $K$ make a positive angle, i.e.,  are ``positively
non-collinear''.   Without that hypothesis, it is hardly surprising
that Brouwer's intersection theorem implies Markov's principle, for 
the hypothesis that $M$ and $K$ are unequal lines is a negative one,
but the conclusion that $M$ meets $L$ is a positive one.

\subsection{Euclid~5 formulated in Tarski's language}
\begin{figure}[ht]
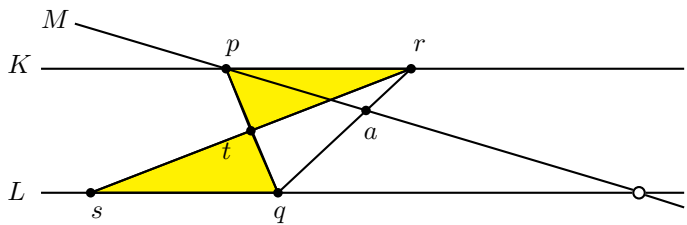

\caption{Euclid~5.  Transversal $pq$ of lines $M$ and $L$ makes corresponding interior angles less than 
two right angles, as witnessed by $a$. The shaded triangles are assumed congruent. Then $M$ 
meets $L$ as indicated by the open circle.}
\label{figure:EuclidParallelRawFigure}
\hskip 2.5cm
\EuclidParallelRawFigure
\end{figure}

 Euclid's version of the parallel postulate mentions angles, and the 
concept of ``corresponding interior angles'' made by a transversal.  
Here we give a formulation of Euclid's parallel postulate, expressed 
without mentioning angles.   Tarski's well-known axiomatization 
of geometry only talks about points;  angles are discussed 
indirectly, as triples of points; hence it is of interest to 
formulate Euclid~5 in that language. 
See Fig.~\ref{figure:EuclidParallelRawFigure}.   In order to eliminate 
the hypothesis about alternate interior angles, 
we  replace it by the hypothesis
that the shaded triangles in the figure are congruent.  Congruence of 
triangles just means that the corresponding sides are congruent.  This 
formulation of the axiom then mentions only points and the relations
of betweenness and congruence,  yet it is conceptually faithful to 
Euclid's formulation.  There are dozens of propositions that are 
classically equivalent to Euclid~5, and Tarski chose a different one 
of those to serve as his parallel axiom; but we wish to follow Euclid closely.

\section{Constructive Euclidean geometry}
\subsection{Is Euclid constructive?}
\label{section:3.1}
In constructive mathematics,  if one proves something exists, one has to show how to construct it. 
 In Euclid's geometry, the means of construction are not arbitrary computer programs, but ruler and compass.  Therefore it is natural to look for 
quantifier-free axioms,  with function symbols for the basic ruler-and-compass constructions.  The terms of such a theory
correspond to ruler-and-compass constructions.  These constructions should depend continuously on 
parameters, since we want to allow an interpretation of geometry in which 
the points are given by successive approximations.    We can see the continuity
of ruler-and-compass constructions  
dramatically in computer animations, in which one can select some of the original points and 
drag them, and the entire construction ``follows along.''
We expect that if one constructively proves that points forming a certain 
configuration exist, then the construction can be done ``uniformly'', i.e., by a single construction 
depending continuously on parameters.  

To illustrate what we mean by a uniform construction, we consider an important example.
There are two well-known classical constructions for constructing a perpendicular to line $L$
through point $p$:  one of them is called ``dropping a perpendicular'', and works when $p$ is not on $L$.
The other is called ``erecting a perpendicular'', and works when $p$ is on $L$.   Classically we 
may argue by cases and conclude that for every $p$ and $L$, there exists a perpendicular to $L$
through $p$.  But constructively, we are not allowed to argue by cases.  If we want to prove
that for every $p$ and $L$, there exists a perpendicular to $L$ through $p$, then we must give 
a single, ``uniform'' ruler-and-compass construction that works for any $p$, whether or not $p$ is 
on $L$.

The other type of argument (besides argument by cases) that is famously not allowed in constructive mathematics is proof
by contradiction.  There are some common points of confusion about this restriction.  The main thing one is not 
allowed to do is to prove an existential statement by contradiction.  For example,  we are not 
allowed to prove that there exists a perpendicular to $L$ through $x$ by assuming there is none,
and reaching a contradiction.   From the constructive point of view, that proof of course 
proves {\em something},  but that something is weaker than existence.  We write it 
$\neg\neg \, \exists$,  and constructively, one cannot delete the two negation signs.

However, when proof by contradiction is used only to prove that two points
are equal, or two segments are congruent,
 or that one point is between two others, the situation seems
qualitatively different.  There is no ``existential information'' missing
from such a proof.  Nothing is being asserted to exist, let alone being
asserted to exist without being constructed.  
The question then arises, whether such proofs are constructively 
acceptable. That is, whether one is or one is not
allowed to prove equalities, congruences, and 
inequality or order relations between points
by contradiction.  Consider for the moment two points $x$ and $y$ on a line.  If we 
derive a contradiction from the assumption $x\neq y$, then $x=y$.  
Formally this becomes the implication
$$\neg\, x\neq y \implies x=y.$$
  In the style of Euclid:
{\em Things that are not unequal are equal.} It has its intuitive grounding in the idea that 
$x=y$ does not make any existential statement.   This principle is named
``the stability of equality.''  (Generally any relation is ``stable''
if it is implied by its double negation.)

The principle that  $x < y$ can be proved by 
by contradiction can be expressed by
$$\neg\neg  \, x < y \implies x < y.$$ 
or equivalently
$$\neg\, y \le x \implies x < y.$$ 
 Since order on a line is 
not a primitive relation in geometry, we consider
 instead the corresponding axiom for the betweenness relation,
namely 
$$ \neg\neg  \, \B(a,b,c) \implies \B(a,b,c).$$ 
This is known as ``the stability of betweenness.''
It is also known as ``Markov's principle'',  since it 
was adopted by Markov as a basic principle of Russian constructive 
mathematics.

Here is a geometric way of thinking about Markov's principle,
rather than as a principle about ordering the continuum.
Markov's principle reduces to the assertion that,
 given two points $s$ and $t$, we can find two circles with centers $s$ and $t$
that separate the two points.  The obvious candidate for the radius is half the distance
between $s$ and $t$, so the question boils down to Markov's principle for numbers: if that 
radius is not not positive, must it be positive?  The stability of betweenness thus has 
the same philosophical status as Markov's principle:  it is self-justifying,  not
provable from other constructive axioms, and leads to no trouble in constructive mathematics,
while simplifying many proofs.  

Brouwer was not a person to accept a principle because it was 
useful, convenient, and apparently harmless, in the sense that it
does not interfere with the constructibility of 
solutions proved to exist with its aid.  Before Brouwer would accept
a principle, he wanted it to be true,
and he saw no reason why Markov's principle has to be true.

\subsection{The form of Euclid's theorems and proofs}
It is helpful to remember that Euclid did not work in first-order logic.
His theorems, and their proofs, have a fairly simple structure:
  Given some 
points, lines, and circles bearing certain relations,  then there exist some further points
bearing certain relations to each other and the original points.   This logical simplicity
implies that (although this may not be obvious at first consideration) 
if we allow the stability axioms for equality and betweenness, then essentially the only differences between classical and 
constructive geometry are the two requirements:
\begin{itemize}
\item You may not prove existence statements by contradiction; you must provide a construction.
\item The construction you provide must be uniform; that is, it must be proved to work without an argument by cases.
\end{itemize}

Sometimes, when doing constructive mathematics, one may use a mental picture in which one imagines
a point $p$ as having a not-quite-yet-determined location.  For example, think of a point $p$ which 
is very close to line $L$.  We may turn up our microscope and we still can't see whether $p$ is 
or is not on $L$.  We think ``we don't know whether $p$ is on $L$ or not.''   Our construction 
of a perpendicular must be visualized to work on such points $p$.  Of course, this is just a 
mental picture and is not used in actual proofs.  It can be thought of as a way of 
conceptualizing ``we don't have an algorithm for determining whether $p$ is on $L$ or not.''

We illustrate these principles with a second example.
Consider the problem of finding the reflection of point $p$ in line $L$.
Once we know how to construct a perpendicular to $L$ through $p$, it is still not trivial to 
find the reflection of $p$ in $L$.  Of course, if $p$ is on $L$, then it is its own reflection,
and if $p$ is not on $L$, then we can just drop a perpendicular to $L$, meeting $L$ at the foot $f$,
and extend the segment {\em pf} an equal length on the other side of $f$ to get the reflection.  But what about 
the case when we don't know whether $p$ is or is not on $L$?   Of course, that sentence technically 
makes no sense;  but it illustrates the point that we are not allowed to argue by 
cases.  The solution to this problem may not be immediately obvious; 
see \cite{beeson2016b}.

The description given above of the form of Euclid's theorems is supported
by Avigad {\em et.~al.} in \cite{avigad2009}: 
\begin{quote}
Euclidean proofs do little more than introduce objects
satisfying lists of atomic (or negation atomic) assertions, and then draw further
atomic (or negation atomic) conclusions from these, in a simple linear fashion.
There are two minor departures from this pattern. Sometimes a Euclidean proof
involves a case split; for example, if $ab$ and $cd$ are unequal segments, then one
is longer than the other, and one can argue that a desired conclusion follows in
either case. The other exception is that Euclid sometimes uses a {\em reductio}; for
example, if the supposition that $ab$ and $cd$ are unequal yields a contradiction
then one can conclude that $ab$ and $cd$ are equal.
\end{quote}
These arguments are constructively acceptable, if we have the 
stability of congruence and betweenness.  
In \cite{beeson2016b},  several examples such arguments 
 in Euclid are examined, including I.6 and I.24. 
If the conclusion
is a congruence or equality statement, even an argument based on
``of two unequal segments, one is longer than the other'' does not 
require Markov's principle, but only the stability of congruence.

\subsection{Euclidean geometry with Markov's principle}
In \cite{beeson2016b}, I advocated adopting the stability 
of equality, congruence, and betweenness
as axioms of constructive geometry.   In particular,  I argued that 
Euclid Books I-IV can be formalized using intuitionistic logic plus
those two stability principles,  using Euclid's version of the 
parallel postulate instead of Hilbert's,  and otherwise making no changes
in Hilbert's axioms.   In \cite{beeson2015b}, a version of constructive
Euclidean geometry is given that is based on Tarski's language and 
axioms, slightly modified.   In addition to using Euclid~5 instead of 
Tarski's parallel axiom, we found it necessary to use strict 
betweenness.   Hilbert used strict betweenness, but Tarski used 
non-strict betweenness.  To avoid confusion, in this paper we use 
$\B(a,b,c)$ for strict betweenness and $\T(a,b,c)$ for non-strict betweenness
($\T$ being the first initial of ``Tarski.'').  For details of the 
axioms see \cite{beeson2015b}; the essential idea is that after replacing
$\T$ by $\B$, we have to ``put back'' some axioms that Tarski originally
used, but later found clever derivations of from the remaining axioms,
using in an essential way the ``degenerate cases'' of those axioms where 
non-strict betweenness holds.

In these two long papers, I showed that Euclid~5 together with the 
stability of congruence and betweenness allow the formalization, not only 
of Euclid Books I-IV,  but also of the construction of coordinates, so 
that one can construct within geometry the field operations defined on 
a fixed line $L$  (taken as the $x$-axis).    

The consistency of Euclidean geometry with stability axioms,  at least
relative to classical Euclidean geometry, is obvious, since all the axioms
are classically valid.    
It may be more instructive to note that Euclidean geometry also has 
a model in the recursive reals, i.e. those real numbers given by 
recursive Cauchy sequences with a specified rate of convergence.  This 
model can be formalized in arithmetic, and the resulting interpretation
is, using standard recursive realizability, consistent with Markov's 
principle and Church's thesis for arithmetic.  Thus, also in the 
sense of recursive mathematics, there is nothing inconsistent about 
the stability of equality and betweenness. 

Julien Narboux pointed out that the stability of equality can be derived from 
the stability of congruence.  The proof is given, in the context of 
Tarski's theories, in Lemma~7.1 of \cite{beeson2016b}. The converse
is also easy to prove.   So there are
really only two stability axioms to consider: stability of congruence
and stability of betweenness.

\subsection{Angles and angle ordering}
As is well-known, betweenness is never explicitly mentioned in Euclid. 
 There are three ways that betweenness  occurs implicitly in the statements of 
 Euclid's propositions:  collinearity, 
ordering of segments, and ordering of angles.  Euclid takes these concepts as 
undefined,  and assumes (in the common notions) the basic properties of 
ordering.  Hilbert also took angles and congruence of angles as primitive notions, but unlike Euclid, he defined angle ordering.
Tarski defined all three notions.
 The  basic properties of angle ordering then become theorems.
These developments are spelled out in 
 \cite{schwabhauser}, Chapter 11,
with attention to constructivity in \cite{beeson2015b}, \S8.11.  
We here review that treatment to see whether and where Markov's principle was used.

The concept ``$x$ lies on $Ray(b,a)$'' is needed to define angles.  In 
\cite{beeson2015b}~\S8.11,  we defined
``$x$ lies on $Ray(b,a)$'' by $\neg (\neg \T(b,x,a) \land \neg \T(b,a,x))$.  
Here $\T$ is non-strict betweenness.  
That definition won't do if we do not have Markov's principle.  Instead, we 
should use the definition
$$ \exists e\,( \B(e,b,x) \land \B(e,b,a) \land eb = ea).$$
That is, $x$ is on the opposite side of $b$ from the
reflection $e$   of $a$ in $b$.  In particular, the reflection
$e$ must exist for $x$ to be on the ray.
Of course, if we assume Markov's principle,
then the two definitions are equivalent.

Then ``$abc$ and $ABC$ are the same angle'' means that the same points lie on $Ray(b,a)$
as on $Ray(B,A)$ and the same points lie on $Ray(b,c)$ as on $Ray(B,C)$.  Then  two angles $abc$
and $ABC$ are congruent if by changing $a$, $c$, $A$, and $C$
to other points on those same rays,  we can make $ab=AB$
and $bc = BC$ and $ac = AC$.

We say ``$f$ lies in the interior of angle $abc$ if 
there is a ``crossbar'' $uv$, with $u$ on $Ray(b,a)$
and $v$ on $Ray(b,c)$, with $u$, $b$, and $v$ distinct,
and for some point $e$ we have $\B(u,e,v)$ and $\B(b,e,f)$.
Then if $f$ lies in the interior of $abc$, it also lies
in the interior of any $a^\prime b c^\prime$ that is the 
``same angle'' as $abc$.

\section{Brouwer's 1949 paper}
\subsection{What Brouwer actually proved}
What  did 
Brouwer mean by asserting the inconsistency of Euclidean geometry? 
In this section we answer that question.
Brouwer worked extensively with ``the continuum'', which we denote
here by $\R$.  We refer to members of the continuum as ``real numbers'',
rather than using Brouwer's terminologies for numbers given by certain
kinds of sequences.  Geometry has a model in which the points are pairs of 
real numbers $(x,y)$.  We refer to this model as ``the model $\R^2$.''

Real numbers are given by sequences of rationals, and order in $\R$ 
is defined in terms of order in the rationals and the concept of sequence,
so it implicitly depends on the natural numbers.  Order in geometry, 
on the other hand, is defined in terms of betweenness on a line, which 
in turn is given by some axioms about betweenness.  In Brouwer's 1949
paper, he is concerned only with the model $\R^2$ of geometry, and not 
with axiomatic geometry.

Here is what Brouwer proved in \cite{brouwer1949b}:

\begin{theorem} [Brouwer]\label{theorem:brouwer} In the model $\R^2$, 
the statement ``the strong parallel postulate implies
Markov's principle'' holds.
\end{theorem}

\noindent{\em Proof}.  Let $L$ be the $x$-axis and $P$ be the point
$(0,1)$.   Assume $\neg\neg \ \epsilon > 0$.  
Let $a = (1,1-\epsilon)$.
 According to the strong parallel
principle, the line through $P$ and $a$ meets
$L$.   The point of intersection is $(x,0)$,
where by similar
triangles $x$ is to 1 as 1 is to $\epsilon$.   (The existence
of $x$ is guaranteed by the strong parallel postulate, not by division
by $\epsilon$, which would not be legal under the 
weak assumption $\epsilon \neq 0$.)  Now Brouwer appeals to 
the fact that $\R$ satisfies the Archimedean axiom:  there is a positive  
integer  $N$, which Brouwer chooses in the form $2^{n_\ell}$,   greater than $\vert x \vert$.  
 Now we have $x < N$, and both sides of 
the inequality have multiplicative inverses.    Brouwer says (line 12 of his paper,
where his $\rho_\ell$ is our $-\epsilon$)
$$ \rho_\ell < -2^{-n_\ell}$$
or in our notation,
$$ \epsilon > \frac 1 N.$$
Brouwer offers no further justification, but we think some is necessary. 
Now $\epsilon = 1/x$, so to justify Brouwer's conclusion we would need to infer
$1/x  > 1/N$  from $x < N$ and the existence of $1/x$ and $1/N$.  We know $N > 0$,
but we only know $\neg\neg \  x > 0$; hence we do not know $xN > 0$ and the 
last step of the following
argument cannot be carried out
$$ \frac 1 x - \frac 1 N = \frac {N-x}{xN} > 0$$
since the last step would require $xN > 0$, which we do not have.
Instead, we get only
$$ \neg\neg \  \ \epsilon > \frac 1 N.$$
We can, however, justify Brouwer's argument using apartness.  The principle 
of apartness, which holds in $\R^2$, says that if $\alpha < \beta$, then 
$$z < \beta \lor \alpha < z.$$ 
It follows that if $\alpha > 0$ and $\neg\neg \  \ \beta > 0$, then $\alpha + \beta > 0$,
since by apartness, $\alpha + \beta >0 \lor \alpha + \beta < \alpha$, and 
the second case, $\alpha + \beta < \alpha$, implies $\beta < 0$, which contradicts
$\neg\neg \  \beta > 0$.   Applying this principle with $\beta = 1/x - 1/N$ and
$\alpha = 1/2N$,  we have 
$$ \frac 1 x - \frac 1 N + \frac 1 {2N} = \frac 1 x - \frac 1 {2N} > 0$$
and hence $\epsilon > 1/2N > 0$ as desired.
   That completes the proof.
 \medskip

 Brouwer used the Archimedean axiom  
to get a positive upper bound on $x$.
That is not necessary,  and neither is the use of apartness.
We now show how to remove these two principles from Brouwer's proof.
First, we observe that from the hypothesis $\neg\neg \  \epsilon > 0$
and the equation $x=1/\epsilon$, we have 
$\neg\neg \  x > 0$.   
From $\neg\neg \  x > 0$ we have $\neg\neg \  x = \vert x \vert$.
By the stability of equality, we have
$x = \vert x \vert$.

Now, instead of using Archimedes's axiom to get
a positive upper bound on $x$, we could just as well use $\vert x \vert + 1$. 
Then we have 
$$ 0 < \frac 1 {2(\vert x \vert + 1)} < \frac 1 {\vert x \vert} = \frac 1 x = \epsilon $$
and hence $\epsilon > 0$.  The use of Archimedes's axiom is a red herring
(for this proof--not necessarily in general for intuitionism).  That
Brouwer used it anyway shows how far from Brouwer's mind 
was any consideration of 
whether his argument was first-order or not.

\subsection{Implications of Brouwer's theorem for \\axiomatic geometry}
Brouwer's theorem is formulated as a theorem about the Euclidean plane
$\R^2$.   But for half a century already at the time of Brouwer's
publication, ever since Hilbert's 1899 book \cite{hilbert1899}, 
geometry had moved on from discussing the one true plane to axiomatic
formulations. 
In Brouwer's paper, he did not consider the question whether Markov's 
principle is contradictory in some axiomatic system for Euclidean 
geometry.  That may have been because of his aversion to axiomatics
in general,  or it may have been the lack of any development at all 
of constructive Euclidean geometry at that time,  or it may have 
been for some other reason entirely.  But now, the question seems 
natural.%
\footnote{The referee pointed out that Brouwer was the thesis advisor, 
two years after the publication
of the article on the contradictority of elementary geometry, of Johanna Adriana Geldof's thesis \cite{geldof}. This thesis has nothing 
specifically intuitionistic in it, except one sentence near the beginning
saying that it assumes any two elements are either equal or not equal,
and also nothing specifically Euclidean, but it does show that Brouwer's
aversion to axiomatics was not absolute.  
}
Can we use Brouwer's proof to show that the axioms of 
Euclidean geometry (as formulated for example in \cite{beeson2015b}, 
but without Markov's principle) allow one to deduce Markov's principle
from Brouwer's version of the parallel axiom?  

We showed above that Brouwer's use of the Archimedean axiom and 
his (implicit) use of apartness are easily eliminated.  The one 
remaining issue is his use of coordinate geometry.
What Brouwer's proof shows is that,  if some theory of Euclidean 
geometry suffices to define coordinates and arithmetic,  then in that 
theory, Brouwer's version of the parallel postulate implies Markov's
principle.  That is, after all,  not too surprising, given that the 
hypothesis of Brouwer's ``intersection theorem'' is the 
negative statement that the two lines $M$ and $K$ are not identical,
but the conclusion is a positive existence statement.  Markov's
principle is ``built-in.''

Later in this paper, we will give a formal theory EG for Euclidean geometry without 
Markov's principle.  The theory EG plus Markov's principle has been
extensively studied in \cite{beeson2015b, beeson2016b}.   The fact that 
an attractive theory of geometry can be formulated without using Markov's
principle as an axiom, and without any obvious way to prove Markov's principle,
leads to the following questions and conclusions:
\smallskip

(1) Does EG (which does not have MP) plus Euclid~5 suffice to define coordinates 
and geometric arithmetic?  For short we say ``EG can define arithmetic'' to 
describe this property.  We claim in this paper that EG can define 
arithmetic.
\smallskip

(2) Then Brouwer's proof shows that EG plus
 Brouwer's intersection theorem  proves MP.    
\smallskip

(3) We know \cite{beeson2016b} that EG + MP   proves Brouwer's
intersection theorem, which is there called the ``strong parallel 
postulate'', or SPP for short.  Hence SPP is 
actually equivalent to MP in EG.
\smallskip
 
HA is ``Heyting's arithmetic'',  the standard formal theory of 
arithmetic with intuitionistic logic.  In the context of HA, 
``Markov's principle'' is the name usually given to  
$$\neg \neg \,\exists x\, P \implies \exists x\, P \qquad \mbox{($P$ primitive recursive).}$$
It is well-known (see e.g. \cite{troelstra})
 that this principle is not provable in HA.
Consider the interpretation of EG in HA 
determined by representing points as 
pairs of (indices of) recursive real numbers.  To verify that the 
interpretations of the axioms
of EG are provable in HA, we do not need Markov's principle, 
since the stability of betweenness is not an axiom of 
 EG.  We also do not need Markov's principle to 
verify any of the other axioms of EG, including Euclid~5, 
since the coordinates of 
the point asserted to exist by Euclid~5 involve only 
positive denominators. 
To show that EG does not prove
the geometric Markov's principle (stability of betweenness) it therefore
suffices to prove that the recursive interpretation of the geometric
Markov's principle is equivalent to the arithmetic Markov's principle.
But that is an easy exercise.  Hence
\smallskip

(4)  EG does not prove Markov's principle. 
\smallskip

  Since Brouwer's ``intersection
theorem'' is equivalent to Markov's principle, EG does not prove 
that ``theorem.''   In particular Brouwer's ``intersection theorem''
is not verifiable in the recursive interpretation (unless
we assume the arithmetic Markov principle). 

Brouwer never mentioned (in the paper under consideration, or anywhere else
as far as I know) 
 Euclid, or Euclid's axioms;  nor did he 
mention Hilbert, or Hilbert's axioms, or any axiomatic system whatever.
He worked simply with the plane $\R^2$ using coordinate geometry.
He showed that if $\R^2$ satisfies the ``intersection theorem'', 
then $\R$ satisfies Markov's principle.  The title of his paper,
however, claims that ``elementary geometry is contradictory.''
To reach that conclusion from his result, one would need to believe
that the intersection theorem is part of elementary geometry, 
and that Markov's principle is contradictory.   Of course 
the intersection theorem is a part of {\em classical} elementary 
geometry, but its proof requires Markov's principle, so it is not a 
part of {\em intuitionistic} elementary geometry, unless one 
assumes Markov's principle or an axiom that implies Markov's principle.
In particular, the intersection theorem in question does not follow
from Euclid's formulation of the parallel postulate.  

\subsection{Why did Brouwer reject Markov's principle?}
  Brouwer had claimed
in that same year (1949) in \cite{brouwer1949a} that Markov's principle
is contradictory.  His proof used real numbers that are limits of 
sequences generated by the ``creative subject'',  who is allowed to 
examine at each stage of mathematical construction, all proofs developed
at earlier stages.  Brouwer regarded this as an improvement over
a paper published the previous year \cite{brouwer1948a}, in which he 
showed that Markov's principle was ``unlikely to be provable.''
In view of these results, Brouwer viewed his version of the parallel
axiom as ``contradictory.''  

Here is a sketch of Brouwer's refutation of MP, in more modern terms.  He used
Kripke's schema (KS), according to which any proposition $\phi$ is 
equivalent to a proposition of the form $\alpha > 0$,  for some real number $\alpha$.
 Taking $\phi$
to be $A \lor \neg A$, and noting that $\neg\neg \ (A \lor \neg A)$ is 
intuitionistically valid, we have $\neg\neg \  \alpha > 0$, so by Markov's principle
$\alpha > 0$; that is, $A \lor \neg A$.  Hence KS implies the law of the excluded
middle.  But the fan theorem (uniform continuity of functions on $2^N$) refutes the 
law of the excluded middle.  Hence KS plus the fan theorem refutes MP.  Brouwer
believed, at the time of writing the papers we are discussing,
 that his theory of the creative subject justified KS, and hence, that MP 
had been refuted.  

Krike's schema has found few (if any) adherents in the 70 years since
Brouwer advocated it.  On the other hand,  Markov's principle is 
consistent with most commonly-studied intuitionistic theories; 
in particular with the theories of intuitionistic analysis including
Brouwer's fan theorem and bar theorem as axioms.  See \cite{troelstra344}
and \cite{kleene-vesley} for proofs of this consistency using variations
of recursive realizability. 

\subsection{Two sides}
In Brouwer's paper \cite{brouwer1949b}, he also considered
 an 
ordering principle that we call ``two-sides'':
$$ x \neq 0 \implies x < 0 \lor x > 0.$$
The name is chosen because the principle can be thought of as saying
that there are ``two sides'' of the $y$-axis:  every point not on 
the $y$-axis lies on the left half of the plane or on the right half.

We digress to show that two-sides can be expressed
in the language of geometry.  
Fix two points $0$ and $1$, and define
$-1$ to be the endpoint of the 
extension of the line segment
from 1 to 0 by itself.
Then $x < 0$ can be defined
as $\B(x,0,1)$. Two-sides can be expressed as
$$ x \neq 0 \implies \B(x,0,1) \lor \B(-1,0,x).$$

Brouwer rejected not only Markov's principle, but also two-sides.
Two-sides is not a theorem of EG, even with the help of Markov's
principle, as shown in \cite{beeson2015b}.  
Proof sketch: the axioms of EG can be expressed, after
introducing some function symbols, without using $\exists$ or $\lor$.
Then cut-elimination can be used to show that
 no disjunctive theorems can be proved (unless one of the disjuncts
can be proved.)   Two-sides does imply Markov's principle, since 
if we assume $\neg\neg \  x > 0$, then $x \neq 0$, so by two-sides,
$x > 0 \lor x < 0$;  but $x < 0$ contradicts $\neg\neg \  x > 0$,
so that case is ruled out, and we conclude $x > 0$.   
Thus two-sides is stronger than Markov's principle.  

\section{Axioms of non-Markovian  geometry }
Brouwer's objection to Markov's principle led us to consider whether 
Markov's principle is really necessary for Euclidean geometry. 
In this paper, we introduce a theory EG  
 of Euclidean geometry, with the stability
of congruence (and hence the stability of equality) but without 
the stability of  betweenness (which is also called Markov's principle).
EG has Euclid~5 for its parallel postulate,  so it 
corresponds closely to Euclid.  The axioms of EG are listed
for reference in \S\ref{section:axiomlist}, but we shall
introduce them gradually, with explanations. 

It probably does not matter whether we take
a Hilbert-type formulation or a Tarski-type formulation as in \cite{beeson2015b},
but for the sake of precision we must pick a specific list of axioms,  and 
it is far simpler to work with the simple language and short list of axioms 
of Tarski's theory.   As a starting point, we consider the constructive version of Tarski's
axioms given in \cite{beeson2015b}, with Markov's principle deleted from the 
list of axioms.  We also need to modify two other axioms (segment extension and Pasch)
to ensure that they do not immediately imply Markov's principle.
 We postpone the details of those modifications of the axioms to the next two sections.
 
EG has line-circle continuity:  a line through a 
point  (strictly) inside a 
circle meets the circle in two distinct points.    We also 
include circle-circle continuity as an axiom.%
\footnote{For our present purposes, there is little to be gained by trying to 
eliminate circle-circle continuity as an axiom, or in general, by 
trying to minimize the number of axioms, since our aim is simply to 
demonstrate the viability of intuitionistic geometry without Markov's principle.}  

What we seek, then, is an axiomatization of EG such that each axiom is equivalent
(using MP) to an axiom of intuitionistic Tarski geometry (as defined in \cite{beeson2015b}),
such that EG does not imply Markov's principle,  and yet suffices for the development 
of Euclidean geometry.   More specifically,  we want EG to satisfy these criteria:

\begin{itemize}
\item  EG suffices to prove the correctness of the uniform constructions 
given in \cite{beeson2015b}, namely, uniform perpendicular and reflection in a point.
This permits the assignment of coordinates $(x,y)$ to each point, 
 given two fixed perpendicular lines 
to serve as the $x$-axis and $y$-axis.

\item  EG suffices to prove that every pair $(x,y)$ occurs as the coordinates of 
some (unique) point.

\item EG suffices to define the (uniform) addition and multiplication of (signed)
points on the $x$-axis, and the construction of square roots of non-negative points.

\item EG suffices to formalize the arguments of Euclid Book I, and probably
Books II-IV as well.
\end{itemize}

The papers \cite{beeson2015b} and \cite{beeson2016b} established these facts for 
a theory including Markov's principle,  so the task here is to find a modified 
version of this theory that does not imply Markov's principle but still satisfies
the criteria listed above.

In particular, we give EG the parallel axiom Euclid~5  (just as in \cite{beeson2016b,
beeson2015b}  rather than the ``strong parallel axiom'' used by Brouwer.   It follows
that, in some sense,  Brouwer was criticizing a ``straw man'',  in that the 
parallel postulate that he found unsatisfactory is not actually Euclid's parallel 
postulate,  and Euclid~5 does not suffer from the flaw (if it is a flaw) that 
Brouwer pointed out.   The strong parallel postulate and Euclid~5 are equivalent 
in Euclidean geometry with Markov's principle (as shown in \cite{beeson2016b}),
but they are {\em not} equivalent if Markov's principle is dropped,  since 
(at least with the aid of the apartness axioms) the strong parallel postulate implies
Markov's principle, while Euclid~5 does not.

There is no philosophical advantage (for the present purposes) in 
trying to choose a minimal set of axioms for EG, and indeed there are 
reasons (discussed below) to be generous in taking more axioms than probably
are necessary.  We  assume two versions of Pasch, and both line-circle and 
circle-circle continuity.  We modify Pasch's axioms to avoid 
degeneracies that imply discontinuous dependence (as we did in 
\cite{beeson2015b}) and also to avoid 
near-degeneracies that imply Markov's principle. 
 The remarkable conclusion is that 
we can then derive Euclid Book I (and probably II-IV), coordinates, and arithmetic,  without Markov's principle, though we do need to add hypotheses
that angles are positive and vertices are distinct.

\subsection{Markov's principle and   betweenness}
Tarski's geometry, and its variant EG, includes a minimal set of 
axioms about betweenness.  
Remember that we use strict betweenness $\B$ rather than non-strict betweenness $\T$
as in Tarski.  If we have Markov's principle, then $\B$ and $\T$ are interdefinable.
But without Markov's principle, there is no apparent way to define $\B$ from $\T$,
so it is good that we took $\B$ as fundamental. 

Tarski's final theory \cite{tarski-givant} had only one betweenness axiom, known as (A6) or ``the identity
axiom for betweenness'':
 $$\T(a,b,a) \implies a=b.$$
In terms of strict betweenness, that becomes $\neg \B(a,x,a)$,  or otherwise expressed,
$\B(a,b,c) \implies a \neq c$.   We also refer to this axiom as (A6).   
The original version of Tarski's theory had  more betweenness axioms (see \cite{tarski-givant}, p. 188).
These were all shown eventually to be superfluous in classical Tarski geometry, through the work of Eva Kallin, Scott Taylor,
Tarski himself, and especially Tarski's student H.~N.~Gupta \cite{gupta1965}.   These proofs 
appear in \cite{schwabhauser}. 
Here we give the axiom numbers from \cite{tarski-givant}, names by which they are known, and also the theorem numbers of their 
proofs in \cite{schwabhauser}:
\smallskip

\hskip-0.3cm
\axioms
$\T(a,b,c) \implies \T(c,b,a)$ &  (A14), symmetry, Satz 3.2 \\
$\T(a,b,d) \land \T(b,c,d) \implies \T(a,b,c)$ & (A15), inner transitivity, Satz 3.5a\\
$\T(a,b,c) \land \T(b,c,d) \land b \neq c \implies $ & \\
\qquad $\T(a,b,d)$ & (A16), outer transitivity, Satz 3.7b \\
$\T(a,b,d) \land \T(a,c,d) \implies$ & \\
\qquad $\T(a,b,c) \lor \T(a,c,b)$ & (A17), inner connectivity, Satz 5.3 \\
$\T(a,b,c) \land \T(a,b,d) \land a \neq b \implies$ & \\
\qquad $\T(a,c,d) \lor \T(a,d,c)$ & (A18), outer connectivity, Satz 5.1
\endaxioms

 Our theory of constructive geometry in 
\cite{beeson2015b}) has three betweenness axioms (numbered as in \cite{tarski-givant},
with an ``i'' added for ``intuitionistic''):
\medskip

\newaxioms{0}
$ \neg \B(a,b,a)$ & (A6-i) \\
$ \B(a,b,c) \implies \B(c,b,a) $ &(A14-i), symmetry of betweenness \\
$ \B(a,b,d) \land \B(b,c,d) \implies \B(a,b,c)$&(A15-i), inner transitivity
 \endaxioms

 In \cite{beeson2015b}, we appealed to G\"odel's double-negation
 interpretation to ``import'' the long, complicated proofs of 
the classically superfluous axioms A17 and A18.  (Of course
the conclusion has to be double negated to eliminate the 
disjunction.)  Since $\T$ is defined negatively,
that still works even without Markov's principle, for the versions
of these axioms using $\T$.  Versions of A17 and A18 in which 
the disjunction is not double-negated of course are not constructively
valid, and if we do double-negate the disjunction, then the double
negation interpretation still works with $\T$ replaced by $\B$.

The following ``connectivity principle'' is
equivalent to A17, and it is what is actually needed for several
very basic lemmas,  such as ``two lines intersect in at most one point.''
\medskip

{\em If $B$ and $C$ are both between $A$ and $D$,  and neither
is between $A$ and the other,  then they 
are equal.}
\medskip

By applying the double negation interpretation as described above
to the complicated proof of A17 in \cite{schwabhauser}. 
Therefore, theoretically, there is no need to add this formula as 
an axiom.  However, few readers are going to check the proof in 
\cite{schwabhauser} and verify the double negation interpretation,
so we might as well just ask them to accept the connectivity principle
as an axiom.  We therefore include it as an axiom of EG.

The  question remains, did we really add enough betweenness
axioms?  The answer is, yes we did, because we can successfully define
coordinates and arithmetic, and with that geometrically defined arithmetic,
the points on a line form a Euclidean field. The axioms for 
Euclidean fields (without assuming Markov's principle) are 
discussed in \S\ref{section:EuclideanFields}. 

\subsection{Distinct points and segment extension}
Tarski's segment extension axiom provides for an extension of segment $ab$ 
by segment $cd$; the result is a point, sometimes written $ext(a,b,c,d)$,
about which the axiom asserts, if $p=ext(a,b,c,d)$, that $bp = cd$ 
and $\B(a,b,p)$.  Tarski's theory has no condition on $a$ and $b$,
but in constructive geometry (with Markov's principle), as discussed in \cite{beeson2015b},
we have to require $a \neq b$; that is, only non-null segments 
can be extended.

We will show next that the (unmodified) 
extension axiom implies Markov's principle.  In preparation 
we review two facts, the uniqueness of extension and the
stability of equality. 
The extension of segment $ab$ by segment $pq$ is unique: if there
are two points $c$  such that $\B(a,b,c)$ and $bc=pq$ then those two 
points are equal, as can be proved using other axioms of geometry.
The principle of ``stability of equality'' was introduced 
in \S~\ref{section:3.1}:  $\neg \neg\, x=y \implies x=y$.
 
Now we prove that the extension axiom implies Markov's principle. 
Suppose $\neg\neg \  B(a,b,c)$; then $a \neq b$, so $p=ext(a,b,b,c)$ is 
defined; but $bc = bp$ and $\B(a,b,p)$ by the extension axiom.
Using the uniqueness of extension (doubly negated), 
from $\neg\neg \  B(a,b,c)$ we can prove $\neg\neg \  p=c$.  Then, 
 by the stability of equality, $p=c$. Then from $\B(a,b,p)$
we have $\B(a,b,c)$, the conclusion of Markov's principle.

Therefore,  if one wishes to do geometry without Markov's principle,
one must modify the extension axiom.  We do so by restricting the applicability of segment extension to 
positive segments.  That concept is defined
as follows:
We say segment $ab$ is {\em positive},  or equivalently, that $a$ and $b$ 
are {\em distinct}, if  
$$ \exists e\,\B(e,a,b) \lor \exists e\, \B(a,b,e) \lor \exists e\, \B(a,e,b).$$
  We write this 
as $a \# b$.   Since the word ``length'' carries 
the idea of measurement by numbers,  we avoid
the phrase ``has positive length'',  but it may be 
helpful to mention it once.  

Once $a \# b$ is defined, the segment extension axiom has additional
subtleties.  We said above that we want the axiom to say that 
positive segments can be extended.  But we also want 
the notion ``positive segment'' to be closed under congruence;
that is, 
$$ab = cd \land a \# b \implies c \# d.$$
We wish this principle to be a theorem. 
We therefore need to have our segment extension axiom say,
``any segment congruent to a positive segment can be extended,''
rather than just ``any positive segment can be extended.''

A second subtlety arises from the distinction between 
$\T$ and $\B$.  In constructive geometry with Markov's principle we could 
just say 
$$ a \neq b \implies \T(a,b,ext(a,b,c,d))$$
and then it can be proved that
$$ a \neq b \land c \neq d \implies \B(a,b,ext(a,b,c,d)).$$
But in non-Markovian geometry, we need to change $\neq$ to $\#$
and take {\em both} axioms.  The following show what we mean,
but they still deal only with the second subtlety and not the first.

$$ a \# b \implies \exists e\,(\T(a,b,e) \land be = cd). $$
$$ a \# b  \land c \# d \implies \exists e\,(\B(a,b,e)\land be=cd). $$
Of course, in the presence of Markov's principle, these are equivalent to 
the form (A4-i) used in \cite{beeson2015b}.

Combining the two solutions, we now give the final form of 
the extension axiom of EG:

$$ A \# B \land ab = AB \implies \exists e\,(\T(a,b,e) \land be = cd). $$
$$ A \# B \land ab = AB  \land C \neq D \land cd= CD \implies \exists e\,(\B(a,b,e)\land be=cd). $$

\begin{lemma}\label{lemma:distinctcongruence}
If $a \# b$ and $ab = cd$ then $c \# d$.
\end{lemma}

\noindent{\em Proof}. Suppose $a \#b$ and $ab = cd$.
There are three cases.  Case 1, for some $e$ we have $\B(a,b,e)$.
Then $b \# e$, so we can extend $cd$ by $be$ to a point $f$ with
$\B(c,d,f)$.  Hence $c\#d$.  Case 2, for some $e$ we have $\B(e,a,b)$,
is treated similarly.  Case 3, for some $e$ we have $\B(a,e,b)$.
Then first extend $ab$ by $ab$ to $e$; so $\B(a,b,e)$, reducing
to case 1.  That completes the proof.
\medskip

The symbol $a \# b$ is traditionally used for apartness, 
a concept introduced by Heyting.  We use the word ``distinct''
instead of ``apart'', because we do not 
include the traditional apartness axiom $a \# b \implies a \# c \lor b \# c$. 
The reasons why not are discussed in \cite{beeson2016b}.
We use ``unequal'' for the negative relation $a \neq b$, and 
``distinct''  or its synonym ``different'' for $a \# b$.
\medskip

{\em Remark}.  One might think that one could adopt one of the three
betweenness conditions as the definition of $a \# b$ and prove the other two 
to be equivalent.  That can be done using some lemmas about betweenness,
for example, to prove the third property, we would 
need   outer transitivity.  But the proof of outer transitivity
begins by extending a segment whose endpoints cannot be seen
to be distinct except by the third property itself; so that 
approach is circular and does not work. Therefore,  
we adopt the definition given instead of a shorter one. 
\medskip

It is often said that ``two points determine a line.''
In intuitionistic geometry, the correct statement is 
``two {\em distinct} points determine a line.''  A line
can be thought of as all the possible extensions of a 
segment.  Intuitively, if we do not know that two points are distinct,
the direction of the segment connecting them is not
clear and a line is not (yet) determined.
\medskip

The following fundamental principle of order on a line would 
have to be taken as an axiom, if it were not provable. However,
in non-Markovian geometry, we note that the hypothesis requires
$b$ and $c$ to be distinct, not just unequal.
\medskip

\begin{lemma}[Outer transitivity] \label{lemma:outertransitivity}
$$\B(a,b,c) \land \B(b,c,d) \land b \# c \implies \B(a,b,d).$$
and also 
$$ \B(a,b,c) \land \B(b,c,d) \land b \# c \implies  \B(a,c,d).$$
\end{lemma}

\noindent{\em Proof.} 
The first form follows from the second, 
 using the symmetry of betweenness. We prove the second form.
 Suppose the hypothesis. 
 Since $\B(a,b,c)$, we have $a \# c$.  Hence segment $ac$ can be extended
by $cd$ to point $e$.  Then by the 5-segment axiom, applied with 
$d$ at the top in Fig.~\ref{figure:TarskiFiveSegmentFigure}, and
$acd$ and $ace$  along the bases of the two figures, we conclude $ed = dd$.  
But  $ed=dd$ implies $e=d$, and by definition of $e$ and the 
extension axiom \hbox{A1-i}, we have $\B(a,c,e)$,
since $c \# d$ by hypothesis.
Hence $\B(a,c,d)$ as desired.
That completes the proof of the lemma.
\medskip

\subsection{The five-segment axiom and SAS}
\begin{figure}[h]
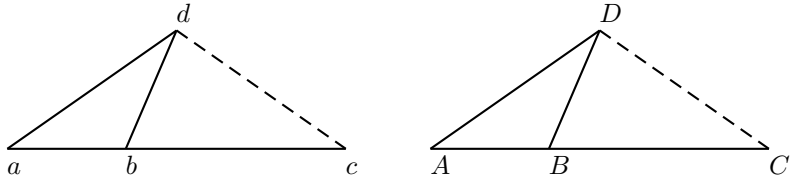
  
\caption{The 5-segment axiom. If $a \# b$ and the
corresponding solid segments are congruent, then $cd = CD$. 
\label{figure:TarskiFiveSegmentFigure}}
\TarskiFiveSegmentFigure
\end{figure}

Tarski replaced Hilbert's fourth and 
fifth congruence axioms (angle transport and SAS)  with an elegant axiom, known as the five-segment axiom.
This axiom is best explained not through its formal statement, but through   
Fig.~\ref{figure:TarskiFiveSegmentFigure}.  
The 5-segment 
axiom says that in Fig.~\ref{figure:TarskiFiveSegmentFigure}, the length of the dashed segment $cd$ is determined by the lengths 
of the other four segments in the left-hand triangle.  Formally, if the four solid 
segments in the first triangle are pairwise congruent to the corresponding segments in the second triangle, then the dashed segments
are also congruent.  

The 5-segment axiom is a thinly-disguised variant of the SAS criterion for triangle congruence.  To see 
this, refer to the figure.  The triangles we are to prove congruent are $dbc$ and
$DBC$.  We are given that $bc$ is congruent to $BC$ and $db$ is congruent to $DB$.   The congruence of angles $dbc$ and $DBC$ is expressed in
the 5-segment axiom by the congruence of triangles $abd$ and $ABD$, whose sides are pairwise equal.   The conclusion,
that $cd$ is congruent to $CD$,  give the congruence of triangles $dbc$ and $DBC$.  In Chapter 11 of \cite{schwabhauser}, one can find a formal proof of the 
SAS criterion from the 5-segment axiom.  It is easily adapted to prove 
the SAS criterion in non-Markovian geometry (which in non-Markovian geometry
requires that the angle in question be positive).

The formal statement
of the axiom uses non-strict betweenness and does not specify that 
$d$ is not collinear with $ab$, though it does require $a\neq b$.
Allowing $d$ to be collinear with $ab$ permits the use of this axiom 
to derive properties of betweenness,  and
 it is not only constructively valid, but very useful.   
 
The question arises whether we ought to strengthen 
the hypothesis $a \neq b$ to $a \# b$, or even  
perhaps require $\B(a,b,c)$ (which would imply $a \# b$, but also $b \# c$).
The latter would weaken the axiom.%
\footnote{The theorem ``all null segments are equal'' would not 
be provable; but adding it back as a new axiom would fix that.
However, there seems no point in doing that.
}  The former seems philosophically
prudent, since without it, perhaps the direction of the 
line $ab$ is uncertain.  We therefore require $a\# b$
in this axiom.%
\footnote{The unmodified 5-segment axiom  
holds in our Kripke model
where Markov's principle fails, so we are not
mathematically compelled to modify it.}

This axiom enables one to replace reasoning about angles with reasoning
about congruence of segments.%
\footnote{The history of this axiom is as follows.
The key idea (replacing reasoning about angles by reasoning
about congruence of segments) was
introduced (in 1904) by J. Mollerup \cite{mollerup1904}.
His system has an axiom closely related to the 5-segment axiom,
and easily proved equivalent.  Tarski's version, however, is 
slightly simpler in formulation.   Mollerup (without comment)
gives a reference to 
Veronese \cite{veronese1891}.  Veronese does have a theorem 
(on page 241) with the same diagram as the 5-segment axiom, and
closely related,  but he does not suggest an axiom related to this 
diagram.}  
We would like to emphasize that the 5-segment axiom is often 
just as easy to use as SAS.  Here is an illustrative example:

\begin{lemma} \label{lemma:verticalangles} Vertical angles are 
congruent.
\end{lemma}

\noindent{\em Proof}.  Let angles $abc$ and $dbe$ be vertical angles;
so $\B(a,b,e)$ and $\B(c,b,d)$.  We may suppose
without loss of generality that $bd = bc = ab = be$.  We 
must show $ac = de$.  Consider the two configurations
$abec$ and $dbce$.  The hypotheses of the 5-segment axiom are
fulfilled, because $ab=db$, $be=bc$, $ec=ce$, and $bc = be$.
Then by the 5-segment axiom, $ac = de$ as desired.  That completes
the proof.  The reader is urged to compare this proof to Euclid's proof.

\subsection{Lines, rays, triangles, and right angles}
The segment extension axiom allows us to create a new point on 
the line containing $a$ and $b$ only when $a$ and $b$ are 
distinct points. 
Without Markov's principle, we have to distinguish between pairs of 
unequal points, and pairs of distinct points.  In some sense only 
distinct points determine a segment, because only then does the 
segment extension axiom apply.  Although Tarski's language speaks
only about points, the import of the axiom is that we are allowed
to construct the line containing $a$ and $b$ only when $a \# b$.
Instead of ``two unequal points determine a line'',
we have ``two distinct points determine a line''.

The relation of collinearity is defined by 
$$ L(a,b,c)\ \leftrightarrow \ \neg\neg \ ( \B(c,a,b) \lor c=a \lor \B(a,c,b) \lor c=b \lor \B(a,b,c)).$$
Pushing one negation sign inwards, $L(a,b,c)$ can be expressed in a negative form.   Similarly, the notion that $x$ is on the ray from $b$ through $c$
can be expressed by $L(b,c,x) \land \neg \B(x,b,c)$.   We write this 
as $Ray(b,c,x)$.  The ray itself is not treated as a mathematical object;
only the relation between three points is formalized.
Nevertheless we speak of $x$ lying on $Ray(b,c)$, which
is the same as saying $Ray(b,c,x)$.

\begin{lemma} \label{lemma:layoff}
Let $a$ and $b$ be distinct points.  Let $c$ and $d$ be 
any points. Then we can lay segment $cd$ off on $Ray(a,b)$,
in the sense that there exists a point $x$ with $Ray(a,b,x)$
and $ax = cd$.  Moreover, point $x$ is unique.
\end{lemma}

\noindent{\em Proof}. Since $a\#b$, segment $ba$ can
be extended to point $e$ with $ae = ab$ and $\B(e,a,b)$.
  (We say $e$ 
is the reflection of $b$ in $a$.)  Then segment $ea$ 
can be extended by $cd$ to point $x$.  Then $ax = cd$ 
and $\T(e,a,x)$, by the segment extension axiom.  But
that is the meaning of $Ray(a,b,x)$, so the first claim 
of the lemma is proved.  It remains to prove the uniqueness
of $x$.  Suppose $y$ is another point on $Ray(a,b)$
with $ay = cd$.  By the definition of $Ray(a,b,x)$, that 
means for some point $e$ (the reflection of $b$ in $a$)
we have $\T(e,a,x) \land \T(e,a,y) \land ax = ay$.
Now apply the 5-segment axiom to the following
 two configurations
in  the form of Fig.~\ref{figure:TarskiFiveSegmentFigure}:
the first has $eax$ on the base and $y$ on the top,
the second has $eay$ on the base and also has $y$ on 
the top.  The hypotheses are fulfilled because $ax=ay$.
The conclusion is $yx = yy$.  By Axiom A3,
$x=y$.   That completes the proof.

\begin{definition} \label{definition:rightangle}
Angle $abc$ is a {\em right angle} if $a$, $b$, and $c$ 
are distinct points, and if $d=ext(a,b,a,b)$ (so $d$ is the 
reflection of $a$ in $b$), then  triangles $abc$ and $dbc$ are congruent
triangles.%
\footnote{We note that
in \cite{schwabhauser}, the predicate $R(a,b,c)$ allows the case $b=c$,
but our definition of ``right angle'' requires three distinct points.}
\end{definition}
It can be proved that this notion respects  
 the congruence relation on angles.
Euclid's Postulate 4, that all right angles are congruent, is a 
theorem in Tarski's geometry.%
\footnote{It is a much {\em simpler} theorem in Hilbert's geometry,
because Hilbert takes as an axioms that an angle can be copied on 
a given side of a given line, and the copy is unique.  These facts
require non-trivial proofs from the axioms of EG or Tarski.}  
It can be proved   
(and without Markov's principle).  We do not have space in this
paper to give the proof, but we will outline it.  The details
can be found in \cite{schwabhauser}, Part I, Chapter 10. 
\medskip

Reflection in a point is an isometry.  That is, if $B$ is 
the midpoint of $AC$ and also of $PQ$, then $AP=CQ$.  
The relation ``$ABC$ is a right angle'' is preserved under
reflection in a point.
Reflection in a line is an isometry.
And the key result (Satz 10.12 in \cite{schwabhauser}):
   Given right angles  $ABC$ and $ABF$ with $BC$ congruent to $BF$,   
   then $AC$ is congruent to $AF$.
(Remember there is no dimension axiom, so picture the two angles in different planes.)  Here is a sketch of the proof: 

Let $M$ be the midpoint of $CF$ (which exists without needing circles since triangle $BCF$ is 
isosceles).  Let $D$ be the reflection of $F$ in $B$, so $ABD$ is congruent to $ABF$.  Here
is the key:   triangle $ABD$ is the reflection of $ABC$ in the line $BM$.   
Since reflection in a line is a  isometry, those triangles are congruent,
proving Satz 10.12.  Now given two right angles, using Euclid's Prop.~23
we can copy one of them into the position of the second angle in Satz 10.12,
and applying Satz 10.12, the angle are congruent.

\begin{lemma} \label{lemma:righttwist}
If $abc$ is a right angle, then $cba$ is a right angle.
\end{lemma}

\noindent{\em Proof}.  The two angles obtained by reflecting
$abc$ first in $ba$ and then in $bc$ are vertical
angles, hence congruent.  In fact the proof of Lemma~\ref{lemma:verticalangles}
works directly, without needing to first prove that an angle congruent 
to a right angle is a right angle.

\subsection{Positive angles}\label{section:NL}

Corresponding to the notion of 
``distinct points'',  there is, intuitively, a notion of ``positive angle'',
which we write $abc > 0$  (it is not necessary to write $\angle abc > 0$).
We also write this as $0 < abc$; both our abbreviations for a statement
involving betweenness (given below).
Saying $abc > 0$ is stronger than simply requiring that $a$, $b$, and $c$ are 
not collinear (with $a$ and $c$ on the same side of $b$)
  in the same way that apartness is stronger than inequality.
It is our immediate aim to define this notion.   We wish this 
notion to be defined in such a way that, after we prove that
coordinates can be geometrically defined,  the angle between $p=(x,y)$
and $q=(b,0)$ with vertex at $(0,0)$ is a positive angle if and only 
if $y > 0$.%
\footnote{Here $y>0$ refers to order on a line, defined axiomatically 
by betweenness, namely $\B(-1,0,y)$, where $-1$ and $0$ are 
two specified points on a specified line.} 
But this cannot be the definition, as many theorems must
be proved before we can construct coordinates.  For example, 
midpoints and perpendiculars are needed.  

We must define ``positive angle'' at the outset, because
(as we shall see in the next section) we need to restrict the 
hypotheses of Pasch's axiom by requiring an angle to be positive.
Our first use of Pasch's axiom will be to repair Euclid's 
constructions in Propositions I.10 to I.10, culminating in the 
theorem that every segment $ab$ with $a \# b$ has a midpoint.
In that proof, we need to apply Pasch's axiom in situations
where the vertex angle is $30^\circ$, $60^\circ$, and $120^\circ$.
Therefore, the definition of ``positive angle'' must
immediately imply that those angles are positive.  More generally,
we need to know that the angles of a right triangle are positive.
(Remember that ``triangle'' implies the vertices are distinct.) 
\medskip

 The following conversations, following the pattern
of ``light bulb'' jokes, illustrate the basic concepts of 
non-Markovian geometry:

\begin{quote}  How many points does it take to determine a line?
{\em Three, two to lay a straightedge on and one to check that 
the other two are distinct.}

How many points does it take to determine a positive angle?
{\em Six, three to determine the vertex and sides, one for 
each side to check the other two are distinct,  and the sixth
to check that the sides don't coincide.}
\end{quote}
\medskip

\begin{definition} \label{definition:positiveangle}
(i)  $abc$ is an {\em apex angle} if there are distinct points $u,v$ on 
$Ray(b,a)$ and $Ray(b,c)$ respectively such that $bu = bv$ and $b \# u$.
\smallskip

(ii) $abc$ is an angle of a right triangle if $bac$ or $bca$ is a 
right angle, or more generally, if angle $abc$ is congruent
to an angle $ABC$ such that $BAC$ or $BCA$ is a right angle.
``Right angle'' is defined in Definition~\ref{definition:rightangle}).

\smallskip

(iii) $abc > 0$ (``$abc$ is a positive angle'') 
 if $abc$ is an apex angle, or a right angle, or an angle 
of a right triangle.
\end{definition}

We will later prove that every positive angle is an apex angle,
 but that can be done only after developing some geometry,
using the definition above.  In other words, if we were to 
take ``apex angle'' as the definition of ``positive angle'',  
we could not prove that an angle of a right triangle is positive, 
because to do so we need some intermediate results that cannot 
be justified until we know that an angle of a right triangle is positive.
Specifically, we will see below that a correct formulation of Pasch's
axiom for non-Markovian geometry requires certain angles to be positive,
and to prove that angles of a right triangle are apex angles, we need
such angles to be positive.  

The formula $abc > 0$, written out in primitive notation, is
existential, since not only is there an explicit $\exists$ 
and an explicit disjunction, but also
distinctness involves an existential quantifier.   This is good, since we 
intend it to express the existence of a positive lower bound on 
the angle $abc$.

We also need to express $abc < \pi$, i.e.,  angle $abc$ is 
positively different from a ``straight angle.''  We just need
to say that the supplement of $abc$ is a positive angle.
\begin{definition} \label{definition:positivesupplement}
Angle $abc$ has a positive supplement, or $abc < \pi$, is 
defined by
$$ abc < \pi \leftrightarrow \exists d\, (\B(a,b,d) \land dbc > 0).$$ 
\end{definition}

\noindent
We emphasize that the notation $abc > 0$ does not imply the assignment
of a measure of any kind to angles.  It is just a statement that 
a point can be found to witness, using the betweenness relation, that 
the angle is not zero.  To express that angle $abc$ is both
positive and has a positive supplement, we abbreviate the 
conjunction of the two statements $0 < abc$ and $abc < \pi$ as  $0 < abc < \pi$.

We can define certain specific angles as follows:
$$ abc = 60^\circ \leftrightarrow ab=bc \land ab = ac \land \#(a,b,c)$$ 
$$ abc = 120^\circ \leftrightarrow \exists d\,
        ( \B(a,b,d) \land ab = bd \land bc = bd \land cd = bc \land \#(a,b,c))$$
$$ abc = 30^\circ \leftrightarrow \exists d\,
        ( \B(a,d,b) \land ad = cd \land ad = db \land \#(a,b,c))$$
$$ abc = 150^\circ \leftrightarrow \exists d\,
        (\B(a,b,d) \land cad = 30^\circ)$$

We would like to prove that each of these four angles is a positive
angle.  That is neither obvious nor easy.  In fact, to do so 
it is necessary to have a strong lower dimension axiom.  That
issue will be discussed below; but the lower dimension axiom 
that we use guarantees the existence of an equilateral triangle
whose sides have midpoints, whose altitudes have distinct endpoints,
and whose medians meet in a central point.  That configuration
directly implies that $60^\circ$ is an apex angle (hence positive),
and $30^\circ$ is an angle of a right triangle.

\begin{lemma} \label{lemma:sixtydegrees}
 If $abc = 60^\circ$ or $abc = 120^\circ$ or $abc = 30^\circ$
or $abc = 150^\circ$ then $abc > 0$.
\end{lemma}

\begin{figure}[ht]
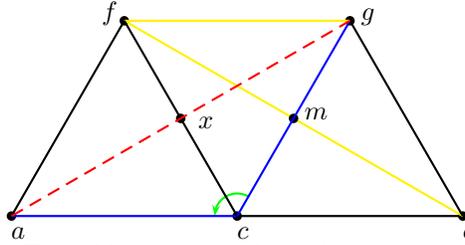

\caption{The $120^\circ$ angle $acg$ is positive, since $\B(a,x,g)$.}
\label{figure:SixtyDegreesFigure}
\center{\SixtyDegreesFigure}
\end{figure}

\noindent{\em Proof}.  The cases of $30^\circ$ and $60^\circ$
follow directly from the lower dimension axiom.  We take 
up the case $abc = 120^\circ$.
 We will show that the ``tiling'' 
shown in Fig.~\ref{figure:SixtyDegreesFigure} can be 
constructed. Then $acg = 120^\circ$ and $acg$ is an apex
angle because $ac = cg$ and $\B(a,x,g)$. 
\smallskip

{\em Remark}. It is 
not trivial to construct the three equilateral triangles 
shown in the figure and show that $\B(a,c,e)$. This has
to be done directly from the axioms,  and Pasch cannot 
be used since we as yet have no positive angles.  Euclid~5
will have to be used, as the construction does not work 
in non-Euclidean geometry. 
\smallskip

We start with the equilateral triangle {\it fac} and the midpoint
$x$ of side $fc$,  whose existence has to be assumed in 
the lower dimension axiom discussed in \S~\ref{section:dimension}
below. 
Triangle {\it fxa} is congruent to triangle $cxa$, so angle $axc$
is a right angle.  Then $cxa$ is also a right angle,
by Lemma~\ref{lemma:righttwist}.
 Hence $gc = ac$, by definition of right angle. 
  Then $ac = ca = fa = gc$. Now $cxg$ is a right angle,
so $gxc$ is also a right angle, so $fg = gc$.    
Let $m$ be the midpoint of $gc$. 

Consider the two lines $ac$ and {\it fg}.  The transversal
{\it fc} makes alternate interior angle equal, as witnessed 
by $f$, $g$, $x$, $a$, and $c$. 
Technically, our formulation of Euclid~5 does not mention angles 
at all; only the congruence of triangles {\it fgx} and $cax$ is needed.
The condition $\B(g,m,c)$ fulfills the hypotheses of Euclid~5
(which would traditionally be expressed as
``the corresponding interior angles {\it mfc} and {\it ecf} make less
than two right angles''). 
  Hence $fm$ meets the line through $ac$
in a point $e$.  That is, $\B(a,c,e) \land \B(f,m,e)$.

Because $m$ is the midpoint of $gc$, angles {\it fmg} and {\it fmc}
are right angles.  Therefore $cme$ and $gme$ are also right angles.
Since {\it cmf} is a right angle, $ce = fc = gc$.  Since $emc$ is a
right angle, $ge = ce$ and triangle $gce$ is equilateral.  
Then angle $acg$ is, by definition, a $120^\circ$ angle. 
But $a \# q$ because $\B(a,x,q)$.   That completes the 
case of a $120^\circ$ angle.
 $fmc$ is a right angle,  and $mf = me$, also $fc = ce$.
 
\begin{figure}[ht]
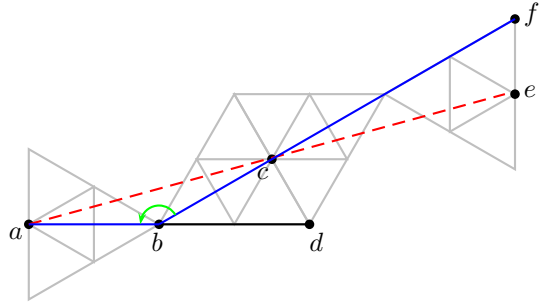

\caption{$abc = 150^\circ$ because $cbd = 30^\circ$, and $a \# c$ 
because $\B(a,c,e)$.}
\label{figure:150}
\center{\OneFiftyFigure}
\end{figure}

Now consider a $150^\circ$ angle.  Fig.~\ref{figure:150}
illustrates the construction.  

The illustrated network of 
equilateral triangles can be proved to exist using the 5-segment
axiom and Euclid~5, starting from point $b$,
as in the $120^\circ$ proof above. Then
$abc$ is a positive angle, since $\B(a,c,e)$ shows that $a \# c$.
But $abc$ is a $150^\circ$ angle, since the witnesses called for 
in the definition of $150^\circ$ angle are provided by $d$ and $e$.
points.  That completes the proof of the lemma.
\smallskip

{\em Remark}.  Fig.~\ref{figure:150} illustrates the 
statement that it takes six points to determine a
(positive) angle in non-Markovian geometry.  Point
$b$ is the vertex; point $d$ witnesses that $ab$ determines
a ray; point $f$ witnesses that $bc$ determines a ray;
point $e$ witnesses that those two rays are distinct.
\FloatBarrier

\begin{definition} \label{definition:triangle}
 A {\em triangle} is an ordered triple
of distinct points $abc$ such that all three angles $abc$,
$bca$, and $cab$ are positive angles with positive supplements.
\end{definition}

\subsection{Non-Markovian inner Pasch}
\begin{figure}[ht]
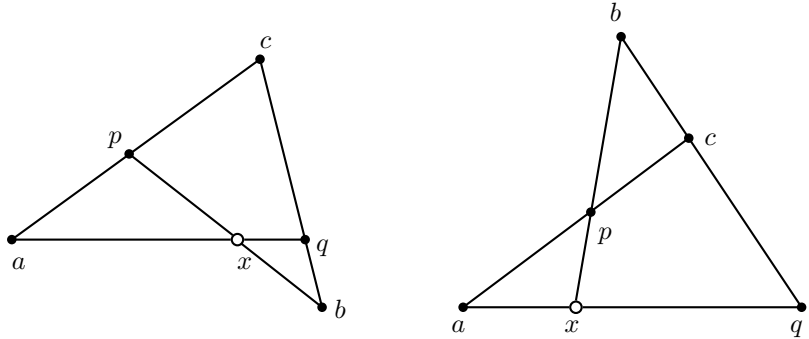
  
\caption{ Inner Pasch (left) and outer Pasch (right).  Line $pb$ meets triangle $acq$ in one side $ac$, and meets an extension of side $cq$.  Then 
it also meets the third side $aq$.
 The open circles show the points asserted to exist. }
\label{figure:InnerOuterPaschFigure}
\InnerOuterPaschFigure
\end{figure}

Tarski used (at various times) ``outer Pasch'' and ``inner Pasch'', both 
formulated using non-strict betweenness $\T$ and allowing various degenerate cases
that are not constructively valid.  These axioms are  
illustrated in   Fig.~\ref{figure:InnerOuterPaschFigure}.
As axioms they have two important advantages over other forms of Pasch:
they are purely existential, and they do not depend on the dimension.
That is, they hold in $\R^n$ for any $n$, in contrast to the version
of Pasch that says, if a line enters a triangle it must exit, which 
fails in $\R^3$.
 
 In \cite{beeson2016b}, we formulated what 
seemed the most general constructively valid version of inner Pasch, replacing $\T$
by $\B$ in the  hypothesis $\T(a,p,c)$ and in the two conclusions, and requiring points $a,b,c$ in Fig.~\ref{figure:InnerOuterPaschFigure}
to be not collinear.  If Markov's principle is dropped, and inner Pasch is not 
modified, Pasch's axiom will imply Markov's principle.

\begin{theorem} \label{theorem:PaschImpliesMarkov}
Inner Pasch without modifications implies Markov's principle.
\end{theorem}

\noindent{\em Remark}.  The theorem refers to the version of 
inner Pasch described above, and used in \cite{beeson2016b}.
\medskip

\noindent{\em Proof}.
Let $p$ and $v$ be any two distinct points, 
and 
assume unrestricted inner Pasch. (In Fig.~\ref{figure:InnerOuterPaschFigure},
$v$ would lie on $pb$ extended; think of $v$ fixed, while $b$ 
can move along the line $vp$ towards $p$).  We will prove 
$\neg\neg \,  \B(p,b,v)$ implies $\B(p,b,v)$.
To that end, assume $\neg\neg \,  \B(p,b,v)$; we have to prove $\B(p,b,v)$.
Let $a$ and $c$ be points with $\B(a,p,c)$ such that $\neg L(a,v,c)$
(for example, $a$ and $c$ can be taken to lie on the perpendicular to $pv$
at $p$.)
We have $\neg L(a,b,c)$, since $L(a,b,c)$ implies $b=p$ (by the 
definition of $L$).  Hence we can apply unrestricted inner Pasch,
obtaining the existence of $x$ with $\B(p,x,b)$.  Now we 
have $\B(p,x,b)$ and $\B(p,b,v)$. By the inner transitivity of 
betweenness, we have $\B(p,b,v)$.  That completes the proof.  
\smallskip

\noindent{\em Remark.} In that proof, we used the inner transitivity of betweenness and 
the existence of points $a$ and $c$ not collinear with $v$, but
these things can be proved from unrestricted Pasch, since as we 
show below, they can be proved even from restricted inner Pasch.  Here 
we only intend to demonstrate the need to restrict inner Pasch in some way.
\medskip

Our formulation of inner Pasch (suitable for non-Markovian geometry)
 uses strict betweenness $\B$ in all four places.
 So the hypotheses include 
$\B(a,p,c) \land \B(b,q,c)$, and the conclusion is 
$\exists x\,(\B(a,x,q) \land \B(p,x,b))$.  In a quantifier-free form,  $x$ is
replaced by a Skolem term.  We also need to make sure the whole 
figure does not degenerate into a line, by requiring the hypothesis
$0 < acb < \pi$.  

For comparison:  Tarski used $\T$ instead of $\B$, and did not care
if the figure degenerated.  For constructive geometry, that is wrong,
since it leads to discontinuous dependence of $x$ on the parameters.
In \cite{beeson2016b},  the hypothesis $\T(b,q,c)$ was retained but 
$\B$ was used for the other hypothesis and the two conclusions.
That is, the degenerate cases $q=c$ and $q=b$ are 
still allowed.  This is constructively sensible as the lines $aq$ and $pb$ still
are transverse, so have a unique intersection point, constructible with a ruler.
To avoid discontinuity, we also required $\neg L(a,b,c)$.   
 
 In non-Markovian geometry, 
the negatively-phrased non-collinearity is not enough, as the 
lemma above shows.  Instead we need $0 < acb < \pi$.
The hypothesis $0 < acb < \pi$ is expressed in the way 
defined in  \S~\ref{section:NL}.
That is an existential statement,  but it occurs only in the hypotheses
of inner Pasch, so inner Pasch is still equivalent to a quantifier-free
 axiom (when a Skolem term is used for $x$). 
 
There are, however, other possible hypotheses that we might also use
in connection with inner Pasch; that is, other possible sufficient 
conditions than $0 < acb  < \pi$ that should support the conclusion 
of inner Pasch without implying Markov's principle.  One such 
hypothesis is $0 < qpa < \pi$.  (Or equivalently, $0 < pqb < \pi$.)
 This may seem strange (after all 
line $pq$ is not even drawn in the diagram for Pasch's axiom), but 
it is the natural condition that we need when we use inner Pasch to 
prove Euclid's exterior angle theorem.  We were not able to prove 
the exterior angle theorem unless we take the form of inner Pasch 
with this hypothesis as an axiom, as well as the more natural $0 < acb < \pi$. 
Of course, after developing the theory of perpendiculars, we can prove
these conditions equivalent, but to get off the ground we need both 
versions of inner Pasch.

\subsection{Non-Markovian outer Pasch}
Please refer to Fig.~\ref{figure:InnerOuterPaschFigure}.
The hypotheses of unrestricted outer Pasch in 
\cite{beeson2016b} 
are $\B(a,p,c)$ and $\T(b,c,q)$, as well as $\neg L(a,b,c)$.  The 
conclusion is that there exists an $x$ such that $\B(a,x,q)$
and $\B(b,p,x)$.

Outer Pasch, like inner Pasch, 
in its unrestricted form implies Markov's principle.
The argument is similar to the one given for inner Pasch;
in the second part of Fig.~\ref{figure:InnerOuterPaschFigure}, we allow $a$ to move towards $q$, as for inner Pasch we allowed $b$ to move towards $p$.
Therefore, outer Pasch
 needs to have its hypotheses strengthened.  As for inner Pasch,
 we replace $\T$ by $\B$,  and replace the non-collinearity 
 hypotheses by requiring an angle to be between 0 and $\pi$. 
 But which angle?  
 
The hypothesis we choose to formulate non-Markovian outer
Pasch is this:
$$0 < baq < \pi \lor 0 < abq < \pi.$$
This choice is a pragmatic one:  it enables the uses of 
outer Pasch that we need (in particular for the crossbar theorem).
Eventually, we will prove that 
if $a$, $b$, and $c$ are 
three distinct points, and one of three angles formed is 
positive, so are the other two.  After that, 
the various possible versions   
of non-Markovian outer Pasch will be equivalent.  

Our theory EG includes both non-Markovian
inner Pasch and non-Markovian outer Pasch as axioms.
 Gupta proved (\cite{gupta1965}; see also \cite{schwabhauser}, Satz 9.3),
that inner Pasch implies outer Pasch; but this is a long development,
and although we did not find a use of Markov's principle,  we are 
not willing to certify that none is necessary.   As mentioned above,
we are not aiming in this paper  to find a minimal set of axioms for EG, 
but instead a sensible set of axioms that provides a smooth non-Markovian
constructive development
of Euclid.  Therefore we include both forms of Pasch. 

\subsection{Line-circle and circle-circle continuity}

We say that point $p$ is {\em inside the circle
with center $a$ passing through $c$}  provided
$p$ lies on a diameter of the circle.%
\footnote{It will not do to say $p$ is equal to the center
or lies on a radius, since we may not know which alternative holds.}
The line-circle continuity axiom says that if 
$p$ is inside the circle with center $a$ passing 
through $c$, then there exist two
points $x$ and $y$ on the circle with $\B(x,p,y)$.   
This axiom is true in $\R^n$ for any $n$, i.e.,
not just in plane geometry.  Without any dimension axiom,
``circles'' become spheres (in $R^3$) or hyperspheres.

We also need a non-strict version of line-circle continuity
in which $\B$ is replaced by non-strict betweenness $\T$.
Classically, the two are easily proved equivalent, but 
the case distinction whether $p$ is strictly inside or
on the circle is not legitimate intuitionistically.  Applied
to Descartes's  geometric construction of square roots, the 
non-strict version implies that non-negative segments have square roots,
and the strict version implies that positive segments have positive
square roots.  We therefore need both versions of line-circle continuity
to establish the existence of coordinates and define arithmetic.

The circle-circle continuity axiom says that if
$C$ and $K$ are circles with distinct centers, and 
 $C$ has a point non-strictly inside  $K$, and a 
point non-strictly outside circle $K$, then there is a point 
lying on both $C$ and $K$.  This immediately implies the 
corresponding strict version, but not conversely, so we 
take the non-strict version as our axiom.%
\footnote{We avoided ``degenerate cases'' of Pasch's axiom, in 
which a two-dimensional picture degenerates to one dimension. 
Using $\T$ in the circle-circle axiom is not a degenerate case
in this sense; when the circles are tangent, the picture is still 
two-dimensional.}

A stronger version of 
circle-circle continuity has the additional hypothesis that  
$p$ is any point not lying on the line $L$ containing
the centers of $C$ and $K$, and the additional 
conclusion that $C$ and $K$ have an intersection point
on the opposite side of $L$ from $p$;  that is, 
there are points $e$ and $x$ with $e$ on both $C$
and $K$ and $x$ on $L$ and $\B(e,x,p)$.   
We do {\em not} assume this stronger version of circle-circle
continuity as an axiom.%
\footnote{
  This means
that Euclid's proof of I.4 (angle bisection) is not directly
formalizable in EG, and his 
proof of the existence of a midpoint is consequently also not 
directly formalizable. 
Instead, we construct midpoints using Gupta's theorem about the existence of a 
midpoint of an isosceles triangle, combined with Euclid's
construction of an equilateral triangle on a given segment.
  Once midpoints are in hand,
perpendiculars can be constructed, and the stronger version of 
circle-circle continuity can
be proved. }

\subsection{Segment-circle continuity}
The segment-circle continuity principle says that if $p$ is a point
inside a circle, and $q$ is a point outside, then there is a point 
on the circle between $p$ and $q$.  This principle has been suggested
as an axiom by many authors, including Tarski (see \cite{tarski-givant}).
But a detailed study shows that it is inadequate; an irremovable
circularity arises in formalizing Euclid without a dimension axiom.
If we try to construct dropped perpendiculars (Euclid~I.12) using 
segment-circle continuity, to check the hypotheses we need
the triangle inequality (I.20).  But  I.19 is needed for I.20, and I.7
for I.19.  In Prop.~I.7, the two triangles that are supposed to coincide may lie
in different planes; that possibility has to be removed by an additional
hypothesis.  Even with its statement thus corrected, I.7 is more difficult to prove than Euclid thought, since he took for granted the fact that an angle
cannot be less than itself, but that principle is actually the essential 
content of I.7.  Ever since Hilbert \cite{hilbert1899}, angle inequality 
has been regarded as a defined concept, and proving I.7 then  
requires dropped perpendiculars (I.12).  But this is circular.
The conclusion is that segment-circle continuity is not a suitable axiom 
to use in formalizing Euclid.%
\footnote{Line-circle continuity does not suffer
from this problem, as the triangle inequality is not required to drop
perpendiculars.   Of course, as Gupta showed, one can construct 
dropped perpendiculars without mentioning circles at all, so there is
no formal result that one continuity axiom is better for I.7 than another,
as none at all is actually needed.  We merely say that Euclid's proof
can be repaired with line-circle, but not with segment-circle.
}

\subsection{Dimension axioms} \label{section:dimension}
When considering plane geometry, one needs a {\em lower dimension axiom}
providing for the existence of  three 
 non-collinear points.  Non-collinearity will not be 
enough in non-Markovian geometry; we need to ensure that there exist
three points forming a positive angle with a positive supplement.
In other words, we need three distinct points $\alpha$, $\beta$,
and $\gamma$ such that $0 < \alpha \beta \gamma < \pi$.%
\footnote{If we 
do not specify this in an axiom, nothing allows to prove any 
betweenness statement,  as the first deduction of a betweenness
statement cannot come from the betweenness axioms (which have 
betweenness hypotheses), nor from Pasch (which has betweenness
in the hypothesis about a positive angle), nor from the extension
axiom, which requires betweenness in the hypothesis about distinct points,
nor from circle-circle or line-circle continuity, where the hypotheses
about ``inside'' are expressed using betweenness.  Nothing would prevent
the plane from collapsing, with all points not not equal to each other.
Hence, it is reasonable to call this a ``lower dimension axiom''.
}

\begin{figure}[ht]
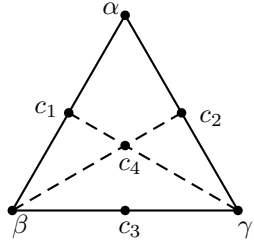

\caption{The lower dimension axiom}
\label{figure:lowerdimension}
\center{\LowerDimensionFigure}
\end{figure}

We take as an axiom that there are three points, given by 
constants $\alpha_1$, $\alpha_2$, $\alpha_3$,  that form an 
equilateral triangle.
Then we add a further axiom, introducing three more constants 
for the midpoints
of the three sides, and one final constant for the common 
intersection point of two of the medians. 
Thus the three vertices are
all distinct, and two of the medians have distinct endpoints.
See Fig.~\ref{figure:lowerdimension}. 
The formal expression of this axiom is given 
in the next section, where all the axioms are listed
for reference. 

The upper dimension axiom for plane geometry says that if $a$ and $b$
are distinct, and three points are each equidistant from $a$ and $b$,
then those three points are collinear.   If we do not use any 
upper dimension axiom,  or we replace 3 by some larger integer, 
then as we have already noted,  inner and outer Pasch still make sense,
and circle-circle continuity becomes sphere-sphere continuity.   
The notions ``same side'' and ``opposite side'' of a line still
make sense if properly defined.  See \cite{beeson2016b} for the definitions.

\subsection{List of axioms for reference}
\label{section:axiomlist}
In this section, we give the complete list of axioms of EG. 
First, we specify the language.  It is first-order predicate
calculus with equality.  There is a 3-ary relation symbol $\B$ 
and a 4-ary symbol $E$.  $E$ is the official version of the 
4-ary relation that we write informally as $ab = cd$; 
first-order predicate calculus as found in textbooks does 
not permit that syntax, so $ab=cd$ must remain as an
informal abbreviation for the more formal syntax $E(a,b,c,d)$.
\smallskip  

A number of other informal abbreviations are used in 
stating the axioms.  The actual axioms are the result 
of replacing these abbreviations by the right-hand side
of their definitions (recursively),
 so the axioms involve only equality,
$\B$, and $E$.  The following is a complete list of 
these abbreviations.  Since each right-hand side involves
only definitions earlier in the list, the recursive 
replacement mentioned does terminate.  That is, there is 
no circularity in this list of definitions.
\medskip

\begin{eqnarray*}
a \# b  &:=&   \exists c\, (\B(c,a,b) \lor \B(a,c,b) \lor \B(a,b,c)) \\
ab = cd &:=&  E(a,b,c,d) \\
R(a,b,c) &:=& \exists d\, (\B(a,b,d) \land ab =bd \land ac = dc \land a \# b \land c \#b \land a \#c)\\
Ray(a,b,x) &:=& \exists d\, (\B(d,a,b) \land \B(b,a,x))\\
apex(a,b,c) &:=& \exists u,v\,(Ray(b,a,u) \land Ray(b,c,v) \land bu=bv \land u \#v)\\
\angle abc= \angle ABC &:=& \exists U,V,u,v\,(Ray(b,a,u) \land Ray(b,c,v) \\
&& \land Ray(B,A,U) \land Ray(B,C,V) \land uv = UV) \\
0 < abc &:=& apex(a,b,c) \lor R(a,b,c) \lor \\
&& \exists A,B,C\,( \angle abc = \angle ABC \land  R(B,A,C) \lor R(A,C,B)) \\
 abc < \pi &:=&  \exists d\, (\B(d,b,a) \land 0 < dbc) \\
0 < abc < \pi &:=&  0 < abc \land abc < \pi \\
\T(a,b,c) &:=& \neg  (a\neq b \land \neg \B(a,b,c)\land b \neq c) 
\end{eqnarray*}
\medskip

Here are the betweenness axioms:
\medskip

\newaxioms{0}
$ \neg \B(a,b,a)$ & (A6-i) \\
$ \B(a,b,c) \implies \B(c,b,a) $ &(A14-i), symmetry of betweenness \\
$ \B(a,b,d) \land \B(b,c,d) \implies \B(a,b,c)$&(A15-i), inner transitivity \\
$\B(a,b,d) \land \B(a,c,d) \land $& \\
$\neg \B(a,b,c) \land \neg \B(a,c,b) 
\implies b = c$ & (A17-i), connectivity
 \endaxioms
 \medskip
 
Here are the axioms concerning congruence and betweenness:
\medskip

\newaxioms{0}
 $A \# B \land ab = AB \implies \exists e\, (\T(a,b,e)\land be = cd)$ & (A4-i1) extension \\
 $A \# B \land ab = AB \land C \# D  \land cd = CD $&\\
 $\implies \exists e\,(\B(a,b,e) \land be = cd)$ &(A4-i2) strict extension\\
 $a \# b \land ab = AB \land bc = BC \land ad=AD$ &\\
 $\land \ bd = BD \implies cd = CD$ & (A5-i) five-segment axiom \\
 $\B(a,p,c) \land \B(b,q,c) \land (0 < acb < \pi \lor 0 < qpa < \pi)$ &\\
 $\implies \exists x\,(\B(p,x,b) \land \B(a,x,q))$& (A7-i1) inner Pasch \\
 $\B(a,p,c) \land \B(b,c,q) \land (0 < baq < \pi  \lor 0 <  abq < \pi)$ &\\
 $\implies \exists x\,(\B(b,p,x) \land \B(a,x,q))$& (A7-i2) outer Pasch 
\endaxioms
\medskip

Here are the lower dimension axioms:
\begin{eqnarray*}
&& \alpha \beta = \beta \gamma \land \alpha\beta = \alpha\gamma \land \alpha \neq \beta  \\
&& \B(\alpha, c_1, \beta) \land \alpha c_1 = c_1 \beta \\
&& \B(\alpha, c_2, \gamma) \land \alpha c_2 = c_2 \gamma \\
&& \B(\beta, c_3, \gamma) \land \beta c_3 = c_3 \gamma \\
&& \B(\beta, c_4, c_2) \land \B(\gamma, c_4, c_1).
\end{eqnarray*} 
We have used constants $c_i$ and  
$\alpha$, $\beta$, $\gamma$; as a result the lower 
dimension axiom is quantifier-free and for convenience 
can be broken into several formulas, written on
different lines without being connected by $\land$.
  This is of no
significance--we might as well have used existentially
quantified variables, in which case the dimension axiom 
would need to be one long formula inside the scope of the 
existential quantifiers.

There is no upper dimension axiom.
\medskip

The first line-circle continuity axiom says 
that if point $a$ is inside
the circle with center $c$ and radius $pq$
(as witnessed by $a$ lying on a diameter $uv$),
and $b$ is a point distinct from $a$ (so that $ba$ 
determines a line),  then there are two points $x$ 
and $y$ on the circle, with $a$ between them.
\begin{eqnarray*}
 &&\B(u,a,v) \land cu = pq  \land cv = pq \land\ \B(u,c,v)  \\
 && \implies \exists x,y\,(cx = pq \land cy=pq \land \B(x,a,y))
\end{eqnarray*}

The second line-circle continuity axiom is similar, but 
with $\T$ instead of $\B$.
\begin{eqnarray*}
 &&\T(u,a,v) \land cu = pq  \land cv = pq \land\ \T(u,c,v)  \\
 && \implies \exists x,y\,(cx = pq \land cy=pq \land \T(x,a,y))
\end{eqnarray*}
 
Here is the circle-circle continuity axiom, 
which says that if the circle with center $c$ and 
radius $pq$ has a point $a$  non-strictly inside the circle with 
center $C$ and radius $PQ$,  and a point $b$ non-strictly outside
that circle, then the two circles meet.
\begin{eqnarray*}
&& ca = pq \land \T(u,a,v) \land \T(u,C,v) \land Cu = PQ \land Cv = PQ\\
&& cb = pq \land \B(C,w,b) \land Cw = PQ \\
&& \implies \exists e\,(ce = pq \land Ce = PQ)
\end{eqnarray*}

Here is Euclid~5.  See Fig.~\ref{figure:EuclidParallelRawFigure} for 
an illustration. 
\begin{eqnarray*}
&&pt = qt \land \B(p,t,q) \land st = rt \land \B(s,t,r) \land pr = qs \\
&&\land \ \B(q,a,r) \implies \exists e\,(\B(p,a,e) \land \B(s,q,e))
\end{eqnarray*}

{\em Logical form of the axioms}.  The occurrences of disjunction in 
the hypotheses of Pasch's axiom can easily be eliminated, using two formulas
for inner Pasch instead of one, and two for outer Pasch.  If the defined
symbols are replaced by their definitions, that creates existential quantifiers
in the hypotheses (which can be eliminated by simple logic) and in the 
conclusions (which can be eliminated by introducing function symbols).
Thus, there is an equivalent formulation in which the axioms are
disjunction-free and quantifier free.  This observation is important for 
certain meta-theorems discussed below.

\section{Development of non-Markovian geometry}
 In this section, we check the development of constructive geometry in  \cite{beeson2016b} and \cite{beeson2015b},  searching for
uses of Markov principle, and try to eliminate them when found.   In \cite{beeson2015b}, 
once having found an axiom system that supports the double-negation interpretation for geometry,
(for which the use of Markov's principle is crucial),  we were able to use it to ``import''
negative theorems from \cite{schwabhauser}.   If Markov's principle is rejected, we can 
no longer do that.
In order to reduce theorems to negative form, we had to find and prove the correctness of 
numerous ``uniform'' constructions, to eliminate arguments by cases;  those proofs necessitated
the direct formal development of a certain amount of geometry.  If we now want to eliminate
Markov's principle, it at first appears that 
the entire development of the two cited papers has to be checked, and in addition,  the parts of \cite{schwabhauser} that could, in the presence of 
Markov's principle, be ``imported'' via the double negation interpretation.

However, we found a way to avoid some of that work.  Namely, as
mentioned above, we include as axioms of EG both inner and outer Pasch,
and both line-circle and circle-circle continuity.
With Markov's principle, this is overkill, as either version of Pasch
implies the other, and either continuity axiom implies the other.
Perhaps these implications can also be proved without Markov's principle,
but verifying that would add nothing to the philosophical points of this paper.   
Gupta proved (constructively and without Markov's principle),
using inner Pasch, that the base of an 
isosceles triangle has a midpoint.  Having circle-circle continuity 
allows us to prove Euclid I.1, so every segment is the base of an isosceles
triangle and hence has a midpoint.  From there we can proceed
to perpendiculars and thence to the main body of \cite{beeson2016b},
eliminating the need to check \cite{schwabhauser}.
In particular, we need to check neither Gupta's proof of outer Pasch
nor his construction of midpoints without using circles, both of which
are long and difficult. Why then, do we need outer Pasch? 
 Aside from its innate interest, outer Pasch 
is the key tool to prove the plane separation theorem, Theorem~2.5 of 
\cite{beeson2015b}.%
\footnote{Gupta's proof of outer Pasch from inner Pasch does not use 
the parallel axiom.  That raises the possibility that it might be 
easy to prove outer Pasch if we allow the use of the parallel axiom.
But then, we have to ask which version of the parallel axiom.  
In \cite{beeson2015b}, Theorem 9.2, we show that Euclid~5 implies 
Tarski's parallel axiom, but we used outer Pasch in the proof.  It is 
not difficult to prove that Tarski's parallel axiom plus inner Pasch
implies outer Pasch.  But we do not have a short proof of outer Pasch
from Euclid~5.  We do have a short proof of inner Pasch from outer Pasch
and Euclid~5.}  

\subsection{Midpoints}
Euclid's own midpoint construction (Prop. I.10)
 is to construct an isosceles triangle on $pq$ and 
then bisect the vertex angle.  But Euclid's proof is defective; a correct
proof has to rely on Pasch's axiom.  Gupta showed in 1965 \cite{gupta1965} 
how to construct midpoints without using circles.  His proof is complicated,
but we are allowed to use circles, so Gupta's construction is irrelevant
for this paper. 

One of Gupta's ``simpler'' theorems enables us to justify the second part of this 
Euclidean midpoint construction, and we present the main construction
of that lemma next, for the 
reader's enjoyment (it is short and beautiful), and because we need to 
check (in the next lemma) that it relies only on a version of Pasch that is acceptable when 
Markov's principle is rejected.

\begin{theorem}   [Gupta]
Assuming Markov's principle, 
the base of an isosceles triangle always has a midpoint.   
 \end{theorem}

\begin{figure}[ht]
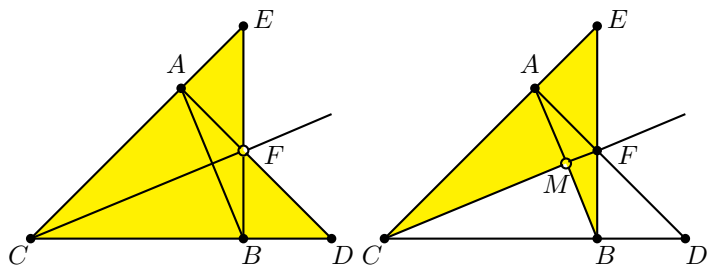

\caption{To construct the midpoint $M$ of $AB$, given $C$ with $AC$ equal to $BC$, constructing first $F$ and then $M$ by inner Pasch.}
\label{figure:GuptaMidpointFigure}
\center{\GuptaMidpointFigureOne \qquad\GuptaMidpointFigureTwo}
\end{figure}
\FloatBarrier

\noindent
{\em Proof}.  See page 56 of \cite{gupta1965}. We repeat the proof here,  
because below we will adapt it to work without Markov's principle.
First we just repeat it as Gupta gave it.  We are not claiming 
that this proof works in EG.  
\medskip

Let $ABC$ be an isosceles triangle with $AC$ equal to $BC$. 
It is desired to find the midpoint $M$ of $AB$. 
Let $\alpha$ and $\beta$ be any two distinct points, and
extend both $CB$ and $CA$ by $\alpha\beta$ to produce points
$D$ and $E$ as shown in Fig.~\ref{figure:GuptaMidpointFigure}.
As shown in the figure, two applications of inner Pasch produce
points $F$ and $M$.  The point $M$ is the desired midpoint.
The proof that $M$ is in fact the desired midpoint is not quite
straightforward, but it appears in \cite{gupta1965}, p. 56,  
so we do not repeat it here.

We now consider whether Gupta's proof still works without Markov's
principle.  The point is that to use non-Markovian Pasch we would
need to prove that the angles in question are positive and less than $\pi$.
In general there seems to be no way to justify that, since we do not 
yet have available midpoints and perpendiculars.  However, 
we do have circle-circle continuity, and hence Euclid I.1, 
so we only need to make Gupta's proof work for equilateral triangles,
not for isosceles triangles in general.   And, in Gupta's proof,
we get to choose the points $\alpha$ and $\beta$, which Gupta needed
only to be any distinct points.  

Here is  Gupta's theorem for equilateral triangles, proved
without  Markov's principle.   

\begin{theorem} \label{theorem:midpoint-helper} [Gupta without Markov]
Intuitionistic Tarski geometry (without any continuity axioms and without Markov's principle)
proves that if
$A$ and $B$ are distinct points, and  $ABC$ is an equilateral
triangle, then $AB$ has a midpoint.  
 \end{theorem}
 
\noindent{\em Proof}.   We choose $\alpha = A, \beta = B$.
Then $A$ is the midpoint of $CE$ and $B$ is the midpoint of $CD$.
 Because $ABC$ is an equilateral triangle, 
 angle $ACB = 60^\circ$ (in the sense
precisely defined in \S\ref{section:NL}).
Therefore by Lemma~\ref{lemma:sixtydegrees},
$ACB > 0$.  
Let $P$ be the 
reflection of $B$ in $C$, so $\B(P,C,B)$ and $PC = CB$.  
Then $PCA = 120^\circ$, so by Lemma~\ref{lemma:sixtydegrees},
$PCA > 0$ and $ACB < \pi$.
That justifies the first application of inner Pasch, so point $F$
in Gupta's proof exists.

To justify the second application of inner Pasch, we need to 
show that $0 < CEB < \pi$. 
That is, $0 < AEB < \pi$. 
 This is where we use $CA=AE$,
which gives us (using points $A,B,C$ and $E$) the conclusion
$AEB = 30^\circ$.  Hence, by Lemma~\ref{lemma:sixtydegrees},
 $AEB > 0$.  Let $Q$ be the reflection
of $A$ in $E$, so $\B(A,E,Q)$ and $EQ=AE$.  
Then $BEQ = 150^\circ$,
so by Lemma~\ref{lemma:sixtydegrees},
 $BEQ > 0$, so $AEB < \pi$.  Then the second application of 
inner Pasch is justified, and the point $M$ exists.
The fact that $M$, once constructed,  is indeed
the desired midpoint can be proved exactly
as in Gupta's thesis.  
That completes the proof.

\begin{theorem} \label{theorem:midpoint}
Every positive segment has a midpoint.
\end{theorem}

\noindent{\em Proof}.
Since circle-circle continuity enables us to construct an equilateral triangle on any 
segment (via Euclid I.1),  and EG has circle-circle continuity,  the corollary 
follows from Gupta's theorem.  We remark that Euclid I.1,  like 
inner Pasch, does not depend on a dimension axiom, since in $R^n$,
so-called circle-circle continuity is really sphere-sphere continuity,
so Euclid I.1 holds without using a dimension axiom.

\begin{lemma} \label{lemma:crossbar} Let $a$, $b$, and $c$ be 
three distinct points with $0 < abc < \pi$, 
and suppose $\B(b,e,c)$.   Then $a \# e$.
\end{lemma}

\begin{figure}[ht]
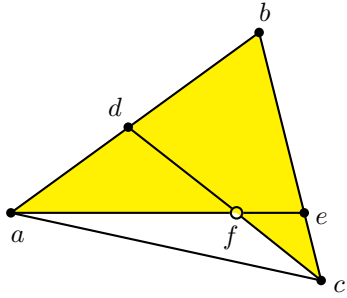

\caption{Given $a \# b$, $a \# c$, $B(b,e,c)$, and $0 < abc < \pi$,
construct $f$ showing $a \# e$.}
\label{figure:lemmacrossbar}
\center{\LemmaCrossbarFigure}
\end{figure}

\noindent{\em Proof}.  See Fig.~\ref{figure:lemmacrossbar}.
 Let $d$ be the midpoint of $ab$.  
Since angle $0 < abc < \pi$,
we can apply inner Pasch to the configuration $adbce$.
The result is a point $f$ with $\B(a,f,e)$.  Hence, by 
the third clause in the definition of $\#$,   $a \# e$.  That 
completes the proof.  
\FloatBarrier

\subsection{Uniform perpendicular} 
In \cite{beeson2015b},  we gave two different constructions of the uniform 
perpendicular (to a line $L$ through a point $p$, without assuming that $p$ is 
or is not on line $L$).   One construction assumes the parallel axiom, but not 
line-circle continuity.  The other assumes line-circle continuity, but not 
the parallel axiom.  (It is an open problem to do it without assuming either of the two.)
In this paper, we assume both those two hypotheses,  so we can disregard the 
comparatively difficult construction that avoids line-circle continuity.  
The other construction is fairly straightforward:  just draw a ``large enough''
circle $C$ about point $p$  (large enough that $C$ meets $L$ in two points $a$ and $b$
with $a \# b$), and then bisect segment $ab$ to find the foot $f$ of the perpendicular,
and then erect the perpendicular to $L$ at $f$,  using the construction of Euclid I.1
to construct an equilateral triangle over $ab$ and connect its vertex to $f$.  
The only tricky part of this construction is to find the radius to use to draw $C$.
Line $L$ is given by two points, say $u$ and $v$, with $u \# v$.  So segment 
$uv$ can be extended;  extend it by $pu$ (which may be null or not, we don't care)
and again by $pv$.  The resulting segment is at least $uv$ longer than $pu$ or $pv$,
so it can serve as the radius of $C$.  

\subsection{Euclid does not need Markov}
 \label{section:MarkovEuclid}
Does Markov's principle, or the stability of congruence, actually play
an important role in Euclidean reasoning?  We searched for a theorem in Euclid Books I-IV
in which Markov's principle is needed.   There is an {\em a priori} constraint, which 
we consider first.  

First, the double negation interpretation.  For simplicity we consider the version of 
EG based on Tarski's language and axioms,  with Skolem symbols $ext$  for the 
segment-extension axiom,  $ip$ for inner Pasch, and $e5$ for Euclid~5.   
As discussed at the end of \S\ref{section:axiomlist}, the axioms are 
then quantifier-free and can be put in the form
 $A \implies B$, where $A$ and $B$
are conjunctions of atomic formulas.   Let $\phi^-$ be the double-negation interpretation 
of $\phi$.  Since double negation commutes with implication and conjunction, each 
axiom $\phi $  satisfies $\phi^- \leftrightarrow \neg\neg \  \phi$.  Each theorem $\psi$ of 
Euclid is also of the form just described.  Hence, if the conclusion of $\psi $ mentions 
only congruence and equality (and not betweenness),  then we will have $\psi^- \implies \psi$,
in EG plus the stability of congruence.   Therefore,  any theorem of Euclid that
actually requires Markov's principle  must mention betweenness in its (formalized) conclusion.  

 We therefore examined
the propositions of Euclid looking for theorems whose conclusions involve
collinearity,  angle ordering, or segment ordering. 
 This investigation is slightly 
complicated, because Euclid  contains well-known errors,
and because Euclid has never yet been formalized faithfully
(i.e., in a theory that permits correcting Euclid's proofs, while 
assuming not too much more than Euclid did).%
\footnote{There have certainly been many correct formal theories of geometry
capable of proving versions of the propositions of Euclid. But that is 
not quite the same thing as formalizing Euclid.  For example, to get
the two circles in Proposition I.1  to intersect, clearly Euclid is 
implicitly assuming some kind of circle-circle continuity.  But Euclid 
has no dimension axiom, so circles are really ``spheres''.  Did Euclid
intend to assume that two circles intersect on a given side of the 
line connecting their centers, or just that they intersect {\em somewhere}?
If only the latter, then his proof that an angle could be bisected is wrong.
If the former, then he would not need to bisect an angle first, but 
could bisect a segment directly by constructing two equilateral triangles
on opposite sides of the segment.  Proposition 7 patently fails in 
three-space, so perhaps Euclid did mean to have a dimension axiom.  But
eventually he works on the Platonic solids, so he needs to {\em not} have
a dimension axiom.  
}
Nevertheless,
those errors can be corrected, so the search for a 
possible use of Markov's principle is possible, even in this rough terrain.

 The first one is Euclid I.14,
which says

\begin{quote}
If with any straight line, and at a point on it, two straight lines not lying on the same
side make the adjacent angles equal to two right angles, the two straight lines will be in 
a straight line with one another.
\end{quote}

\noindent The collinearity statement in the conclusion is naturally expressed using 
betweenness.  For a diagram, refer to your copy of Euclid.  Our point about this example 
is that stability of betweenness is not actually needed.  Let $E$ (in Euclid's figure)
lie on $CB$ extended, with $BE = BC$.  
Let $PB \perp BC$; then $CBP$ and $PBE$ are right angles, so $CBA$ and $ABE$ are equal
to two right angles.  Hence angles $CBA$ and $ABE$ together equal angles $CBA$ and $ABD$.
Hence angle $ABE$ equals angle $ABD$.  It follows that $E$ and $D$ lie on the same 
line $CBE$ on the same side of $B$.   (That follows in Tarski's system from the definition 
of equality of angles.)  We do not need the argument using Markov's principle  that
Euclid gives in I.14.

\subsection{The exterior angle theorem and its consequences}
The next potential example is I.16, the exterior angle theorem.
Euclid's proof constructs a crucial point $F$, but he
left a gap in failing to prove that $F$ lies in the 
interior of a certain angle.

The proof is also instructive, for those not yet
completely familiar with Tarski's treatment of angles, 
in that it shows how to ``unwind'' a theorem about angles, to see what it 
``really states'' 
when angles are eliminated in favor of a points-only formulation.

 \begin{lemma}[Exterior angle theorem, Euclid I.16] 
 \label{lemma:exteriorangle}
 Suppose $\B(B,C,D)$ and $A \# B$ and $0 < BAC < \pi$.
Then  $ACD >BAC$.  
\end{lemma}

\noindent{\em Remark}. The hypotheses are weaker than 
``$ABC$ is a triangle'', which by definition
 requires all three  
angles to be positive and have positive supplements. 
\medskip

\begin{figure}[ht]
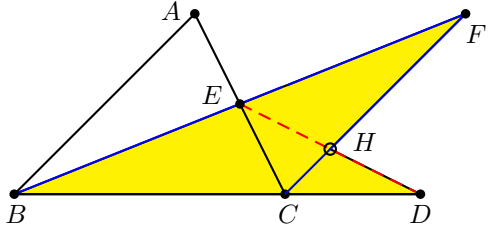

\center{\OneSixteenInnerFigure}
\caption{Triangle $ABC$ has exterior angle $ACD$.
$E$ is the midpoint of $AC$ and also of $BF$.
$H$ exists by inner Pasch.}
\label{figure:I.16}
\end{figure}

\noindent{\em Proof}.  Euclid's diagram is extended by 
 Fig.~\ref{figure:I.16}. 
  We show how to complete Euclid's proof.
To prove $ACD > BAC$, we must construct a point $F$ in the interior of $ACD$
such that $ACF = BAC$.  Euclid knew what angle ordering means:  he
 constructs $F$ as shown in the figure.
We can construct $E$ (and hence $F$) without Markov's principle,
since by Lemma~\ref{lemma:crossbar}, $B \# E$, so segment $BE$ can 
be extended as required.  
The lemma is applicable by the hypothesis $0 < BAC < \pi$.

Euclid neglects to  prove that 
$F$ lies in the interior of $ACD$.  To fill this gap in Euclid,
we need point $H$ in the figure to exist.  
We wish to obtain $H$ from inner Pasch, applied
to the five-point configuration shaded in 
Fig.~\ref{figure:I.16}.  But to use inner Pasch 
in Markov-free geometry, we need an appropriate
angle to be positive.  It would suffice to 
show that $0 < FBC < \pi$, but there seems to be 
no way to do that until the theory of positive 
angles is developed, which cannot precede the 
exterior angle theorem. 
However, our formulation of inner Pasch allows
it to be applied if $0 < ECF < \pi$. 

Since $E$ is the midpoint of both $AC$ 
and $BF$,  and vertical angles are equal (Lemma~\ref{lemma:verticalangles}),
triangles $AEB$ and $CEF$ are congruent, and triangles
$AEF$ and $CEB$ are congruent, by the SAS congruence theorem
(Euclid I.2); that is, $AF = CB$.
  Then triangles $ABC$ and $CFA$ are also congruent;
hence angle $BAC$ is equal to angle $FCA$.
Since $0 < BAC < \pi$, we have  $0 < FCA < \pi$.
But angle $FCA$ is equal to angle $FCE$, since $\B(C,E,A)$.
Angle $FCE$ is equal to angle $ECF$.  Therefore 
$0 < ECF < \pi$, justifying the desired application of inner Pasch.
That completes the proof that $F$ is in the interior of angle $ACD$,
which in turn completes the proof of the exterior angle theorem.
\medskip

\noindent{\em Remark.}  It is an open question whether the exterior 
angle theorem is still provable if inner Pasch is formulated with 
only the hypothesis $0 < acb < \pi$ instead of the hypothesis
$0 < acb < \pi \lor 0 < pqb < \pi$.  We tried various ways of applying
outer Pasch without success.
\medskip

The immediate corollary (Euclid I.17) is that any two angles of a triangle, taken together,
are less than two right angles.  In particular, no triangle contains two right angles.  

\begin{lemma} \label{lemma:legsmallerhypotenuse}
In a right triangle, the hypotenuse is greater than either leg.
\end{lemma}

\noindent{\em Proof}.  First prove Euclid I.18 and I.19.
Then apply them as indicated in Exercise~22, p.~198 of \cite{greenberg}.
Alternately, see Satz~11.46 of 
 \cite{schwabhauser}; the  proof there uses nothing but elementary 
betweenness and congruence, and the existence of perpendiculars.
\medskip

\subsection{Crossbar theorem}
The ``crossbar theorem'' proved in this section
is an easy consequence of outer Pasch,
but I do not know how to derive it from inner Pasch
(without Markov's principle). 

\begin{theorem}\label{theorem:crossbar} Let $a$, $b$, and $c$ be distinct points 
with $0 < abc < \pi$ and $0 < buv < \pi$.
  Suppose $\B(a,e,c)$.  Let $u$ and $v$ be
points with $\B(b,a,u)$ and $\B(b,c,v)$.  Then there exists
a point $w$ with $\B(u,w,v)$ and $\B(b,e,w)$.
\end{theorem}

\begin{figure}[ht]
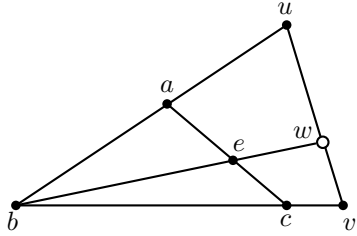

\caption{Crossbar theorem.  $Ray(b,e)$ meets the crossbar $uv$.}
\label{figure:crossbar}
\center{\CrossbarFigure}
\smallskip
\end{figure}

\begin{figure}[ht]
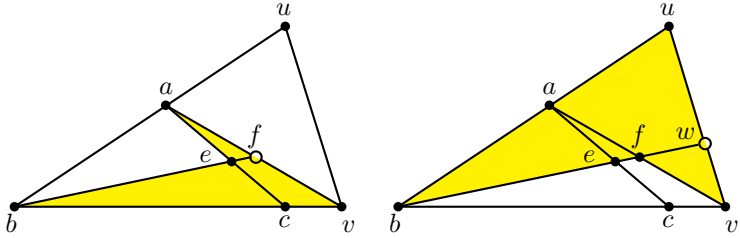

\caption{Proof of the crossbar theorem. Apply outer Pasch twice.}
\label{figure:crossbarproof}
\center{\CrossbarProofFigure \CrossbarProofFigureTwo}
\smallskip
\end{figure}

\noindent{\em Proof}.  Because we are given $0 < abc < \pi$,
we can apply non-Markovian outer Pasch to the configuration $bcvae$.  
The conclusion is the existence of a point $f$
with $\B(a,f,v)$ and $\B(b,e,f)$. 
 Now we apply 
outer Pasch again, this time to the configuration
$bauvf$.  Again, the required hypothesis is fulfilled,
since  $0 < abv < \pi$, and $abv$ is the same 
angle as $abc$.  The result is a point $w$ such 
that $\B(u,w,v)$ and $\B(b,f,w)$.  Now we have
$\B(b,f,w)$ and $\B(b,e,f)$.  By the inner transitivity
of betweenness, we have $\B(b,e,w)$ as desired.
That completes the proof.

\subsection{Further properties of angle ordering}

Euclid took the congruence of all right angles 
 as his Postulate 4.  Hilbert (\cite{hilbert1899}, p. 20) remarks 
that this was ``unjustified'', and says that the proof of it goes back to Proclus.  

 \begin{lemma} \label{lemma:allrightanglescongruent} All right angles are congruent.   In other words, 
if $abc$ and $ABC$ are right angles with $ab = AB$ and $bc = BC$ then $ac = AC$.
\end{lemma}

\noindent{\em Proof}.  This is Satz 10.12 in \cite{schwabhauser}.  However, the proof
appeals only to the definition of angle congruence and simple theorems, such as the 
fact that reflections in points and in lines are isometries.

\begin{lemma} \label{lemma:triangleinequality}
(i) If $a$, $b$, and $c$ are not collinear, the triangle inequality holds: $ac < ab + bc$.
\smallskip

(ii)  Whether or not $a$, $b$, and $c$ are collinear, we have $ac \le ab + bc$.
\end{lemma}

\noindent{\em Proof.}  The proof we gave in \cite{beeson2015b}, Lemma~8.14,
used Markov's principle. But Euclid's proof in I.20, relying on I.5 and I.19,
is perfectly constructive as it stands, not requiring Markov's principle or 
any argument by cases.   Part (ii) has a negative conclusion, since $\le$ can be 
expressed negatively; hence its provability without Markov's principle follows from 
the double negation interpretation.

 \subsection{A lemma about two perpendiculars}
The next lemma 
is used in \cite{beeson2015b} to prove that every Lambert quadrilateral
(plane quadrilateral with three right angles) is a rectangle,  which 
is a key step in establishing a coordinate system on a geometric basis.
In \cite{beeson2015b}, Lemma~8.15,  we used the exterior angle theorem 
and Markov's principle to prove the lemma.   Here we give 
a proof without Markov's principle, relying on inner Pasch and Euclid~5.
The proof in \cite{beeson2015b} does not use Euclid~5, i.e., it is a proof 
in neutral geometry, but it does use Markov's principle.  
  
\begin{lemma}\label{lemma:saccheri-helper}   Let $a,b,c$, and $d$ be   
distinct points (pairwise distinct)
 with $ab \perp ad$ and   $cd \perp ad$.  If $bc$ meets the line 
containing $ad$, 
then the intersection point $m$ is between $a$ and $d$.
\end{lemma}

\begin{figure}[ht]
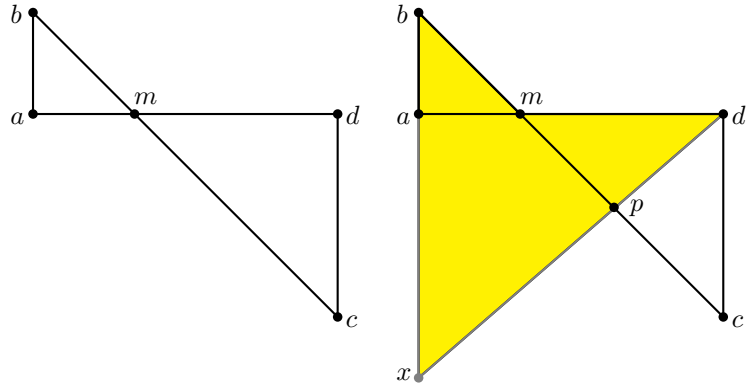

\caption{To show $\B(a,m,d)$, construct $p$,  then construct $x$ by Euclid~5,
then apply inner Pasch to the shaded configuration.}
\label{figure:saccheri-helper}
\SaccheriHelperFigure
\end{figure}

\noindent{\em Proof}.  Please refer to Fig.~\ref{figure:saccheri-helper}.  The 
first part of the figure illustrates the theorem, and the second part shows the 
construction used for the proof. 
 The lines containing $ab$ and $cd$ are parallel,
since both are perpendicular to $ad$.  Let $p$ be any point between $m$ and $c$,
for example, the midpoint of $mc$.  We have $\B(b,m,p)$ and $\B(m,p,c)$.
Hence $\B(b,p,c)$, by the outer transitivity of betweenness (which has been 
discussed above).

Then Euclid~5 applies, and yields the existence of a point $x$ such that
$\B(d,p,x)$ and $\B(b,a,x)$.   Then $bxd > 0$, because it 
is an angle of the right triangle $axd$, which has $a \# d$ by 
hypothesis and $a \# x$ because $\B(b,a,x)$.
Now we can apply inner Pasch to the 
configuration $dpxab$.  That yields a point $j$ such that $\B(b,j,c)$ 
and $\B(a,j,d)$.   Then $j$ is the intersection of $ad$ and $bc$.  Since
those segments are not collinear, their intersection is unique, so $j=m$.
Hence $\B(a,m,d)$ as desired.  That completes the proof.
\FloatBarrier

\subsection{Parallelograms}
We have not yet defined the word ``parallel.''  It may seem 
obvious that two lines are parallel just when they do not intersect,
but that is only true in plane geometry;  skew lines in space 
are not considered parallel, and ``parallel'' can be defined in a 
way that does not depend on the presence of a dimension axiom.
Intuitively, two lines are parallel if they lie in the same plane
and do not meet.  Precisely, we define $ab$ to 
be parallel to $cd$ if $a$ and $b$ are on the same side
of $cd$ and there is no point $x$ such that $L(a,b,x)$
and $L(c,d,x)$.  This definition assumes that $a \# b$ and $c \# d$.
Then it can be proved, without an upper dimension axiom, that 
if $ab$ is parallel to $cd$ then $cd$ is parallel
to $ab$.   The proof requires the plane separation theorem of Hilbert.
A version of that theorem is proved in \cite{beeson2016b}, and it
can also be derived in non-Markovian geometry using the crossbar
theorem (Theorem~\ref{theorem:crossbar}).

The traditional Euclidean results about parallel lines
and parallelograms offer no difficulties in non-Markovian
intuitionistic geometry, provided we define a {\em transversal}
of two lines $L$ and $K$ to be a line $J$ that makes angles
with $L$ and $K$, both of which are between $0$ and $\pi$.

We omit the proofs of the following facts:

\begin{lemma} \label{lemma:parallel}
If lines $L$ and $K$ are cut by a transversal meeting
them at distinct points (and ``transversal'' is defined as above), then 
\smallskip

(i) if $L$ and $K$ are parallel, 
then alternate interior angles are equal, corresponding
angles are equal, and interior angles on the same side 
of the transversal are together equal to two right angles, and
\smallskip 

(ii) if any of the conditions in (i) hold, then $L$ and $K$
are parallel.
\end{lemma}

\begin{definition} A {\em parallelogram} is a quadrilateral
with opposite sides parallel and adjacent vertices distinct,
whose diagonals each are transversals of each pair of opposite sides.
\end{definition}

In other words, the diagonals form angles between 0 and $\pi$ with 
the sides.
It is not required by this definition that diagonally opposite
vertices be distinct.

 The following lemmas can all be proved
without Markov's principle.    

\begin{lemma} \label{lemma:oppositesides}
Opposite sides of a parallelogram are congruent.
\end{lemma}

\begin{lemma} \label{lemma:parallelogram2}
If opposite sides of a quadrilateral are congruent,
and adjacent vertices distinct, and the diagonals meet,
then the quadrilateral is 
a parallelogram.
\end{lemma}

\begin{lemma} \label{lemma:diagonalsofparallelogram}
The diagonals of a parallelogram meet, and
they bisect each other.
\end{lemma} 

\begin{figure}[ht]
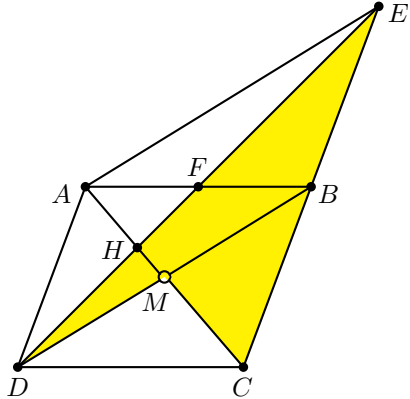

\center{\DiagonalsOfParallelogramFigure}
\caption{The diagonals of parallelogram $ABCD$ 
meet at their common midpoint $M$. $F$ exists by the plane separation theorem. 
$H$ exists by the crossbar theorem.  $M$ exists by 
inner Pasch.}
\end{figure}

\noindent{\em Proof}.  Let $ABCD$ be a parallelogram.
Extend $BC$ to $E$ with $CB=BE$.  Then $E$ is on the opposite
side of $AB$ from $C$.  By definition of parallel, $C$ and $D$
are on the same side of $AB$.  By the plane separation theorem,
$E$ is on the opposite side of $AB$ from $D$.  Therefore
there is a point $F$ collinear with $AB$ with $\B(D,F,E)$.
Then $F$ is in the interior of angle $AEB$.
Note that $0 < DCB < \pi$, since $ABCD$ is a parallelogram.
  By Lemma~\ref{lemma:parallel}
and Lemma~\ref{lemma:oppositesides}, and the SAS congruence theorem,
triangle $ABE$ is congruent to triangle $DCB$.  Hence angle
$AEB$ is equal to angle $DBC$.  Hence $0 < AEB < \pi$.  Then
the crossbar theorem can be applied to angle $AEB$, and since
$\B(A,F,B)$,  the ray $FE$ meets the crossbar $AC$.  Let $H$ be the 
intersection point.  Now we can apply inner Pasch to the 
configuration $DHECB$ (shaded in the figure), because $0 < BDC < \pi$.
The result is a point $M$ such that $\B(C,M,A)$ and $\B(D,M,B)$.
That point $M$ is the desired intersection of the two diagonals of 
$ABCD$.  Now by vertical angles (Lemma~\ref{lemma:verticalangles})
and Lemma~\ref{lemma:parallel}, the four triangles into which the 
diagonals divide $ABCD$ are each congruent to their reflections in $M$.
Hence $M$ is the midpoint of each diagonal.  That completes the proof. 

\begin{corollary} \label{lemma:diagonalvertices}
All four vertices of a parallelogram are distinct.
\end{corollary}

\begin{lemma}\label{theorem:midpointsofoppositesides}
The lines connecting the midpoints of opposite sides
of a parallelogram meet at their common midpoint,
which is also the point of intersection of the diagonals.
\end{lemma} 

\noindent{\em Proof}. This follows from Lemma~\ref{lemma:diagonalsofparallelogram} and the fact that 
reflection in a point preserves betweenness and congruence.

\subsection{Angle copying}
Hilbert took as an axiom the proposition that a given
angle can be copied to a new position, specified by 
a line $L$, a point $p$ on $L$, and a point $q$ not on $L$;
the angle should be copied so that its vertex is at $p$
and its third point $t$ is on the opposite side of $L$ from $q$.
That is, a point $x$ such that $\B(t,x,q)$ is asserted to 
exist.

In our development,  angle copying is a theorem.  We will 
prove that any angle can be copied, not just any positive angle.
Then evidently the conclusion $\B(t,x,q)$ needs to be 
weakened to $\T(t,x,q)$.

To specify a side of a line $L$ given by two distinct
points $p$ and $s$, we have to give a point
$q$ which is distinct from the foot of the perpendicular
to $L$ passing through $q$.  It is not enough just 
to assume $\neg L(p,s,q)$.  We abbreviate this 
property by ``$q$ is distinct from $L$.''

\begin{theorem}[Angle copying]  \label{theorem:anglecopying}
  Let $p$ and $s$ be any two distinct points,
and let $q$ be any point distinct from the line 
containing $p$ and $s$.  Assume $a$ and $c$ are 
distinct from $b$ (but not necessarily $a \# c$).
Then angle $abc$ 
can be copied to the line containing $ps$, with 
vertex and $p$ and the third point lying on the
opposite side of $ps$ from $q$.
 More precisely,
 there exist points $a^\prime$, 
$c^\prime$, and $x$ such that 
angle $abc$ is congruent to angle $a^\prime pc^\prime$, 
$L(p,s,c^\prime)$, and $\T(q,x,a^\prime)$.  Moreover,
if $a \# c$, then $\B(q,x,a^\prime)$.
\end{theorem}

\noindent{\em Proof}.  
Let $c^\prime = ext(s,p,b,c)$. Then $pc^\prime = bc$ 
and $L(p,s,c^\prime)$.  Let $c^\prime$ be a point
with  $bc = pc^\prime$ on $Ray(s,p)$.
Erect a perpendicular line $K$ to $pc^\prime$ at $c^\prime$,
on the opposite side of $ps$ from $q$. 

Then there are two points $a^\prime$ on $K$ such that 
$a^\prime c^\prime = ac$.  Either choice will make
triangle $a^\prime p c^\prime$ congruent to $abc$,
since they are right triangles with corresponding legs
equal.  We need to choose the correct $a^\prime$ to 
be on the opposite side of $pc^\prime$ from $q$.
 By construction of $K$,
one of the two rays of $K$ emanating from $p$ is 
on the opposite side of $ps$ from $q$.  Choose $a^\prime$
on that side of $ps$.
That completes the proof.

\subsection{Positive angles revisited}

Recall that we defined a positive angle by two cases:
a positive angle is either an angle of a right triangle
(with positive sides) or an apex angle. 

\begin{theorem}\label{theorem:positivehypotenuse}
Suppose $abc$ is any angle with $0 < abc < \pi$.
Then $a \# c$.
\end{theorem}

\noindent{\em Remarks}. (1) The special case when $abc$
is a right angle is already interesting.
That special case can be stated 
loosely as, ``every right angle has a positive hypotenuse.''
Note that $abc$ does not count as a right triangle until 
all its sides are positive; to be a right 
{\em angle} it only needs positive legs.  That is why 
we speak of the hypotenuse of an angle here, rather than
the hypotenuse of a triangle.

(2) The proof makes use of the parallel postulate.
But our theory of non-Markovian geometry has already
depended heavily on the parallel postulate.  Our point
in this paper is simply that non-Markovian geometry 
is a coherent conception; we make no attempt to develop
it without the parallel axiom, in contrast to \cite{beeson2015b},
where we did develop constructive neutral geometry to some extent
(with Markov's principle).  
\medskip

\noindent{\em Proof}.  Since $abc$ is an angle, we 
have $a \# b$ and 
$c \# b$.   
We construct a network of 
copies of $abc$ that ``tile'' part of the plane
near that triangle, 
as shown in Fig.~\ref{figure:positivehypotenuse}. 
The idea is that we can draw lines parallel to 
$ab$ and $cb$, because $a \# b$ and $c \# b$, but 
we cannot necessarily draw lines parallel to $ac$ 
because we only know $a\neq c$.  Nevertheless we
can construct the tiling shown in the figure.  Here
are the details:

\begin{figure}[ht]
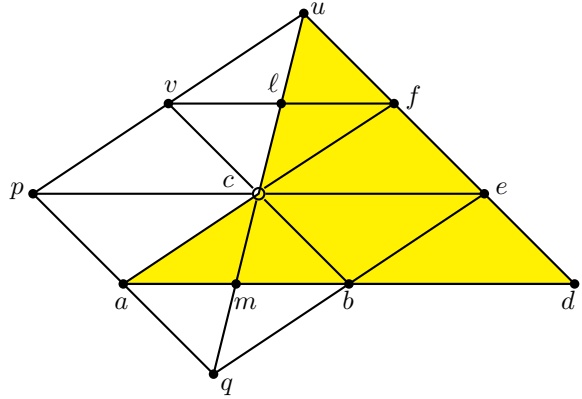

\caption{Starting with $abc$ with $a\#b$ and $b \# c$,
construct the midpoint $m$ of $ab$ and apply Pasch
to prove $\B(a,c,f)$.  Hence $a \# c$.}
\label{figure:positivehypotenuse}
\center{\PositiveHypotenuseFigure}
\end{figure}
\medskip

Extend $ab$ by $ab$ to 
point $d$. 
Copy angle $abc$ so $ab$ goes to $bd$,
and $a$ goes to a point $e$ on the same side of $ab$ 
as $a$.  Then triangle $abc$ is congruent to triangle
$bde$.   Therefore angle $cab$ is congruent to 
angle $ebd$.  Therefore segments $ac$ and $be$ are 
parallel (in the sense that if there are lines through
those segments, they cannot meet, because that would 
contradict the exterior angle theorem).  Angle $cba$ is 
congruent to angle $edb$, so $ed$ is parallel to $cb$.
Since $c \# b$ and $ed = cb$, we have $e \# b$, so we can 
extend $de$ by $de$ to a point $f$ with $ef = bc$. 

We have $cb = ed$ and  $eb = be$. Since $ed$ is parallel
to $bc$, the alternate interior angles $cbe$ and $bed$
are equal.  Hence, by SAS, 
triangle $ceb$ is congruent to triangle $dbe$.
Then $cedb$ is a quadrilateral with opposite sides equal.
Hence, by Lemma~\ref{lemma:parallelogram2}, it is 
a parallelogram.  
Hence, by Lemma~\ref{lemma:oppositesides},
$cd = bd = ab$. By construction of $f$, $ef = ed = bc$.
Hence, by SAS, triangle {\it cef} is congruent to triangle 
$abc$.  Note, however, that we haven't proved 
that $a$, $c$, and $f$ lie on a line, and indeed
the betweenness relation $\B(a,c,f)$ is what we are 
trying to prove, as that is one way to show $a \# c$.
We construct two more copies of triangle $abc$,
by constructing a line {\it vf} through $f$ parallel to $ad$ and 
$ce$, with $vf = ab$ and $v$ on the same side of $f$
as $a$,  and extending $df$ through $f$ by $bc$ to point $u$.
  Let $m$ be the midpoint of 
segment $ab$, which exists since $a \# b$.
Let $\ell$ be the midpoint of $vf$, which exists since
$vf = ab$ and $a \# b$, so $v \# f$. 
 
 Angle $bde$ is congruent to angle $abc$, because 
triangle $cba$ is congruent to triangle $edb$.  Since $0 < abc < \pi$
 by hypothesis, we have $0 < bde < \pi$.  Hence we can
apply inner Pasch to the 
configuration {\it ufdam}.  The result is a point $k$
such that $\B(u,k,m)$ and $\B(a,k,f)$.  In fact
$k = c$, although we haven't proved that.
\ Extending segments 
$eb$ and $uv$, we get a large parallelogram  
$pueq$ as shown in the figure.  The intersection
of its diagonals is at the midpoint of each diagonal,
namely $c$.  Hence $\B(u,c,q)$.  Now $a$ and $f$
are midpoints of opposite sides of the parallelogram,
so by Lemma~\ref{theorem:midpointsofoppositesides},
$\B(a,c,f)$.  Then by definition, $a \# c$.
That completes the proof.

\begin{theorem} \label{theorem:positiveimpliesapex}
Every positive angle is an apex angle.
\end{theorem}

\noindent{\em Proof}. A positive angle, by definition,
is either an apex angle, or a right angle, or an angle of a right triangle.
If it is an apex angle, there is nothing to prove.

A right angle is an apex angle, since given a right 
angle, we can lay off the same distance on each of its
sides, constructing a right isosceles triangle. By
Theorem~\ref{theorem:positivehypotenuse}, the base of 
this triangle has distinct endpoints, and hence, 
it has a midpoint.  That makes the right angle an apex triangle.
  
\begin{figure}[ht]
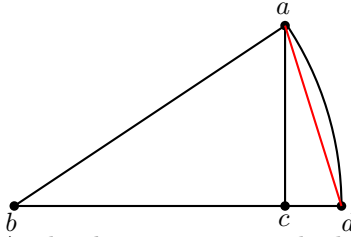

\caption{Angle $abc$ is an apex angle, because 
$a \# d$ as the hypotenuse of right angle $acd$.}
\label{figure:positiveapex}
\center{\PositiveApexFigure}
\end{figure}

It remains to prove that every 
angle of a right angle is an apex angle.  Let $acb$
be a right triangle with the right angle at $c$. 
Fig.~\ref{figure:positiveapex} illustrates the proof.

\medskip

According to Lemma~\ref{lemma:legsmallerhypotenuse},
$ac < ab$.  By definition of triangle, $a \# c$.
Let $e$ be the reflection of $c$ in $b$, so 
$\B(e,b,c) \land eb = ec$.  
($e$ is not shown in the figure.)
Then by Lemma~\ref{lemma:layoff}, there is a unique
point $d$ with $\B(e,b,d)$ and $bd = ab$.  By
Lemma~\ref{lemma:legsmallerhypotenuse}, $bc < ab$.
Hence $bc < bd$.  By definition of $<$, there exists
a point $c^\prime$ with $\B(b,c^\prime,d)$ and $bc^\prime = bc$.
We have $\B(b,c^\prime,d)$ and $\B(e,b,d)$.
Applying inner transitivity (Axiom A15-i), we have $\B(e,b,c^\prime)$.
Now we have
$\B(e,b,c^\prime) \land bc^\prime = bc$, and also
$\B(e,b,c) \land bc = bc$.  By the uniqueness part of 
Lemma~\ref{lemma:layoff}, $c=c^\prime$. Since $\B(b,c^\prime,d)$
we have $\B(b,c,d)$.
 Therefore $c \#d$. 
Now $acd$ is a right triangle since $acb$ is a right angle,
and both legs $cd$ and $ac$ have distinct endpoints. 
Applying Theorem~\ref{theorem:positivehypotenuse} to 
right angle $acd$, we obtain $a \#d$.  Since $ba = bd$,
this proves $abc$ is an apex angle, as desired.
That completes the proof.
\medskip

\begin{theorem}[Angle bisection] 
Every positive angle can be bisected.
\end{theorem}

\noindent{\em Proof}.  Let $abc$ be a positive angle.  Then 
by Theorem~\ref{theorem:positiveimpliesapex},  $abc$ is an apex
angle.  That is, there are points $u,v$ on $Ray(b,a)$ and 
$Ray(b,c)$ respectively such that $bu = bv$ and $u \#v$.
Since $u \# v$,  Theorem~\ref{theorem:midpoint} implies that
$uv$ has a midpoint, say $m$.  By Theorem~\ref{theorem:positivehypotenuse},
since $b \# u$ and $u \# m$, we have
$b\#m$.  Hence $bmu$ and $bmv$ are triangles.  Since $um = vm$,
angle $ubv$ is bisected by $Ray(b,m)$.  Since angle $abc$ is the 
same angle as $ubv$,  angle $abc$ is bisected by $Ray(b,m)$.
That completes the proof.
\medskip

{\em Remark}.  Euclid's Proposition I.9 purports to prove that every
angle can be bisected; but the proof has more than one flaw.  First,
if the angle to be bisected is equilateral, the two points determining 
Euclid's bisector might coincide, and thus fail to determine a line.
Second, without any dimension axiom,  the equilateral triangle produced
by applying I.1  in the proof of I.9 might lie in a different plane 
from the angle to be bisected, thus the ``bisector'' would fail to lie
in the plane of the angle.  These difficulties are not easy to repair,
especially without using a dimension axiom; but it is done as above,
by using Gupta to find the midpoint of the base of an equilateral 
triangle, and then using inner Pasch.  In other words,  the above 
approach, which we use here in non-Markovian geometry, is (as far as I know)
actually the simplest way to repair Euclid I.9.  Of course, if we 
do not care about avoiding Markov's principle, some of the complications
fall away.

\section{Geometric arithmetic without Markov's principle}
In \cite{beeson2016b, beeson2015b}, we showed that in EG plus Markov's principle,
one can define coordinates, and one can define addition and multiplication of points
on a fixed line.  The constructions differ from the classical ones 
because we need to avoid a case distinction according to the signs of the 
arguments.  Two ways of defining addition and multiplication 
constructively are studied in the cited paper.  One is purely geometrical,
in which addition is defined using uniform rotation and reflection, 
and multiplication is defined using a circle, following Hilbert. 
We have checked that Markov's principle is not used in these constructions
and proofs in any essential way.  Actually, we found only one place
where it was used, and we fixed that by providing in this paper a 
non-Markovian proof of Lemma~\ref{lemma:saccheri-helper}.  In the discussion
below, we will point out how that lemma is used. 

\subsection{Coordinates} 
Once we have the uniform perpendicular, it is possible to define coordinates.
Pick two perpendicular lines to serve as the $x$-axis and the $y$-axis.
Then to each point $p$ we can assign   the feet $a$ and $b$ of 
the perpendiculars through $p$ to the $x$-axis and $y$-axis, respectively.  
Then $a$ is the $x$-coordinate of $p$.  To get the $y$-coordinate,
however, we must rotate $b$ by $90^\circ$,  because the $y$-coordinate is supposed
to be a point on the $x$-axis,  the same distance from the origin as $b$.  Therefore
we need to define uniform rotation;  ``uniform'' in the sense that it works without 
a case distinction as to the sign of $b$.  A rotation by $90^\circ$ is equivalent
to two reflections in the line $y=x$; so it suffices to define uniform reflection.

The final step in the construction of coordinates is to show that, given points
$x$ and $y$ on the $x$-axis,  we can construct a point $P(x,y)$ whose coordinates
are $(x,y)$.  Namely, erect a perpendicular $J$ to the $x$-axis at $x$
on the opposite side of the $x$-axis from $(0,-1)$.  (The theorem about 
erecting perpendiculars requires a point not on the line, and then the 
perpendicular is constructed on the opposite side of the line from that point.
This guarantees that $J$ lies in the plane of the two axes. 
 Rotate
$y$ counterclockwise about the origin
by $90^\circ$ (in the plane of the two axes, i.e.,
using the two axes to define the angle of rotation)
 to a point $q$.  
Drop a perpendicular from $q$ to $J$.  The foot of that perpendicular is $P(x,y)$.
But why are the coordinates of $P(x,y)$ equal to $(x,y)$?   Let $p=P(x,y)$ and 
let 0 be the origin.
 By construction 
the quadrilateral $qpx0$ has right angles at $p$,$x$, and 0; it is therefore a 
Lambert quadrilateral (three right angles).  It lies in the plane determined
by the two axes since the points $p$ and $q$ are on the opposite side of 
the $x$-axis from $(0,-1)$.

 Using the parallel postulate, 
we prove every Lambert quadrilateral in a plane 
is a rectangle.  Hence the angle at $q$ 
is also a right angle; hence $q$ is the foot of the perpendicular from $p$ to 
the $y$-axis; hence $y$ is the $y$-coordinate of $p$, as claimed.   

Thus the construction of coordinates depends on the theorem that every planar
Lambert quadrilateral is a rectangle (a  
quadrilateral with four right angles).
We claim that this theorem can be proved,
 without using Markov's principle, from the weakest
form of the parallel postulate (Playfair's axiom).  It is mentioned 
in Lemma~5.24 of 
\cite{beeson2016b},  but the proof is not given there.   In 
\cite{beeson2015b}, Lemma 10.2,  it is proved in intuitionistic
Tarski geometry that a planar Saccheri quadrilateral is a rectangle.
(A Saccheri quadrilateral is a quadrilateral $abcd$ with right
angles at $b$ and $c$ and $ab=cd$.)  Half a Saccheri quadrilateral is 
a Lambert quadrilateral, and a Lambert quadrilateral can be doubled to 
make a Saccheri quadrilateral, so the cited lemma implies that a planar
Lambert triangle is a rectangle, as desired.  The proof of that theorem
given in \cite{beeson2015b} is almost free of Markov's principle,
but it does appeal to Lemma~\ref{lemma:saccheri-helper} (in the numbering
of the present paper), which in the cited paper
is proved using Markov's principle.  However, in this paper we 
  proved Lemma~\ref{lemma:saccheri-helper} without using Markov's principle.
Hence, the uniqueness of coordinates is proved without using Markov's 
principle. 

\subsection{Euclidean fields without Markov's principle}\label{section:EuclideanFields}
In \cite{beeson2016b}, we gave axioms for (intuitionistic) Euclidean fields. 
The language has constants 0 and 1, function symbols $+$, $-$, $\cdot$, $1/x$ or $x^{-1}$ 
for multiplicative inverse, and $\sqrt{}$, 
and a predicate $P(x)$  for ``$x$ is positive.''  The axioms are the field axioms,
the axioms for order (involving $P$), and the axiom that non-negative
 elements have square roots.
In addition we took Markov's principle $\neg\neg \,  P(x) \implies P(x)$ as an axiom.

Now we want to drop Markov's principle.   
When Markov's principle is allowed,  non-zero elements have multiplicative
inverses if and only if positive elements have inverses, since 
$1/x = x/\vert x \vert$, where $\vert x \vert = \sqrt{x^2}$.
This argument works if $\vert x \vert$ is positive when $x$ is nonzero,
but that is equivalent to Markov's principle.  In the absence
of Markov's principle, then, the axiom that positive elements have 
multiplicative inverses does not imply that all nonzero elements have 
multiplicative inverses, and we leave that axiom unmodified.

 Models of the resulting theory will be called
``weakly Euclidean fields'', since we already used ``Euclidean fields'' for models 
that satisfy Markov's principle.   To be precise, the axioms are the usual axioms of commutative ring theory, plus 
the following:
\begin{eqnarray*}
&\neg x \neq y \implies x = y&\mbox{\qquad stability of equality} \\
&0 \neq 1&  \mbox{\qquad EF0} \\
&P(x) \implies \exists y\,( x \cdot y = 1 \land P(y)) &\mbox{\qquad EF1}  \\
&P(x) \land P(y) \implies P(x+y) \land P(x\cdot y)  &\mbox{\qquad EF2} \\
&x+y = 0 \implies \neg(P(x) \land P(y))   &\mbox{\qquad EF3}\\
&x+y=0 \land \neg P(x) \land \neg P(y) \implies x = 0  & \mbox{\qquad EF4} \\
& x + y = 0 \land \neg P(y) \implies \exists\, z( z \cdot z = x) &\mbox{\qquad EF5}
\end{eqnarray*}

Equivalently one may express EF1  and EF5 using the function symbols $1/x$ and $\sqrt x$,
achieving thereby a quantifier-free axiomatization; and if one does not like using a 
function symbol for a not-everywhere-defined function, one may use the Logic of Partial Terms
(LPT), as explained in \cite{beeson2015b}.
Note that EF1 requires that the inverse positive element not only exist
but also be positive itself.   If Markov's principle is allowed, one 
can prove that the inverse of a positive element is positive, but since
we are not assuming Markov's principle, that cannot (as it turns out)
be proved, and we assume it as part of EF1.

Let $\F$ be an ordered field.  Then ``the plane over $\F$'',  denoted by $\F^2$, is a geometrical structure
determined by defining relations of betweenness and equidistance in $\F^2$, using the given order of $\F$.
Namely, $b$ is between $a$ and $c$ if it lies on the interior of the segment $ac$.  That
relation can be expressed more formally in various ways, for example (using the cross product) 
by  $(c-b) \times (b-a) = 0$  and $(b-c) \cdot(a-b) > 0$.   

 The equidistance relation $E(a,b,c,d)$,
which means that segment $ab$ is congruent to segment $cd$,  can be defined by $(b-a)^2 = (d-c)^2$.
Note that no square roots were used, so these definitions are valid in any ordered field. 

The cross product $u \times v$ is defined using components as
usual. 
If the lines through $uv$ and $st$ intersect (and do not coincide), then
using Cramer's rule, 
 the formula for the 
 intersection point has denominator given by the cross product $D=(t-s) \times (v-u)$.  (Details are in \S9.3 of \cite{beeson2016b}.) 
 So, in verifying inner Pasch and outer Pasch, we will need to 
 know that the cross products for the lines that need to intersect
 are positive.  Thus what has to be proved is the 
 following lemma:
 
 \begin{lemma} \label{lemma:crossproduct}
 In $\F^2$,  the cross product (in the sense defined above) of 
 the vectors giving the sides of angle $abc$ with $0 < abc < \pi$ is
  positive.  In particular, if $d$ lies on $ab$, then $\vert d- c \vert > 0$,
  and if $ab$ lies on the $x$-axis, then $\vert c_2 \vert > 0$. 
 \end{lemma}
 
 \noindent{\em Remark}.  At the end, the proof requires the 
 modified EF1,  which says that inverses of positive elements are positive.
 \medskip
 
 \noindent{\em Proof}.    By a rotation and translation, we can assume that 
 $b = (0,0)$ and $a$ lies on the positive $x$-axis.  (All three cases
 in the definition of ``positive angle'' are invariant under rotation and translation.) Then the cross product in question is $a \times c =c_2a_1$.
 Since $a_1 > 0$,  the cross product has positive absolute value
 if and only if $\vert c_2 \vert > 0$.  The penultimate statement follows from $c_2 > 0$, since $\vert d-c \vert$ is the hypotenuse a right triangle with 
 a vertical leg of length $c_2$, and in $\F^2$ the theorem that 
 a leg of a right triangle is shorter than the hypotenuse holds.

  Assume $0 < abc < \pi$.    If angle $abc$ is a right 
 angle or an angle of a right triangle,  $\vert v_2 \vert > 0$ is immediate.  
 So we may assume 
$abc$ is an apex angle.  Then there are distinct points $u = (u_1,0)$
 and $v$ on the ray from $(0,0)$ through $c$, 
 such that $u_1 > 0$ and $\vert u \vert = \vert v \vert$.  We have 
 $\vert u- v \vert > 0$ since $u$ and $v$ are distinct.
 The cross product of the two sides of the angle is $u \times v = 
 u_1v_2$, since $u_2 = 0$.  We must show $\vert u_1v_2\vert > 0$.
 Since $u_1 > 0$, that is equivalent to $v_2 > 0$.  We have not 
 yet used the assumption $abc < \pi$,  which we will now need.
 Let $d = (-a_1,0)$.  Then $0 < dbc$, by definition of $abc < \pi$.
 There are three cases:  $dbc$ is an apex angle, a right angle, or 
 an angle of a right triangle.  If either of the latter two then $c_2 > 0$,
 so $v_2 > 0$ and we are done.  So we may assume $dbc$ is an apex angle.
 Then there are distinct points $p$ and $q$ on the negative $x$-axis and   
 $bc$, respectively, such that $\vert p \vert = \vert q \vert$ and 
 $\vert p \vert > 0$.  Extending the positive-length vectors $u$ and
 $v$ by $\vert p \vert$ and the vectors $p$ and $q$ by $\vert u \vert$,
 we can assume without loss of generality that $\vert u \vert = \vert p \vert
 = \vert v \vert = \vert q$.  Then $v=q$ and angle $pqu=pvu$ is inscribed
 in the semicircle $pvu$ with center at $b = (0,0)$.   Since classically
 $\F^2$ satisfies the theorem that an angle inscribed in a semicircle
 is a right angle,  and being a right angle is equationally described,
 $pvu$ is a right angle.  $\vert p_2 \vert$ is the altitude of that right triangle.
 By Lemma~\ref{lemma:crossproduct}, $\F^2$ satisfies 
 that the altitude of a right triangle 
 (with positive legs) is positive.  Hence $\vert p_2 \vert >0$.
 But $c_2/\vert c \vert = p_2/\vert p \vert$, and $1/\vert c \vert $
 and $1/\vert p \vert$ are both positive
 by EF1, since $\vert c \vert > 0$.  Hence $c_2$ is a positive
 quantity times $p_2$.  Hence $\vert c_2 \vert > 0$, as desired.
 That completes the proof of the lemma.

The following theorem verifies that we have the 
axioms for weakly Euclidean fields right.

\begin{theorem} \label{theorem:geometrytofields}
Every plane $\F^2$ over a weakly Euclidean field $\F$ is a model of Euclidean
geometry EG.
\end{theorem}

{\em Remark.} This theorem depends on the 
modifications we made to the extension axiom and to inner and outer Pasch; had we not changed the formulations of 
inner and outer Pasch to require an angle to be positive and have a positive 
supplement, we would need Markov's principle to verify them.  Although 
these calculations are elementary,  we need to verify that Markov's
principle is not necessary, and that 
we only need to divide by elements that are known to be positive, 
and we  need 
square roots only of positive elements (as opposed to nonzero elements).
The fact that these calculations succeed without Markov is a strong
indication that we have the axioms of EG right.
\medskip

\noindent{\em Proof}. 
 We define $x > y$ to abbreviate $P(y-x)$ and 
$x \ge y$ to abbreviate $\neg x < y$, and similarly for $<$ and $\le$.
We prove the transitivity of $<$.  Suppose $x < y$ and $y < z$.
Then $P(y-x) \land P(z-y)$.  By EF2,  $P((y-x)+(z-y))$.  
Hence $P(z-x)$.  Hence $x < z$.  

We define absolute value in $\F^2$:
$$\vert (x,y) \vert := \sqrt{x^2 + y^2}.$$
We have $x^2 \ge 0$ for all $x$ (for example, since $\le$ is stable,
we can prove that by contradiction).  Hence $\vert u \vert$ is always
defined.   

The interpretation of betweenness in $\F^2$ has to be given without using
cases.  Collinearity $L(u,v,w)$ can be defined by 
$$  L(u,v,w) :=  (w-u) \times (w-v) = 0. $$
Then betweenness is defined by 
$$ \B(u,v,w) := L(u,v,w) \land \vert v-u \vert + \vert w-v \vert = \vert w-u \vert \land \vert v-u \vert > 0 \land \vert w-v \vert > 0. $$
Since translations and rotations are given by linear maps with
determinant one, they preserve cross and dot products, and hence 
they preserve distance and collinearity, and hence they preserve
betweenness as well.

We claim that two points $u$ and $v$
 in $\F^2$ are distinct if and only if $\vert v-u \vert > 0$. 
According to the definition of $u \# v$,  $u$ and $v$ are distinct
if there is a point $x$ bearing any of the three possible betweenness
relations with $u$ and $v$.  For example, suppose $\B(u,x,v)$.
Then by the definition of betweenness, $\vert x-u \vert > 0$ and 
$\vert v-x \vert > 0$ and $L(u,x,v)$.  We must show $\vert v-u \vert > 0$.
By a rotation and translation, we may assume $u=(0,0)$ and $x = (x_1,0)$
and $y = (y_1,0)$.  Then we have $x_1 > 0$ and $y_1 > x_1$, and we
must show $y_1 > 0$; but that is the transitivity of $>$, which we 
proved above.   Similarly
for the other two possible  betweenness relations.  Thus 
$u \# v$ if and only if $\vert v- u \vert > 0$.

We turn to the betweenness axioms. Consider the axiom $\neg \B(a,b,a)$.
In $\F^2$, we  have $\B(a,b,c) \implies L(a,b,c)$, by the 
definition of $\B(a,b,c)$.  Suppose $\B(a,b,a)$.  Then $L(a,b,a)$.
By a rotation and translation, we can assume $a$ and $b$ lie on 
the $x$-axis.  Then $\B(a,b,a)$ becomes $\vert b_2-a_2 \vert 
+ \vert a_2 - b_2 \vert  = \vert a_2-a_2 \vert$.  But the last
expression is 0.  With $x = b_2 - a_2$, we have $2 \vert x \vert = 0$.
Hence $x = 0$; that is, $b_2 = a_2$.  Hence $a=b$.  But that 
contradicts $\B(a,b,a)$.  Hence our assumption $\B(a,b,a)$ is
contradictory. That is, $\neg \B(a,b,a)$, which is what had to be proved.

Consider the symmetry of betweenness, $\B(a,b,c) \implies \B(c,b,a)$.
Suppose $\B(a,b,c)$; after a rotation and translation,
we can assume $a$, $b$, and $c$
lie on the $x$-axis.  For simplicity we drop subscripts, writing just $a$ 
instead of $a_1$, etc.  Then we have $\vert c - b \vert + \vert b-a \vert 
= \vert c-a \vert$ and $\vert c-b\vert > 0$
and $\vert b -a \vert> 0$.   What has to be proved is an equation representing
non-strict betweenness $\T(c,b,a)$ and the two inequalities $\vert a-b \vert > 0$
and $\vert b-c \vert > 0$.
By the stability of equality, we can argue classically for the 
non-strict betweenness, which we take as proved. The two inequalities 
follow from the easy lemma $\vert -x \vert = \vert x \vert$.  That 
completes the verification of the symmetry of betweenness.

Consider the inner transitivity of betweenness, axiom A15-i:
$$\B(a,b,d) \land \B(b,c,d) \implies \B(a,b,c).$$  
Suppose the hypothesis.  Again we can apply a rotation and translation to 
bring all the points to the $x$-axis.    As before
the non-strict inequality $\T(a,b,c)$ is expressed by an equation,
so by the stability of equality, we can prove it classically, which
we assume done.  
For simplicity we drop
the subscripts, writing just $a$ instead of $a_1$, etc.
Then $a \le b \le c \le d$.  We have $a < b < d$ and 
$b < c < d$ by hypothesis.  We need to prove $a < b < c$.
But we have both inequalities already; there is no existential 
assertion to prove.  That completes the verification of the 
inner transitivity axiom.

Next we turn to the inner Pasch axiom. Please refer to 
Fig.~\ref{figure:InnerOuterPaschFigure}.  By a rotation and translation,
and possibly renaming some of the points, 
we may assume $aq$ lies on the $x$-axis and $c_2 \ge 0$ and $q_1 \ge a_1$.
  We have to show that 
$\vert (q-a) \times (b-p) \vert > 0$.  That is, 
$\vert (b_2-p_2)(q_1-a_1) \vert > 0$.  The hypothesis of 
inner Pasch includes the specification that either angle $bca$
or angle $pqb$ is between $0$ and $\pi$.  First assume $0 < bca < \pi$.
By Lemma~\ref{lemma:crossproduct},
applied to angle $bca$,
any point on $bc$ is at a positive distance from $a$.   
Hence $\vert a-q \vert > 0$.  
Hence $q_1-a_1 > 0$.  Therefore 
by EF1, $1/(q_1-a_1) > 0$.  
Hence it suffices to show $\vert b_2-p_2 \vert > 0$.
Because $\B(c,q,b)$ holds by hypothesis, 
we have $\vert c - q \vert > 0$, so $acq$ is an angle, and 
the same angle as $acb$.  Since being positive respects the 
same-angle relation, and even angle congruence,  $0 < acq < \pi$.
Then by Lemma~\ref{lemma:crossproduct}, $c$ lies
at positive distance from $aq$.  That is, $c_2 > 0$.
Therefore $b_2 < 0$, since $\B(b,q,c)$ and $q_2 = 0$.
Because $\B(a,p,c)$ holds by hypothesis and $c_2 > 0$ and $a_2 = 0$, we have 
$p_2 > 0$.  Now we have $b_2 < 0 < p_2$.  Hence $\vert b_2-p_2 \vert > 0$,
as desired.  That completes the verification of inner Pasch,
in the case $0 < bac < \pi$.
We note that it would not work without Markov's principle 
if we had only $\T(c,q,b)$ for a hypothesis 
instead of $\B(c,q,b)$.

It remains to verify inner Pasch under the assumption $0 < qpa < \pi$
instead of $0 < bca < \pi$.  Then Lemma~\ref{lemma:crossproduct} 
tells us that $\vert (q-p)\times(a-p) \vert > 0$. 
We are free to choose the origin of coordinates so that
$p_1= 0$ and $a_2=q_2 = 0$. Then
\begin{eqnarray*}
\vert (q_1-p_1)(a_2-p_2)-(q_2-p_2)(a_1-p_1) \vert &>& 0 \\
 \vert (-q_1p_2 - (-p_2)(a_1) \vert &>& 0 \\ 
  \vert p_2(a_1-q_1) \vert &>& 0.
\end{eqnarray*}
By Lemma~\ref{lemma:crossproduct}, $\vert a_1-q_1) \vert > 0$.
Hence $p_2 > 0$.  Since $\B(a,p,c)$, that 
implies $c_2 > 0$.  That in turn implies $\vert (a-c) \times (q-c) \vert > 0$,
which in turn implies $\vert (a-c) \times (q-\ell) \vert > 0$, for any
point $\ell$ on the ray $cq$; in particular for the point $\ell$ such 
that $\vert \ell-c \vert = \vert a-c \vert$.  
Hence $0 < acb < \pi$.  But we have already verified that 
under that hypothesis, inner Pasch holds, so that completes
the verification of inner Pasch. 

Now we turn to outer Pasch.  We may assume $a=(0,0)$
and $q = (q_1,0)$, so $aq$ lies on the $x$-axis.
 As with inner Pasch, the crucial
issue is to verify that the lines that need to intersect have 
positive cross product.  For outer Pasch that is $\vert (q-a)\times(b-p)\vert$,
which is $\vert q_1-a_1\vert\, \vert b_2-p_2\vert$.  Since
$0 < baq < \pi \lor 0 < abq < \pi$,   
$q$ and $a$ are distinct points, so $q_1 > a_1$,
and $b_2 > 0$.  
Hence $\vert q_1-a_1 \vert > 0$, and to finish the verification
it suffices to prove $\vert b_2-p_2 \vert > 0$.
Since $\B(b,c,q)$, we have $0 < c_2 < b_2$, and since $\B(a,p,c)$
we have $0 < p_2 < c_2$.  Hence $0 < p_2 < b_2$. That
completes the verification of outer Pasch.

Any axiom whose analytic-geometry interpretation does not involve
a positive occurrence $P(x)$
will be automatically verifiable, by G\"odel's double-negation 
interpretation.  That is, starting with the classical verification,
we just push double negations inwards, and when they reach an 
equality, we drop them by the stability of equality, and the result is 
an intuitionistic verification.  If the statement had 
a positive occurrence of $P(x)$
we would need Markov's principle to drop the double negation, but 
a negative occurrence causes no problem.
If the statement does not involve $P(x)$, the need does not arise.
On these grounds, the 5-segment axiom, line-circle continuity
and circle-circle continuity
(which have occurrences of $P(x)$, counting as 
negative by virtue of being in the hypotheses)
hold in $\F^2$.  (See \cite{beeson2016b} for a more detailed discussion.)

The dimension axioms (once properly formulated) are easily 
verified in $\F^2$.  That completes the proof of the theorem.

\section{Geometric arithmetic}
We have already discussed the construction of coordinates relative
to a pair of perpendicular lines chosen as axes.  The definitions
of multiplication and addition for positive arguments then offers
no difficulty.  For example, the multiplicative inverse of a 
point $(a,0)$ on the $x$-axis is found by erecting a perpendicular to 
the $x$-axis at the point 1, and using Euclid 5 to prove that it 
meets the line containing the line through $(0,0)$ and $(a,1)$.
That point has coordinates $(1,1/a)$. To use the non-Markovian
version of Euclid 5, we need the hypothesis that point $(a,0)$ is 
distinct from $(0,0)$.  The passage from signed to unsigned multiplication
is discussed in \cite{beeson2016b}, where two solutions are given.
For example, the solution by using Hilbert's definition of multiplication
using a circle works in EG (but it does require the non-strict version
of line-circle continuity, which we included partly for that reason). 
It follows that, loosely speaking, ``every model of EG is 
a plane over a Euclidean field.''    
 
\section{Some results about EG}
In this section we present some metamathematical results about 
non-Markovian geometry EG.
 
\subsection{A Kripke model of non-Markovian geometry}
We want to demonstrate that non-Markovian geometry is ``about something''
by giving a model of it in which Markov's principle fails.
  Of course, we cannot give a model in the classical 
sense, since the double negation of Markov's principle is a theorem.
But we can give a two-node Kripke model.
(For an 
introduction to Kripke models and a proof of the completeness theorem, see 
\cite{troelstra}, Part V, pp.~324{\em ff}. Since our model is 
particularly simple, you may be able to understand it without 
a background in Kripke models.)

In view of the fact that non-Markovian geometry is interpretable in 
the theory of Euclidean fields (without Markov's principle),  it 
suffices to give a Kripke model that does not satisfy Markov's 
principle of Euclidean field theory 
in the form $\neg\neg \,  P(x) \implies P(x)$.  Then a corresponding
model of EG can be constructed by replacing each field $\F_\alpha$ 
in the 
first Kripke model by the plane $\F_\alpha^2$ over $\F_\alpha$,
with betweenness and congruence interpreted as in the previous section.

\begin{theorem} Euclidean field theory does not imply Markov's 
principle, and geometry EG does not imply Markov's principle.
\end{theorem}

\noindent{\em Proof}.  The second claim follows from the first 
and Theorem~\ref{theorem:geometrytofields}.  We will construct a Kripke
model of Euclidean field theory in which MP fails.
\def\M{{\mathcal M}}
Let $\F$ be any non-Archimedean field (so it contains some ``infinite''
element, larger than any positive integer).   
Let $\F_0$ be the set of ``finitely bounded'' elements of $\F$, i.e.,
those elements bounded by some positive integer. 
$\F_0$ contains ``infinitesimal'' elements
$x$, by definition those whose reciprocals are not finitely bounded.
 We interpret
$P(x)$ in $\F_0$ to mean that $x$ is positive in $\F$ and not infinitesimal. 
Then $\F_0$ is not a classical model of Euclidean ordered field theory,
but it is still admissible as a node in a Kripke model.

 The definition of a Kripke model includes these two clauses: 
$\M_\alpha \satisfies \neg \phi$ means 
that $\phi$ fails everywhere above $\alpha$, 
  and $\M_\alpha \satisfies A \implies B$ means that if 
$\M_\beta \satisfies A$  for $\beta > \alpha$,
then   $\M_\beta \satisfies B$.
The soundness theorem (easy to check) 
for Kripke models says that the satisfaction relation preserves 
intuitionistic logical consequence.  Hence, to prove that 
Markov's principle is not a theorem of Euclidean field theory,
it suffices to exhibit a Kripke model in which Markov's principle
is not satisfied.   That Kripke model, call it $\M$, will be a two-node 
model. The root node $\M_0$ has universe $\F_0$, but $P(x)$ is 
interpreted in $\M_0$ to mean that $x$ is positive in $\F$ and
not infinitesimal.  The other node, lying ``above'' $\M_0$, is 
$\M_1$, which is the classical model $\F$.
 
We now show that Markov's principle fails in $\M$.
Let $x$ be infinitesimal in $\F$.
Then the root node $\M_0$  satisfies $\neg\neg \,  P(x)$, since 
$x$ is positive in $\F$, and satisfaction in $\M_1$ is just classical
satisfaction. 
 Suppose,
for proof by contradiction, that $\M$ satisfies Markov's principle.
Then the root node satisfies the hypothesis $\neg\neg \,  P(x)$ 
of Markov's principle, so it must satisfy the conclusion $P(x)$.
But, by definition of $P$ in $\M_0$,  that implies 
$x$ is not infinitesimal, contradiction.  Hence Markov's
principle does not hold in $\M$.

We still need to check that the axioms of Euclidean field theory 
hold in $\M$.  Since satisfaction at node $\M_1$ is just classical
satisfaction, we only need to check $\M_0$. 

Consider EF0: $\neg x \neq y \implies x = y$.  Suppose $\M_0$
satisfies $\neg x \neq y$.  Then $\M_1$ does not satisfy $x \neq y$.
Hence $x=y$ in $\F$.  Hence $x=y$ in $\F_0$.  Hence $\M$ satisfies EF0.

Consider EF1: $P(x) \implies \exists y\,( x \cdot y = 1)$.
Suppose $\M_0 \satisfies P(x)$. Then $x$ is positive and not
infinitesimal.  Hence $y =1/x$ in $\F$
 is finitely bounded, so it belongs
to $\F_0$.  Hence $x$ has an inverse in $\F_0$.  EF1 also requires
that inverse to be positive, which it is, since it is positive 
in $\F$. 
Therefore $\M_0$ satisfies EF1.

Consider EF2: $P(x) \land P(y) \implies P(x+y) \land P(x\cdot y)$.
Suppose $\M_0 \satisfies P(x) \land P(y)$.  Then $x$ and $y$ are 
positive and not infinitesimal.  Then their sum and product are
positive.  Since $x$ is not infinitesimal, by definition $1/x$
is finitely bounded.  Let $n$ be an integer such that $1/x < n$. Then
$$ \frac 1 {x+y} < \frac 1 x < n.$$   Hence $1/(x+y)$ is finitely
bounded.  Hence $x+y$ is positive in $\M_0$.   Also 
$$ \frac 1 {xy} = \frac 1 x  \frac 1 y <  nm$$
provided $1/y < m$.  Hence $xy$ is also positive in $\M_0$. 
Hence $\M_0$ satisfies EF2.

Consider EF3:  $x+y = 0 \implies \neg(P(x) \land P(y))$. 
Suppose $x+y = 0$ holds in $\M_0$ and $P(x)\land P(y)$ holds in $\M_0$ or 
in $\M_1$.  Then $x$ and $y$ are positive in $\F$, contradicting 
the fact that $\F$ satisfies EF3.  Hence $\M$ satisfies EF3.

Consider EF4: $x+y=0 \land \neg P(x) \land \neg P(y) \implies x = 0 $.
Suppose $\M_0 \satisfies \neg P(x)$. Then $\M_1$ does not satisfy $P(x)$,
so $x$ is not positive in $\F$.  If the same is true of $y = -x$,
then since $\F$ classically satisfies EF4,  $x=0$.  Hence $\M_0$
satisfies EF4. 

Consider EF5, namely $x + y = 0 \land \neg P(y) \implies \exists\, z( z \cdot z = x)$.  Suppose $\M_0$ satisfies $x+y=0 \land \neg P(y)$.  
Then $y = -x$ is not positive in $\F$, so by EF5 in $\F$, $x$ has a 
square root $z$ in $\F$.  We claim $z$ is finitely bounded, and hence 
lies in $\F_0$.  If $\vert z \vert < 1$ we are done.  Otherwise, let $n$ be an integer greater that $\vert x \vert $; then $\vert z \vert \le \sqrt n < n+1$.
That completes the proof of the theorem.
\medskip

{\em Remark}.  In the corresponding geometrical Kripke model,
whose two nodes are $\M_0 = F_0^2$ and $\M_1 = F_1^2$, with 
betweenness and congruence defined in the obvious way (details in \cite{beeson2016b}),  two points of $\M_0$ are distinct if and only 
if the difference of both their coordinates is infinitesimal,  and 
an angle is positive if and only if the rays that form its sides
contain points at equal distances from the vertex that are 
distinct; loosely speaking, the angle is not infinitesimal.  
This model is in accord with the intuition that if Markov's
principle does not hold, one has to imagine that every point has 
a ``cloud'' of nearby points, infinitesimally close to it but not 
equal;  so near that the two points cannot determine a line.
\medskip

{\em Remark}. This model, and another more complicated model,
were given in \cite{beeson2016b}, to show that Markov's principle
does not follow from EF0-EF5.  But it was not appreciated at the
time that Markov's principle {\em does} follow from the axioms 
of constructive geometry,  unless the extension axiom and Pasch's axiom(s)
are modified by requiring distinct points and positive angles;
geometry without Markov's principle was purposely not dealt with 
in \cite{beeson2016b}, since ``sufficient to the day are the 
tribulations thereof.''  In other words, the difficulties of 
non-Markovian geometry then seemed intimidating, and were postponed. 

\subsection{EG is not too weak}
We created the theory EG by starting with constructive geometry
from \cite{beeson2016b},  removing Markov's principle, and modifying
the segment extension and Pasch axioms.  We now show that
 we have not made EG too weak, in the sense that if 
we add Markov's principle back to EG, we also automatically get
back the unmodified axioms of constructive geometry. 

\begin{lemma} \label{lemma:nottooweak}
EG plus Markov's principle proves that any two unequal points
are distinct, and if $a,b,c$ are not collinear, then 
$0 < abc < \pi$.
\end{lemma}

\noindent{\em Proof}.
Let $a \neq b$.
Let $\alpha\# \beta$ (two such points exist by 
the lower dimension axiom of EG).  Extend $\alpha\beta$ by $ab$ to 
point $e$. Then $T(\alpha,\beta,e)$.  Since $a \neq b$ and $ab = \beta e$,
we have $\beta \neq e$. 
With Markov's 
principle,  we have $\B(x,y,z)$ equivalent to
 $\T(x,y,z) \land x \neq y \land y \neq z$.  
 Hence, with Markov's principle, $\B(\alpha,\beta,e)$.
Hence $\beta \# e$.  Then $ab$ is congruent to the positive segment
$\beta e$.  Hence it can be extended by $cd$ to a new point $f$, 
using the segment
extension axiom of EG.  Then $\B(a,b,f)$.  Then by definition
$a \# b$.

Now suppose $ \neg L(a,b,c)$.  We will prove
$0 < abc < \pi$.
Since $\neg L(a,b,c)$, we have $a \neq b$ and $b \neq c$, and
using Lemma~\ref{lemma:layoff}, we 
can assume without loss of generality that
 $ab = bc$.  Since $\neg L(a,b,c)$ we have $a \neq c$.
As shown above, Markov's principle then implies 
$a \# c$, so $0 < abc$.  Let $d$ be the reflection of $c$ in $b$;
then similarly, $0 < abd$, so $abc < \pi$. 
That completes the proof of the lemma.

\begin{theorem} \label{theorem:nottooweak}
EG plus Markov's principle proves the axioms of constructive 
Tarski geometry given in \cite{beeson2016b}.   
\end{theorem}

\noindent{\em Proof}.  We will prove the segment extension axiom A4, 
which says that if $a \neq b$ then $ab$ can be extended by any
segment $cd$.  By Lemma~ref{lemma:nottooweak}, with Markov's
principle
$a \neq b$ implies $a \#b$, so the segment extension 
axiom of EG applies.  That completes the proof of axiom A4.  

Next, we will verify the (unmodified) inner Pasch axiom.  In 
constructive geometry, the hypotheses include $\neg L(a,b,c)$.
By Lemma~\ref{lemma:nottooweak}, $0 < abc < \pi$.  Then
the inner Pasch axiom of EG can be applied. 
 That completes the 
verification of inner Pasch.   

The dimension axiom of EG immediately implies the dimension axiom 
of constructive geometry, and no other axioms were modified.
That completes the proof of the theorem.

\subsection{EG is strong enough}
In this section we advance the following principle
\begin{quote}
   If, in any theorem of Euclidean geometry of the type 
   found in Books I-IV,  we replace ``$a$ and 
   $b$ are unequal'' by ``$a$ and $b$ are distinct'',  and we 
   require that all angles be positive,  then the resulting theorem
   is provable in non-Markovian geometry EG.
\end{quote}
The first difficulty here is ``of the type found in Books I-IV.''  
The principle surely fails if the theorem in question is Markov's
principle itself.  But Markov's principle is not a proposition 
of the kind found in Euclid.  As discussed above, those propositions
are ``Horn formulas''; that is, their hypotheses are atomic 
statements (possibly involving defined concepts) and some given points, and the conclusions
assert the existence of some points bearing positive relations to 
the original points. 

 A second difficulty is that the language
of EG does not mention ``angles'' directly, so what exactly is meant
by ``require that all angles be positive''?  The phrase ``angle $abc$''
means (in ordinary geometry) that the three points are distinct and 
not collinear.   We call a formula $\phi$ ``of Euclidean form'' if it 
is expressed as a Horn formula in a language extending EG by a 
3-ary predicate $Angle$ and a 2-ary predicate $D$, where $D$ 
occurs only negatively (that is, in the antecedent of the implication).
  Then the 
``ordinary interpretation'' of $\phi$ is obtained by replacing 
$D(a,b)$ by $a \neq b$ and $Angle(a,b,c)$ by $\neg L(a,b,c)$.
The ``non-Markovian interpretation'' is obtained by replacing 
$D(a,b)$ by $a \# b$ and $Angle(a,b,c)$ by $0 < abc < \pi$,
both of which are abbreviations for existential formulas of EG.
The existential quantifiers in the antecedent are then pulled out, so 
that the non-Markovian interpretation of $\phi$ has additional free
(or universally quantified) variables.   Now we are
able to state a theorem:

\begin{theorem} Let $\phi$ be a formula of Euclidean form 
whose ordinary interpretation is provable in constructive geometry.
Then the non-Markovian interpretation of $\phi$ is provable in EG.
\end{theorem}

\noindent{\em Remark}. By Theorem~\ref{theorem:nottooweak}, 
``provable in constructive geometry'' is equivalent to ``provable 
in EG + MP.''
\medskip

\noindent{\em Proof.}    The non-Markovian
interpretation of any Euclidean formula is a Horn formula, so 
it suffices to prove that Markov's principle is conservative
over EG for Horn theorems.  Let $\phi$ be a Horn formula
provable from a list of axioms $\Gamma$ and formulas $\neg\neg P \implies P$,
where $P$ is atomic (Markov's principle is of that form).  
We write one instance of Markov's principle explicitly and 
suppose $\Gamma$ may contain more instances.  We will show how 
to eliminate one instance at a time.  By cut-elimination, 
there is a cut-free derivation in Gentzen's sequent calculus G3 
\cite{kleene1952} of 
$$\Gamma, \neg \neg P \implies P \seq \phi.$$
Since $\phi$ and $\Gamma$ do not contain quantifiers or disjunction, by \cite{kleene1951},
the inferences may be permuted so that the implication on the left is the 
last-introduced connective, introduced by rule $\implies \seq$.
Then there are cut-free derivations of $P \seq \phi$ and of
$ \Gamma \seq \neg \neg P$.
But again we can permute the inferences so that the last inference
introduces $\neg \neg P$.  Hence there is a cut-free derivation of 
$\Gamma, \neg P \seq \perp$.  Permuting again if necessary, there is 
a cut-free derivation of $\Gamma \seq P$.    Recall that we 
have a cut-free derivation of $P \seq \phi$.  Applying the cut rule 
to these two derivations, we obtain a derivation (with one cut) 
of $\Gamma \seq \phi$.  Hence we have eliminated one instance of 
$\neq \neg P \implies P$.  Applying cut elimination, we 
now have a cut-free proof of $\Gamma \seq \phi$, so we have 
eliminated one instance of Markov's principle.  Continuing in 
this fashion, or more formally,  by induction on the number of 
instances of Markov's principle in the original proof, we eventually
eliminate all instances of Markov's principle.  That completes
the proof of the theorem.
\medskip

In \cite{beeson2016b} we concluded, by a combination of geometry 
and metamathematics, that the essential content of Euclid Books I-III
is formalizable in constructive geometry.  Now, the above metatheorem
allows us to conclude that
 Euclid Books I-III are provable without Markov's
principle,  if the hypotheses of unequal points and nonzero
angles are strengthened as above to distinct points and positive angles.

\section{Conclusions}
Much of Brouwer's work was focussed on the nature of the continuum.
Although most, but not all, of his work preceded the development of recursive analysis and recursive realizability 
(which he never mentioned, even in later papers), he was clearly 
aware that ``lawlike'' or computable real numbers are not enough to 
``fill out'' the geometric continuum;  therefore choice sequences 
or even lawless sequences were needed.  Brouwer
pointed out the non-constructive nature of the ``intersection theorem'',
and by implication, the classical view of elementary geometry.

Our point in this paper is that Brouwer's criticisms of 
``elementary geometry'' definitely do not apply to geometry
as Euclid wrote it.   We have shown that
there is a coherent and beautiful  
theory of intuitionistic 
Euclidean geometry that neither assumes nor implies 
Markov's principle.   

The consistency and coherence of this theory of non-Markovian
 Euclidean geometry
shows that Euclidean geometry, the law of the excluded
middle,  and Markov's principle are three separate issues.
One can choose to accept Markov's principle (but 
not the law of the excluded middle),  in which case one
gets a very nice constructive 
geometry \cite{beeson2015b,beeson2016b}; 
or one can choose not to accept Markov's principle, in which case one 
is forced to make some finer distinctions, but still gets a coherent
intuitionistic geometry, as shown here.

\end{document}